\documentclass[a4paper,fleqn]{cas-sc}

\usepackage[square,sort,comma,compress,numbers]{natbib}
\usepackage{graphicx}
\usepackage{graphics}
\usepackage[section]{placeins}
\usepackage{caption}
\usepackage{subcaption}
\usepackage{epsfig}
\usepackage{booktabs}
\usepackage{color}
\usepackage{amsmath}
\usepackage{mathtools}
\graphicspath{ {./figures/} }
\usepackage{booktabs}
\usepackage{multirow}
\usepackage{pifont}

\newcommand{\ary}[1]{\boldsymbol{\mathsf{#1}}}
\newcolumntype{L}{>{$}l<{$}}
\newcolumntype{C}{>{$}c<{$}}
\usepackage{lineno}
\def\tsc#1{\csdef{#1}{\textsc{\lowercase{#1}}\xspace}}
\tsc{WGM}
\tsc{QE}
\tsc{EP}
\tsc{PMS}
\tsc{BEC}
\tsc{DE}
\begin{document}
\let\WriteBookmarks\relax
\def\floatpagepagefraction{1}
\def\textpagefraction{.001}
\shorttitle{ Piezoelectric Shell Vibration Analysis Using Catmull-Clark Subdivision Surfaces}
\shortauthors{\today}

\title [mode = title]{Vibration Analysis of Piezoelectric Kirchhoff-Love Shells based on Catmull-Clark Subdivision Surfaces}                      
% \tnotemark[1]

% \tnotetext[1]{This research project is funded by the Engineering and Physical Sciences Research Council (EPSRC), Strategic Support Package: Engineering of Active Materials by Multiscale/Multiphysics Computational Mechanics - EP/R008531/1.}

%\tnotetext[2]{The second title footnote which is a longer text matter
%   to fill through the whole text width and overflow into
%   another line in the footnotes area of the first page.}

\author[1]{Zhaowei Liu}[orcid=0000-0002-0572-7415]
\cormark[1]
%\fnmark[1]
\ead{zhaowei.liu@glasgow.ac.uk}
%\ead[url]{www.cvr.cc, cvr@sayahna.org}

%\credit{Conceptualization of this study, Methodology, Software}

\address[1]{Glasgow Computational Engineering Centre, University of Glasgow, Glasgow, G12 8LT, United Kingdom}

\author[1]{Andrew McBride}[orcid = 0000-0001-7153-3777]

\author[1]{Prashant Saxena}[orcid = 0000-0001-5071-726X]

\author[2]{Luca Heltai}[orcid = 0000-0001-5514-4683]
\address[2]{SISSA (International School for Advanced Studies), Via Bonomea 265, 34136 Trieste, Italy}

\author[3]{Yilin Qu}[]
\address[3]{State Key Laboratory for Strength and Vibration of Mechanical Structures, Xi’an Jiaotong University, Xi’an 710049, Shaanxi, China}

\author[1,4]{Paul Steinmann}[orcid = 0000-0003-1490-947X]

\address[4]{Institute of Applied Mechanics, Friedrich-Alexander Universit\"at Erlangen-N\"urnberg, D-91052, Erlangen, Germany}

\cortext[cor1]{Corresponding author}

\begin{abstract}
%This template helps you to create a properly formatted \LaTeX\ manuscript.
%
%\noindent\texttt{\textbackslash begin{abstract}} \dots 
%\texttt{\textbackslash end{abstract}} and
%\verb+\begin{keyword}+ \verb+...+ \verb+\end{keyword}+ 
%which
%contain the abstract and keywords respectively. 
%Each keyword shall be separated by a \verb+\sep+ command.
An isogeometric Galerkin approach for analysing the free vibrations of piezoelectric shells is presented. The shell kinematics is specialised to infinitesimal deformations and follow the Kirchhoff-Love hypothesis. Both the geometry and physical fields are discretised using Catmull-Clark subdivision bases. It provides the required $C^1$ continuous discretisation for the Kirchhoff-Love theory. The crystalline structure of piezoelectric materials is described using an anisotropic constitutive relation. Hamilton's variational principle is applied to the dynamic analysis to derive the weak form of the governing equations. The coupled eigenvalue problem is formulated by considering the problem of harmonic vibration in the absence of external load. The formulation for the purely elastic case is verified using a spherical thin shell benchmark. Thereafter, the piezoelectric effect and vibration modes of a transverse isotropic curved plate are analysed and evaluated for the Scordelis-Lo roof problem. Finally, the eigenvalue analysis of a CAD model of a piezoelectric speaker shell structure showcases the ability of the proposed method to handle complex geometries.
\end{abstract}
%\begin{graphicalabstract}
%\includegraphics{figs/grabs.pdf}
%\end{graphicalabstract}

%\begin{highlights}
%\item Research highlights item 1
%\item Research highlights item 2
%\item Research highlights item 3
%\end{highlights}

\maketitle

\section{Introduction}
\label{intro}
% \subsection{Piezoelectric effect}
Piezoelectricity is a reversible two-way coupling effect resulting from electromechanical interactions in certain crystalline materials. In 1880, \citet{curie1880piezoelectric} discovered the direct piezoelectric effect whereby a mechanical excitation generates an electrical potential. Shortly thereafter, \citet{lippmann1881principle} derived the converse piezoelectric effect from fundamental thermodynamic principles. In 1881, \citet{curie1881contractions} proved its existence as a strain that occurs when an electric field is applied. Shortly after the piezoelectric phenomenon was discovered, Langevin and Rutherford independently applied the piezoelectric effect for submarine detection devices~\cite{katzir2012knew}. In the subsequent hundred years, the piezoelectric effect has been extensively studied and a wide range of novel piezoelectric materials and devices invented and applied to engineering applications. The direct piezoelectric effect is used in sensors/transducers~\cite{redwood1961transient,jaffe1965piezoelectric,abboudnedelec1995,tzou1990distributed,ng2005sensitivity,safari2008piezoelectric} and energy harvesters~\cite{erturk2011piezoelectric,kim2011review}, while the converse piezoelectric effect is used in resonators~\cite{benes1979piezoelectric, nowotny1987general, hollkamp1994multimodal} and actuators~\cite{hagood1990modelling,dosch1992self,hwang1993finite}.

Piezoelectric sensors and actuators are often constructed from films, plates and shells as they can generate large strains under small loads. Early studies of piezoelectric structures focused on simple geometries such as rods~\cite{ragland1967piezoelectric}, plates~\cite{tiersten1962forced, tiersten1963thickness} and cylindrical shells~\cite{nowinski1963nonlinear,paul1966vibrations}. Laminated piezoelectric plates~\cite{heyliger1994static,dash2009nonlinear} are also well studied. With the development of active, adaptive and smart structures, piezoelectric materials are widely used because of their ability to achieve a precise and complex mechanical response to electrical loads. This motivates the requirement for analysis of piezoelectric structures with complex geometries. The finite element method is the ideal modelling framework to analyse such complex structures and to deal with inherent nonlinearities. \citet{allik1970finite} proposed a three-dimensional finite element method for electroelastic analyses, focussing mainly on piezoelectric vibrations. The early works of the piezoelectric finite element method have been reviewed by~\citet{benjeddou2000advances}. \citet{tzou1990distributed} evaluated the performance of intelligent piezoelectric thin plates using a finite element approach. \citet{hwang1993finite} developed a finite element model of laminated plates with piezoelectric sensors and actuators. A nonlinear finite element approach to phase transition in piezoelectric materials was proposed by~\citet{ghandi1996nonlinear}, while \citet{lam1997finite} analysed piezoelectric composite laminates. A static and dynamic analysis of a piezoelectric bimorph was undertaken by~\citet{wang2004finite}.

% \subsection{Isogeometric shell formulation}
Although many three-dimensional finite element approaches for piezoelectric structures have been proposed, work on piezoelectric Kirchhoff-Love shells is limited. 
Kirchhoff-Love and Reissner-Mindlin shell theories categorise shells into "thin" and "thick" according to the ratio of curvature radius to thickness. The Kirchhoff-Love shell theory, also called "classical shell model", is tailored to thin shells. The Reissner-Mindlin shell theory is an extension of the Kirchhoff-Love theory, which can be applied to both thin and thick shells since it accounts for shear deformations. However, Reissner-Mindlin shells theory require additional rotational degrees of freedom, resulting in a larger system matrix than the Kirchhoff-Love shell theory.
Kirchhoff-Love shells require only three translational degrees of freedom, which is computationally more efficient. However, the Kirchhoff-Love finite element method requires $C^1$ continuity of the basis functions while a conventional Lagrangian interpolation only provides $C^0$ continuity. 

\citet{hughes2005isogeometric} presented the framework for isogeometric analysis (IGA) in 2005. IGA provides higher-order continuity by using splines as interpolation functions and it allows for exact geometric representation which completely eliminates geometry error in the numerical solution. However, volume parameterisation of a Computer Aided Desgin (CAD) model is the most challenging problem for IGA~\cite{cottrell2009isogeometric}. Shell formulations are well suited for IGA since they only require a discretisation of the mid-surfaces of the shell.
\citet{kiendl2009isogeometric} developed an isogeometric approach for Kirchhoff-Love shells using Non-Uniform Rational B-Splines (NURBS). Isogeometric Reissner-Mindlin shell has also been well studied in ~\cite{benson2010isogeometric, benson2011large}.
\citet{cirakortiz2000} developed a $C^1$ conforming discretisation based on Loop subdivision surfaces for an elastic Kirchhoff-Love shell formulation and applied it to hyperelastic thin shells~\cite{Cirak:2001aa}. Subdivision surfaces is an alternative to NURBS surfaces and is a mature geometry modelling method widely used in the animation and gaming industry. An attractive feature of subdivision surfaces is that they can be evaluated using spline functions, while retaining a simple polygonal mesh data structure able to represent complex geometries. Extraordinary vertices in the mesh allow for local refinement and patch conforming, both challenges faced by NURBS. Subdivision surfaces shell formulations have been extended to applications including shell fracture~\cite{cirak2005cohesive}, shape optimisation~\cite{BANDARA201862, chen2020acoustic}, fluid-structure interaction~\cite{cirak2007large}, non-manifold geometry~\cite{cirak2011subdivision} and structural-acoustic analysis~\cite{liu2018isogeometric}.  The ability of subdivision surfaces to analyse thin shells underpins the analysis of the electromechanical coupled thin shells presented here.

Applications for piezoelectric shells, such as resonators, actuators and energy harvesters, often involve the structural dynamics. Thus, understanding the effect of electroelastic coupling on the vibration mode of piezoelectric structures is critical. The coupling effect will influence the lattice structure of the piezoelectric material and enhance the stiffness of such structure via the so-called "piezoelectric stiffening" effect~\cite{johannsmann2015piezoelectric}. Thus, the natural frequencies of vibration modes increase. This effect is used in laminated beams~\cite{waisman2002active} and plates~\cite{donadon2002stiffening} with piezoelectric actuators to enhance their stiffness. However, the "piezoelectric stiffening" effect of piezoelectric thin shells with complex geometry is seldom studied. This work provides a numerical analysis tool for understanding these effects in piezoelectric thin shells.
 
The proposed method adopts Catmull-Clark subdivision surfaces to formulate a novel isogeometric Galerkin approach to analyse piezoelectric thin shells with arbitrary geometries. The formulation for analysing the thin shell considers three different electric conditions, which are no electrodes, prescribed voltage with electrodes and short-circuited electrode. These are comprehensively derived and summarised. A method to tailor the natural frequency of piezoelectric curved plate by changing its curvature is also presented. In addition, the "piezoelectric stiffening" effect of piezoelectric thin shells with complex geometry is also examined. 

This contribution is organised as follows. Section~\ref{sec:notations} introduces the notation and defines various coordinate systems used throughout the manuscript. Section~\ref{sec:shell} illustrates the kinematics of Kirchhoff-Love shells. Section~\ref{sec:subd} briefly reiterates the theory of Catmull-Clark subdivision surfaces. Energy considerations required for piezoelectric thin shells are presented in Section~\ref{sec:energy_density} and~\ref{sec:energy}. Hamilton's variational principle is applied and the resulting weak form of the governing equations of piezoelectric shells derived in Section~\ref{sec:variational_setting}. Section~\ref{sec:system_equation} discretises the weak form of the governing equations using Catmull-Clark subdivision bases resulting in the discrete system of equations. Finally, Section~\ref{sec:numerical} presents three numerical examples to demonstrate the ability of the proposed piezoelectric thin shell method to deal with various geometries and a range of mechanical and coupled problems.

\section{Notations}
\label{sec:notations}
\paragraph{Brackets:}
Two types of brackets are used. Square brackets $[ ]$ are used to clarify the order of operations in an algebraic expression. Circular brackets $( )$ are used to denote the parameters of a function. If brackets are used to denote an interval then $( )$ stands for an open interval and $[ ]$ is a closed interval.
\paragraph{Symbols:}
A variable typeset in a normal weight font represents a scalar. A bold weight font denotes a first or second-order tensor. An overline indicates that the variable is defined with respect to the reference configuration and if absent, the variable is defined with respect to the current (deformed) configuration. A scalar variable with superscript or subscript indices normally represents the components of a vector or second-order tensor. Upright font is used to denote matrices and vectors.

Indices $i,j,k,\dots$ vary from $1$ to $3$ while $a, b, c,\dots$, used as surface variable components, vary from $1$ to $2$. Einstein summation convention is used throughout.

The comma symbol in a subscript represents partial derivative, for example, $A_{,b}$ is the partial derivative of $A$ with respect to its $b^{\text{th}}$ component. $\nabla(\bullet)$ is the three-dimensional gradient operator. 
\paragraph{Coordinates:}
$ \mathbf{c}_i$ represent the basis vectors of an orthonormal system in three-dimensional Euclidean space and $x,y$ and $z$ are its components. $\boldsymbol{\xi}_i$ denote the orthonormal basis vectors in the local element space and $\xi,\eta$ and $\zeta$ are its coordinate components. The three covariant basis vectors for a surface point are denoted as $\mathbf a_i$, where $\mathbf a_1, \mathbf a_2$ are two tangential vectors and $\mathbf a_3$ is the normal vector. 
% $\mathbf{t}_i$ are the transformation vectors used to calculate the covariant components of the material tensors.

\begin{nolinenumbers}
\section{Kirchhoff-Love shell kinematics}
\label{sec:shell}
% \subsection{Motion}
The Kirchhoff-Love hypothesis can be applied to three-dimensional structures in which one dimension is much smaller than the other two. Important examples include plates and shells. It is assumed that lines perpendicular to the mid-surface remain straight and perpendicular to the mid-surface after deformation (see Figure~\ref{fig:shell_definition}).
\begin{figure}[h]
\centering
  \includegraphics[width=0.6\linewidth]{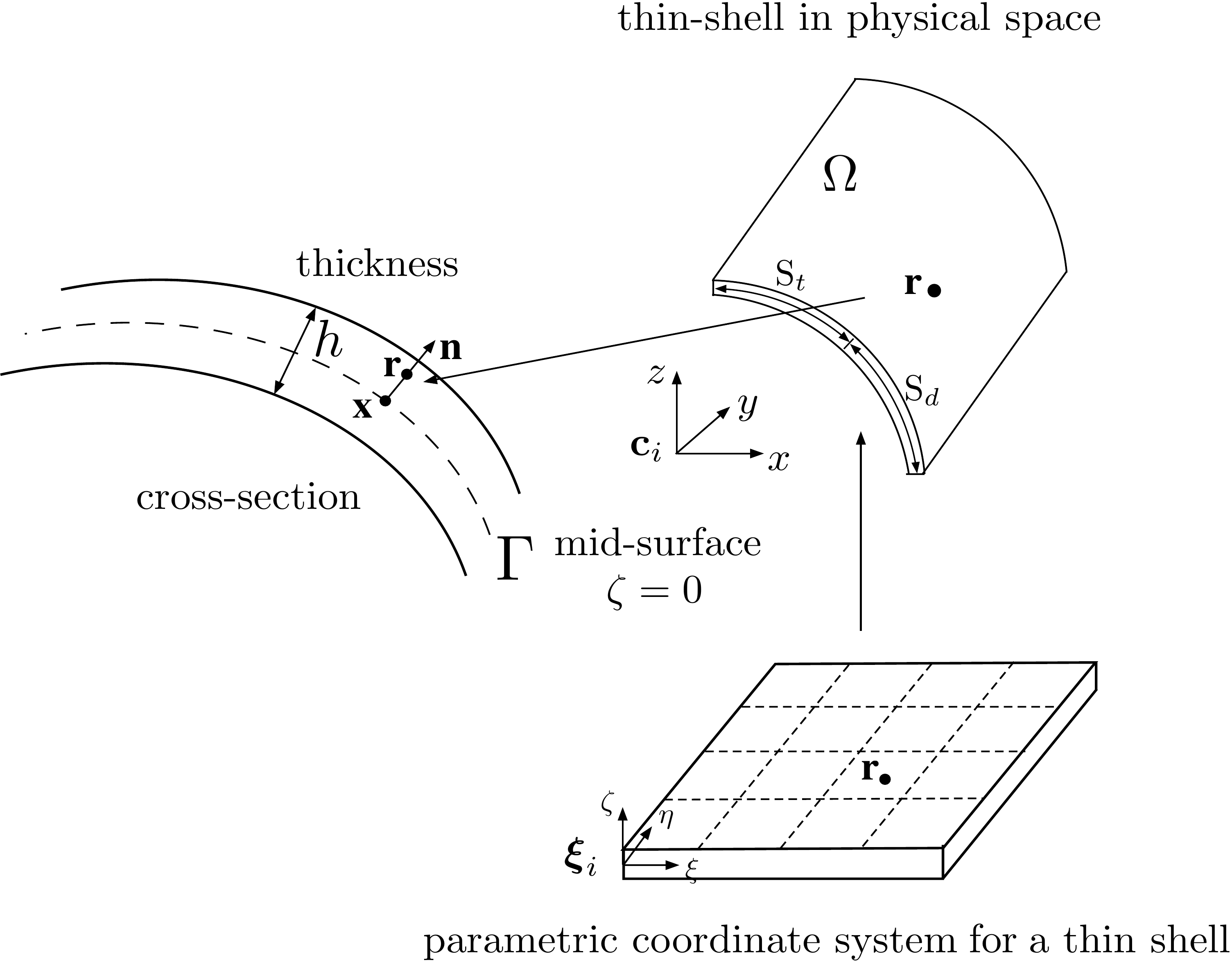}
\caption{Kirchhoff-Love shell coordinates}
\label{fig:shell_definition}
\end{figure}
The shell occupies the physical domain $\Omega$ and has a uniform thickness $h$. The thickness does not change upon deformation. The mid-surface of the shell is denoted by $\Gamma$. 
\begin{figure}[h]
\centering
  \includegraphics[width=0.6\linewidth]{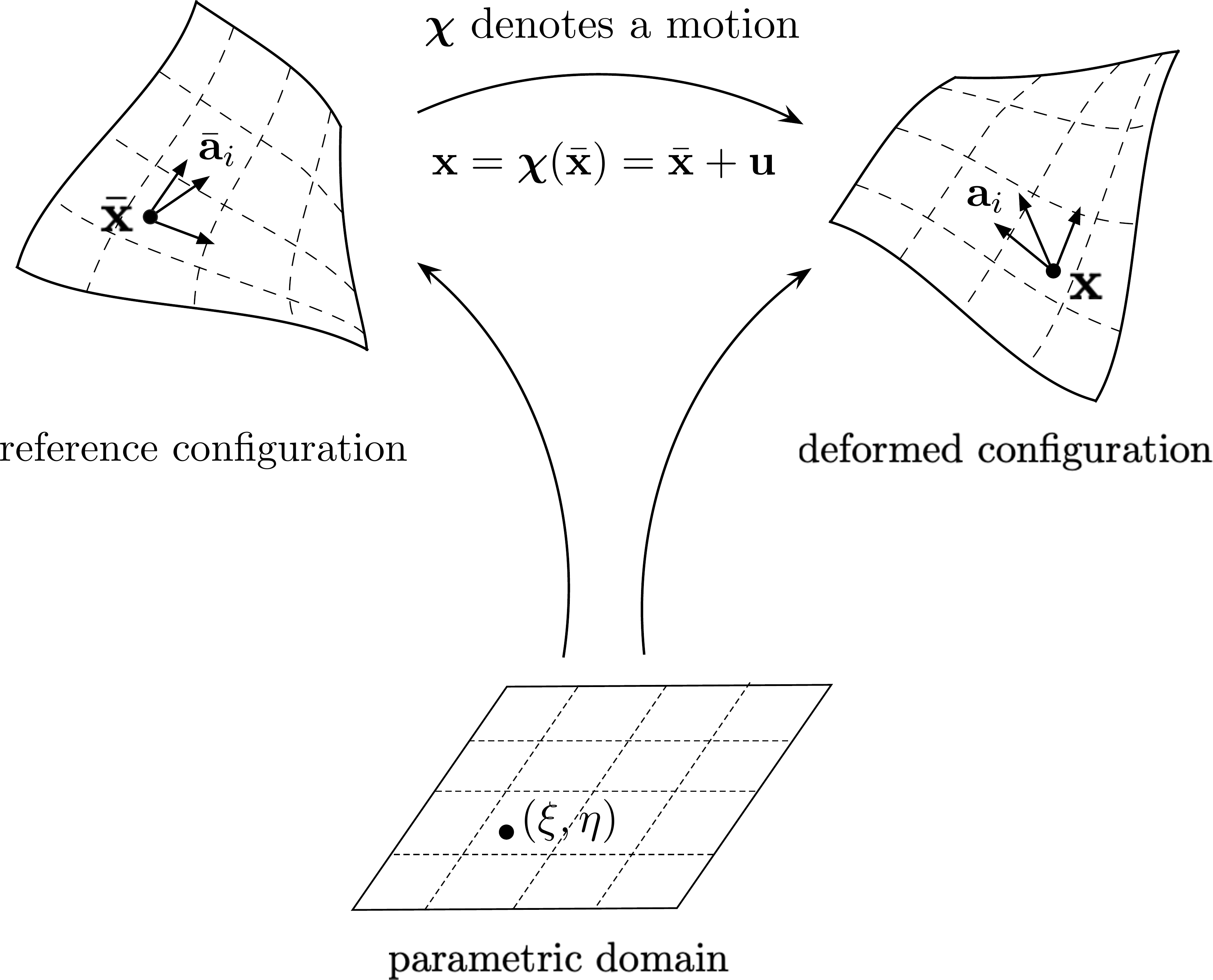}
\caption{Reference and deformed configurations for the mid-surface of a Kirchhoff-Love shell.}
\label{fig:shell_deformed_reference}
\end{figure}
Figure~\ref{fig:shell_deformed_reference} shows the reference and deformed configurations of the mid-surface. Points on the mid-surface in the reference and the deformed configurations are denoted by $\bar{\mathbf x}$ and $\mathbf x$. respectively, and are obtained as map from the parametric coordinates $\xi$ and $\eta$. The position vector of a point in the deformed configuration $\mathbf r$ is computed using the mid-surface point $\mathbf x$ and the normal vector $\mathbf n$ as
\begin{equation}
\mathbf r (\xi,\eta,\zeta) = \mathbf x (\xi,\eta) + \zeta \mathbf{n} (\xi,\eta),
\end{equation}
where $\zeta \in [-{h}/{2}, {h}/{2}]$. A mid-surface point in the deformed configuration ${\mathbf x}$ can be expressed as
\begin{equation}
\mathbf x = \bar{\mathbf x} + \mathbf u,
\label{eq:reference_to_deformed}
\end{equation}
where $\mathbf u$ denotes the displacement.
\subsection{Green-Lagrangian strain tensor}
The covariant basis vectors of the tangent plane of the mid-surface in the reference and the deformed configurations are defined by
\begin{equation}
\bar{\mathbf a}_1 = \frac{\partial \bar{\mathbf x}}{\partial \xi}, \quad  \bar{\mathbf a}_2 = \frac{\partial \bar{\mathbf x}}{\partial \eta},
\quad \text{and} \quad
{\mathbf a}_1 = \frac{\partial {\mathbf x}}{ \partial \xi}, \quad  {\mathbf a}_2 = \frac{\partial {\mathbf x}}{\partial \eta}.
\end{equation}
Thus, the normal vectors in the two configurations can be computed as 
\begin{equation}
\bar{\mathbf n} = \bar{\mathbf a}_3 = \frac{\bar{\mathbf a}_1 \times \bar{\mathbf a}_2}{\bar{J}}, \quad \text{and} \quad  {\mathbf n} = {\mathbf a}_3 = \frac{{\mathbf a}_1 \times {\mathbf a}_2}{{J}},
\end{equation}
where $\bar J$ and $J$ are the respective Jacobians given by
\begin{equation}
\bar{J} = |\bar{\mathbf a}_1 \times \bar{\mathbf a}_2|,\quad \text{and} \quad J = | {\mathbf a}_1 \times {\mathbf a}_2|.
\end{equation}
Thus, the covariant components of the metric tensor for the mid-surface points $\bar{\mathbf x}$ and $\mathbf x$ are respectively given by
\begin{equation}
\bar{a}_{ij} = \bar{\mathbf a}_i \cdot \bar{\mathbf a}_j, \quad \text{and} \quad {a}_{ij} = {\mathbf a}_i \cdot {\mathbf a}_j.
\end{equation}
The contravariant metric tensors are defined by
\begin{equation}
\bar{a}^{ik}\bar{a}_{kj} = \delta^{i}_{j}, \quad \text{and} \quad a^{ik}a_{kj} = \delta^{i}_{j},
\label{eq:co_and_contra_metric}
\end{equation}
where $\delta^{i}_{j}$ denotes the Kronecker Delta. 
The three-dimensional covariant basis vectors for the shell in the reference and the deformed configurations are respectively given by
\begin{equation}
\bar{\mathbf g}_1 = \frac{\partial \bar{\mathbf r}}{\partial \xi} = \bar{\mathbf a}_1 + \zeta \bar{\mathbf a}_{3,1}, \quad \bar{\mathbf g}_2 = \frac{\partial \bar{\mathbf r}}{\partial \eta} = \bar{\mathbf a}_2 + \zeta \bar{\mathbf a}_{3,2}, \quad \bar{\mathbf g}_3 = \frac{\partial \bar{\mathbf r}}{\partial \zeta} = \bar{\mathbf a}_3,
\label{eq:cov_tensors_r}
\end{equation}
and 
\begin{equation}
{\mathbf g}_1 = \frac{\partial {\mathbf r}}{\partial \xi} = {\mathbf a}_1 + \zeta {\mathbf a}_{3,1}, \quad 
{\mathbf g}_2 = \frac{\partial {\mathbf r}}{\partial \eta} = {\mathbf a}_2 + \zeta {\mathbf a}_{3,2}, \quad 
{\mathbf g}_3 = \frac{\partial {\mathbf r}}{\partial \zeta} = {\mathbf a}_3,
\label{eq:cov_tensors_d}
\end{equation}
where $(\bullet)_{,1}$ and  $(\bullet)_{,2}$ represent the partial differentials with respect to $\xi$ and $\eta$, respectively. The components of the covariant metric tensors are defined by
\begin{equation}
\bar{g}_{ij} = \bar{\mathbf g}_i \cdot \bar{\mathbf g}_j\quad \text{and} \quad g_{ij} = \mathbf g_i \cdot \mathbf g_j,
\end{equation}
which allows one to define the Green-Lagrange strain tensor $\mathbf {S}_{n}$ as
\begin{equation}
\mathbf {S}_{n} \coloneqq \frac{1}{2} [g_{ij} - \bar{g}_{ij}]\bar{\mathbf{g}}^i \otimes \bar{\mathbf{g}}^j,
\label{eq:green-lagrange}
\end{equation}
where $\bar{\mathbf{g}}^i$ denote the contravariant basis vectors defined by
\begin{equation}
\bar{\mathbf{g}}^i \cdot \bar{\mathbf{g}}_j = \delta^i_j.
\end{equation}
%\begin{equation}
%S_{ij} = \frac{1}{2} [g_{ij} - \bar{g}_{ij}].
%\label{eq:green-lagrange}
%\end{equation}
\subsection{Linearisation and simplification of the strain tensor}
On substituting equations~\eqref{eq:cov_tensors_r} and~\eqref{eq:cov_tensors_d} into~\eqref{eq:green-lagrange} and ignoring higher-order terms, the Green-Lagrange strain tensor linearised in $\zeta$ follows as
\begin{equation}
\mathbf S = \mathbf A + \zeta \mathbf B.
\label{eq:strain_decomposation}
\end{equation}
The components of the tensors $\mathbf A $ and $\mathbf B$ are
$\alpha_{ij}$ and $\beta_{ij}$, respectively, with $\alpha_{13}$ and $\alpha_{23}$ measuring the shearing in the normal direction $\bar{\mathbf a}_3$, and are zero under the Kirchhoff-Love assumption. The stretching in normal direction is given by $\alpha_{33} = 0$ and it vanishes due to the assumption that the thickness does not change with deformation. Similarly, $\beta_{i3} = 0$ as the normal vector is perpendicular to the two basis vectors.
Thus, the two tensors $\mathbf A$ and $\mathbf B$ reduce to two-dimensional tensors in the subspace defined with two contravariant basis vectors as
\begin{equation}
\mathbf A \coloneqq \alpha_{ab} \,\mathbf{\bar{g}}^{a} \otimes \mathbf{\bar{g}}^{b}\quad \text{and} \quad \mathbf B \coloneqq \beta_{ab} \,\mathbf{\bar{g}}^{a} \otimes \mathbf{\bar{g}}^{b},
\end{equation}
where their components are computed as
\begin{align}
%\alpha_{\alpha \beta} &= \frac{1}{2} [\mathbf a_{\alpha} \cdot \mathbf a_{\beta} - \bar{\mathbf a}_{\alpha} \cdot \bar{\mathbf a}_{\beta}], \\
%\beta_{\alpha \beta} &= \mathbf a_{\alpha} \cdot \mathbf a_{3,\beta} -  \bar{\mathbf a}_{\alpha} \cdot \bar{\mathbf a}_{3,\beta},
{\alpha}_{a b} = \frac{1}{2} [\mathbf a_{a} \cdot \mathbf a_{b} - \bar{\mathbf a}_{a} \cdot \bar{\mathbf a}_{b}] \quad \text{and} \quad {\beta}_{a b} = \mathbf a_{a} \cdot \mathbf a_{3,b} -  \bar{\mathbf a}_{a} \cdot \bar{\mathbf a}_{3,b}.
\end{align}
The membrane strain components are denoted as  $\alpha_{a b}$ while the bending strain components $\beta_{a b}$ measure the change in the curvature of the shell. 
In order to compute the bending strain tensor, the product rule of differentiation is applied and the components expressed as
\begin{equation}
\beta_{a b} = \bar{\mathbf{a}}_{a,b} \cdot \bar{\mathbf a}_3 - {\mathbf{a}}_{a,b} \cdot {\mathbf a}_3.
\end{equation}
On substituting Equation~\eqref{eq:reference_to_deformed} into the membrane and bending strains, the components can eventually be computed to first order in $\mathbf u$ as
\begin{align}
\alpha_{a b} = &\frac{1}{2} [\bar{\mathbf a}_{a} \cdot \mathbf u_{,b} + {\mathbf u}_{,a} \cdot \bar{\mathbf a}_{b}], \\
\beta_{a b} = &-\mathbf u_{,{a b}} \cdot \bar{\mathbf a}_3 + \frac{1}{\bar{J}} \big[ \mathbf u_{,1} \cdot [\bar{\mathbf a}_{a , b} \times \bar{\mathbf a}_2 ] + \mathbf u_{,2} \cdot [\bar{\mathbf a}_{1} \times \bar{\mathbf a}_{a , b} ] \big] + \frac{\bar{\mathbf a}_3 \cdot \bar{\mathbf a}_{a , b}}{\bar J}\big[ \mathbf u_{,1} \cdot [\bar{\mathbf a}_{2} \times \bar{\mathbf a}_3 ] + \mathbf u_{,2} \cdot [\bar{\mathbf a}_{3} \times \bar{\mathbf a}_{1} ] \big] .
\end{align}
Thus, the linearised strain tensor $\mathbf S$ is computed using the covariant basis vectors along with the first and second derivatives of the displacement $\mathbf u$. 

\section{Catmull-Clark subdivision surfaces}
\label{sec:subd}
Kirchhoff-Love shells require that the test and trial functions of the Galerkin method are in the Hilbert space $H^2(\Omega)$~\cite{cirakortiz2000}. Hence a $C^1$ continuous discretisation is required. Conventional Lagrangian bases only provide $C^0$ continuity. Catmull–Clark subdivision surfaces~\cite{stam1998exact}, which adopt cubic B-splines as interpolating functions, display $C^2$ continuity everywhere except at the surface points related to extraordinary vertices~\cite{peters1998analysis}, where continuity is only $C^1$. Figure~\ref{fig:element_bases} shows an example of cubic B-splines for one dimensional elements. The Catmull-Clark subdivision surfaces adopt a tensor-product structure of two cubic B-splines to interpolate points on a two-dimensional surface. Figure~\ref{fig:subd_surface} shows a smooth surface constructed by successive subdivision from a coarse polygonal mesh using the Catmull-Clark subdivision scheme~\cite{catmull1978recursively}. 
\begin{figure}[]
\centering
  \includegraphics[width=0.5\linewidth]{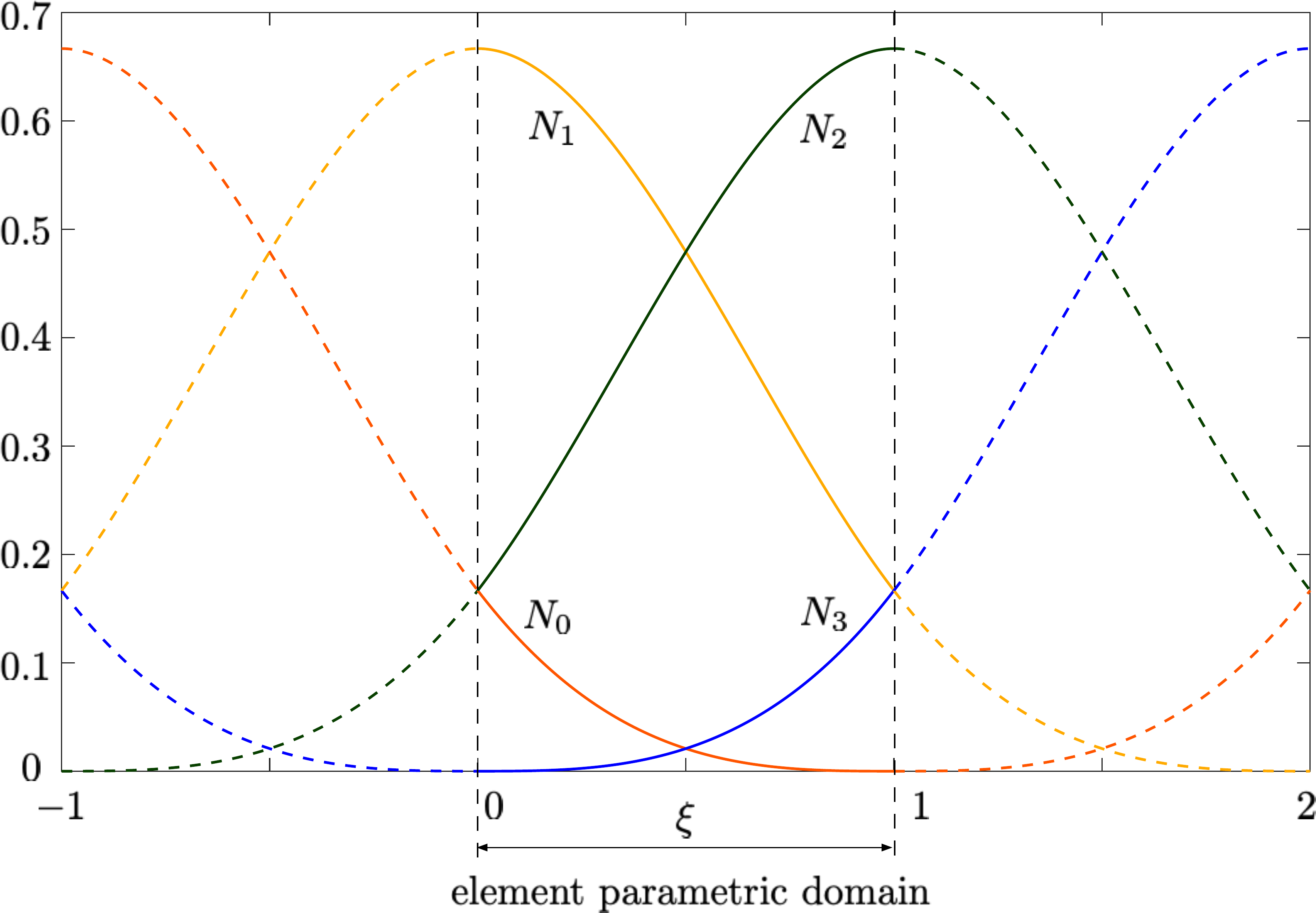}
\caption{An example of cubic B-splines in a one-dimensional parametric domain. Spline functions span multiple elements.}
\label{fig:element_bases}
\end{figure}
\begin{figure}[]
\centering
  \includegraphics[width=0.5\linewidth]{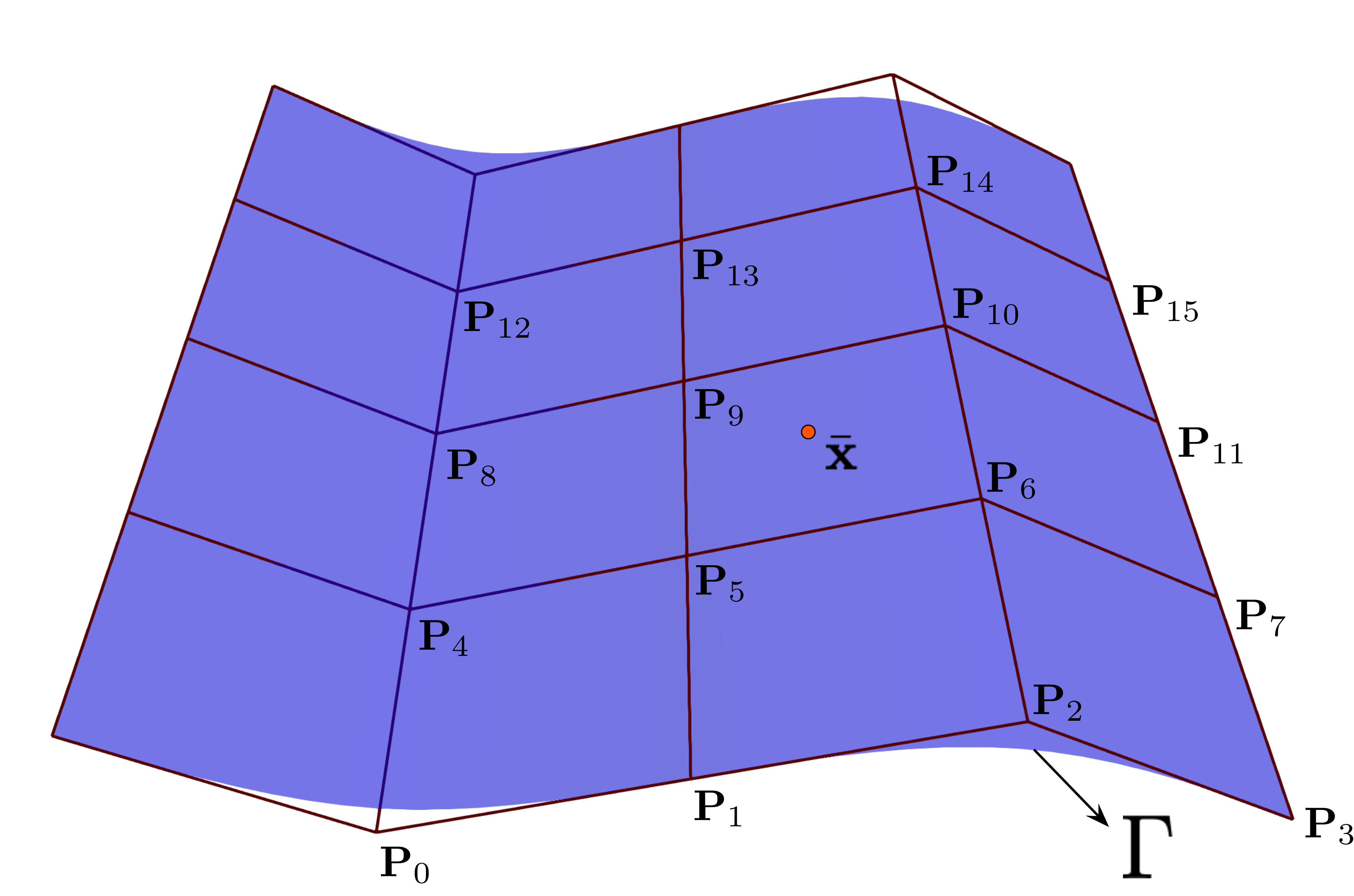}
\caption{The mid-surface of the shell $\Gamma$ is a Catmull-Clark subdivision surface constructed from a control polygonal mesh.}
\label{fig:subd_surface}
\end{figure}
The surface, composed of points $\bar{\mathbf x} \in \Gamma$, can be interpolated using the basis functions (cubic B-splines) and control points as
\begin{equation}
\bar{\mathbf x} =   \sum_{A=0}^{n_{b}-1} N^A \mathbf{P}_A,
\end{equation}
where $n_b$ is the number of basis functions. The $A^\text{th}$ basis function is denoted as $N^A$ and $\mathbf P_A$ denotes the $A^\text{th}$ control point. An element of a regular patch with $16$ basis functions is shown in Figure~\ref{fig:subd_surface}. We note that the control points are not necessarily on the surface $\Gamma$. Further details of an isogeometric Galerkin method using Catmull-Clark subdivision surfaces can be found in~\cite{liu2020assessment}.

\section{Piezoelectric shell formulation}
\label{sec:formulation}
\subsection{Energy densities}
\label{sec:energy_density}
The electric enthalpy density per unit volume for a coupled piezoelectric problem~\cite{mason1951piezoelectric, tiersten1967hamilton} is most generally defined by
\begin{align}
\mathcal{H}(\mathbf S, \mathbf E) = W_{\text{ela}} (\mathbf S) - W_{\text{piezo}}(\mathbf S, \mathbf E) - W_\text{elec}( \mathbf E).
\end{align}
The electric enthalpy density contains the elastic energy density $W_\text{ela}$, the piezoelectric energy density $W_\text{piezo}$ and the electric energy density $W_\text{elec}$. The electric field is denoted as $\mathbf E$.
The piezoelectric and electric energy densities are expressed as
\begin{equation}
W_\text{piezo} (\mathbf S, \mathbf E) = \mathbf E \cdot [\mathbf e : \mathbf S] = e^{ijk}E_i [\alpha_{jk} + \zeta \beta_{jk}],
\end{equation}
and

\begin{equation}
W_\text{elec} ( \mathbf E) =  \frac{1}{2}[\boldsymbol{\kappa} \cdot \mathbf E]\cdot \mathbf E = \frac{1}{2}\kappa^{ij}E_iE_j,
\end{equation}

respectively.
The components of the third-order piezoelectric tensor $\mathbf e$ are $e^{ijk}$ while $\kappa^{ij}$ are the components of the second-order dielectric tensor $\boldsymbol{\kappa}$. 
Since the structure is thin and has uniform thickness, we introduce the quadratic elastic strain energy density per unit area $\widetilde{W}_\text{ela}$ for the Kirchhoff-Love shell as
\begin{equation}
\widetilde{W}_\text{ela}( \mathbf S) = \int_{-\frac{h}{2}}^{\frac{h}{2}} W_{\text{ela}} (\mathbf S) \, \mathrm d \zeta.
\end{equation}
%If the material is assumed as isotropic, the elastic strain energy density per unit area consists of the membrane and bending parts~\cite{cirakortiz2000} as
%\begin{equation}
%\widetilde{W}_\text{ela}( \mathbf A, \mathbf B) = \frac{1}{2} \frac{Eh}{1-\nu^2}\left[ [\mathbf A : \mathbf H : \mathbf A] + \frac{h^2}{12} [\mathbf B : \mathbf H : \mathbf B] \right]= \frac{1}{2} \frac{Eh}{1-\nu^2} H^{abcd}\alpha_{ab}\alpha_{cd} + \frac{1}{2} \frac{Eh^3}{12[1-\nu^2]} H^{abcd}\beta_{ab}\beta_{cd} ,
%\end{equation}
%where $E$ and $\nu$ are the Young's Modulus and Poisson's ratio, respectively. %Here $i,j,k,l$ vary from $1$ to $2$. 
%$H^{abcd}$ denotes the component of the fourth order tensor $\mathbf H$ computed from the contravariant metric tensors as
%\begin{equation}
%H^{abcd} = \nu\bar{a}^{ab}\bar{a}^{cd} + \frac{1}{2} [1-\nu][\bar{a}^{ac}\bar{a}^{bd} + \bar{a}^{ad}\bar{a}^{bc}].
%\end{equation}

The piezoelectric material is normally anisotropic due to the interaction between the mechanical and electrical states in crystalline materials with no inversion symmetry. Thus, with $\mathbf S = \mathbf A + \zeta \mathbf B$, one defines a general formulation for the elastic energy density per unit area by
\begin{equation}
\widetilde{W}_\text{ela}( \mathbf S) = \widetilde{W}_\text{ela}( \mathbf A, \mathbf B) = \frac{h}{2} \left[ [\mathbf A : \mathbf C : \mathbf A] + \frac{h^2}{12} [\mathbf B : \mathbf C : \mathbf B] \right],
%= \frac{h}{2} C^{abcd}\alpha_{ab}\alpha_{cd} + \frac{1}{2} \frac{h^3}{12} C^{abcd}\beta_{ab}\beta_{cd} ,
\end{equation}
where $\mathbf C$ is the fourth-order elastic tensor which can be defined using the covariant base vectors by
%If the material properties are defined in the curvilinear coordinate system and the material is assumed to be transverse isotropic, the components of $\mathbf C$ can be expressed as
\begin{equation}
\mathbf C =  {C}^{ijkl} \bar{\mathbf g}_i \otimes \bar{\mathbf g}_j \otimes \bar{\mathbf g}_k \otimes \bar{\mathbf g}_l = \tilde{C}^{mnop} {\mathbf t}_m \otimes {\mathbf t}_n \otimes {\mathbf t}_o \otimes {\mathbf t}_p.
%= \tilde{C}^{ijkl} \boldsymbol{\xi}_i \otimes \boldsymbol{\xi}_j \otimes \boldsymbol{\xi}_k \otimes \boldsymbol{\xi}_l.
\end{equation}
The preferable anisotropy directions of the piezoelectric material is denoted as $\mathbf t _m$. Therefore, the components of the elasticity tensor are related by
\begin{equation}
{C}^{ijkl} = \tilde{C}^{mnop} [\bar{\mathbf g}^i \cdot \mathbf t_m] [\bar{\mathbf g}^j \cdot \mathbf t_n] [\bar{\mathbf g}^k \cdot \mathbf t_o] [\bar{\mathbf g}^l \cdot \mathbf t_p].
\end{equation}

%The transformation between two basis system is defined by
%\begin{equation}
%\boldsymbol{\xi}_i = g_i^j \bar{\mathbf g}_j,
%\end{equation}
%where $g^i_j$ can be considered as the $j^\text{th}$ component of contravariant basis $\bar{\mathbf g}^i$. Thus

%\begin{equation}
%{C}^{abcd} = \tilde{C}^{abcd} \cdot  \bar g^{ab} \cdot \bar g^{cd}.
%\end{equation}

%where
%\begin{equation}
%\mathbf t_i = \frac{\bar{\mathbf g}_i}{|\bar{\mathbf g}_i|}.
%\end{equation}

\subsection{Kinetic energy}
\label{sec:energy}
The kinetic energy of a Kirchhoff-Love thin shell is defined by
\begin{equation}
\Pi_{\text{kin}} = \frac{\rho h}{2} \int_\Gamma  \left[\frac{\partial u_i}{\partial t}\right]^2 \,\mathrm d \Gamma
\end{equation}
where $\rho$ denotes the mass density per unit volume which is here assumed constant.
\subsection{Electric enthalpy}
The total electric enthalpy of the system is composed of three parts:
\begin{align}
\mathfrak{E} ( \mathbf S, \mathbf E) = \Pi_\text{ela}( \mathbf S)  - \Pi_\text{piezo} ( \mathbf S, \mathbf E) - \Pi_\text{elec}( \mathbf E),
\end{align}
where  $\Pi_\text{piezo}$ is the piezoelectric energy. The dielectric energy is denoted as $\Pi_\text{elec}$ and the elastic energy is defined by
\begin{equation}
 \Pi_\text{ela}( \mathbf{S}) = \int_{\Gamma} \widetilde{W}_\text{ela} ( \mathbf S) \, \mathrm d \Gamma.
\end{equation}
% \begin{figure}[]
% \centering
%   \includegraphics[width=0.4\linewidth]{piezo_shell_phi.pdf}
% \caption{Assumption of the piezoelectric shell. The top and bottom surfaces of the shell are connected with the ground.}
% \label{fig:piezo_shell_phi}
% \end{figure}
%  {\color{red} 
To consider the piezoelectric and the dielectric energy for a thin shell formulation, a power series expansion is applied to the electric potential with respect to the thickness coordinate $\zeta$~\cite{tiersten1993equations} and the first three terms are retained, that is
 \begin{equation}
\phi(\mathbf{r}(\xi,\eta,\zeta)) \approx \phi^{(0)}(\mathbf x(\xi,\eta)) + \zeta  \phi^{(1)}(\mathbf x(\xi,\eta)) +   \left[\zeta^2 - \left[\frac{h}{2}\right]^2\right]   \phi^{(2)}(\mathbf x(\xi,\eta)).
\label{eq:phi_expansion}
 \end{equation}
 \begin{figure}[]
\centering
  \includegraphics[width=\linewidth]{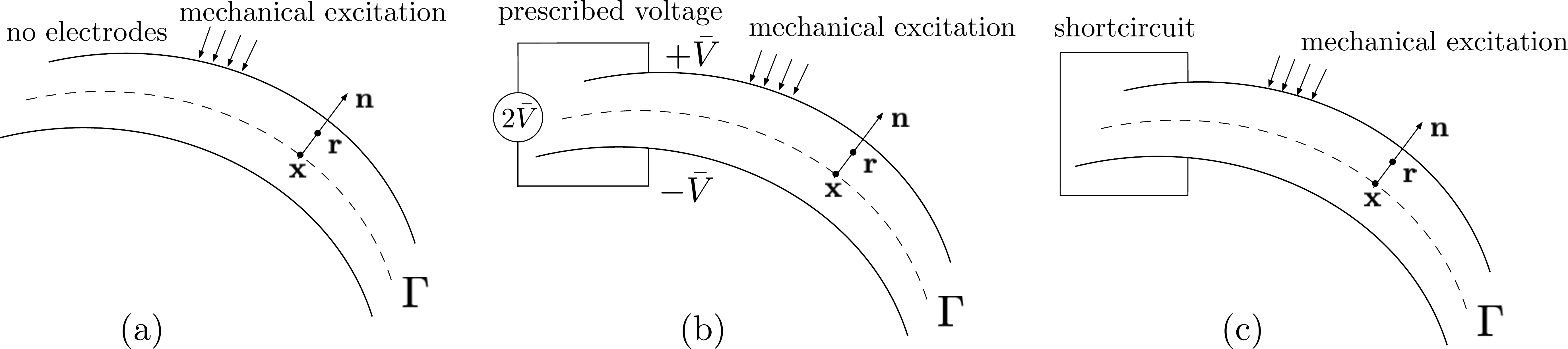}
\caption{Three electric setups for piezoelectric shells. (a) A shell in free space with no electrodes. (b) An electroded shell with symmetrically prescribed voltage. (c) A special case in which the electrodes are short-circuited.}
\label{fig:three_shells}
\end{figure}
%  {\color{red}
The electric field is computed as
\begin{equation}
\mathbf E = - \nabla \phi,
\end{equation}
and it can be expressed using contravariant basis vector as
\begin{equation}
\mathbf E =  E_i \bar{\mathbf{g}}^i,
\end{equation}  
Due to the large relative permittivity of piezoelectric materials, the electric field in the surrounding free space is neglected.
The energy contributions and hence the coupling effect depends on the configuration of the piezoelectric shell structure. Unelectroded and electroded shells along with a special short-circuited case, as displayed in Figure~\ref{fig:three_shells} are three options considered here.
\paragraph{$\bullet$ Shell with no electrodes }
In this case, the shell structure is assumed to be embedded in free space, thus $\phi^{(1)} \neq 0$ and $\phi^{(2)} \neq 0$. 
Upon substituting expression~\eqref{eq:phi_expansion}, the contravariant coefficients of the electric field are calculated as
\begin{align}
{E}_1 &= -\frac{\partial \phi}{\partial \xi} = - \phi^{(0)}_{,\xi} - \zeta\phi^{(1)}_{,\xi} -\left[\zeta^2 - \frac{h^2}{4}\right]  \phi^{(2)}_{,\xi}, \quad {E}_2 =-\frac{\partial \phi}{\partial \eta} =  -\phi^{(0)}_{,\eta} - \zeta\phi^{(1)}_{,\eta} -\left[\zeta^2 - \frac{h^2}{4}\right]  \phi^{(2)}_{,\eta}, \nonumber \\ 
\quad {E}_3 &=-\frac{\partial \phi}{\partial \zeta} = -\phi^{(1)} -2\zeta \phi^{(2)}.
\label{eq:E_components}
\end{align}
The piezoelectric energy is expressed as
% {\color{red}
\begin{equation}
\Pi_\text{piezo}(\mathbf A, \mathbf B, \mathbf E) =\int_{\Omega} e^{ibc} E_i [\alpha_{bc} + \zeta\beta_{bc}] \,\mathrm d \Omega.
\label{eq:piezo_energy_1}
\end{equation}
On substituting expressions~\eqref{eq:strain_decomposation} and~\eqref{eq:E_components} into~\eqref{eq:piezo_energy_1}, the piezoelectric energy can be expressed as
\begin{align}
\Pi_\text{piezo}(\mathbf A, \mathbf B, \phi^{(0)}, \phi^{(1)}, \phi^{(2)}) = &-h\int_{\Gamma} e^{abc} \phi^{(0)}_{,a} \alpha_{bc} \,\mathrm d \Gamma -  \frac{h^3}{12}\int_{\Gamma} e^{abc} \phi^{(1)}_{,a} \beta_{bc} \,\mathrm d \Gamma - h \int_{\Gamma} e^{3bc} \phi^{(1)} \alpha_{bc} \, \mathrm d \Gamma \nonumber \\
&+\frac{h^3}{6}\int_{\Gamma} e^{abc} \phi^{(2)}_{,a} \alpha_{bc} \,\mathrm d \Gamma - \frac{h^3}{6} \int_{\Gamma} e^{3bc} \phi^{(2)} \beta_{bc} \, \mathrm d \Gamma.
\end{align}
The third order piezoelectric tensor is expressed either in the covariant basis or the local coordinate system as
%\begin{equation}
%e^{ibc} = \tilde{e}^{ibc}a^{bc}.
%\end{equation}
%\begin{equation}
%e^{ibc} = \tilde{e}^{ibc} |\bar{\mathbf{g}}^i| \bar{a}^{bc}.
%\end{equation}
\begin{equation}
\mathbf e = {e}^{ijk} \bar{\mathbf g}_i \otimes \bar{\mathbf g}_j \otimes \bar{\mathbf g}_k  = \tilde{e}^{lmn} {\mathbf t}_l \otimes {\mathbf t}_m \otimes {\mathbf t}_n,
%= \tilde{e}^{lmn} \boldsymbol{\xi}_l \otimes \boldsymbol{\xi}_m \otimes \boldsymbol{\xi}_n，
\end{equation}
%The transformation between two basis system is defined by
%\begin{equation}
%\boldsymbol{\xi}_i = g_i^j \bar{\mathbf g}_j,
%\end{equation}
%where $g^i_j$ can be considered as the $j^\text{th}$ component of contravariant basis $\bar{\mathbf g}^i$. 
with components related via
\begin{equation}
e^{ijk} = \tilde{e}^{lmn} [\bar{\mathbf g}^i \cdot \mathbf t_l] [\bar{\mathbf g}^j \cdot \mathbf t_m]  [\bar{\mathbf g}^k \cdot \mathbf t_n].
\end{equation}
In the present work, the piezoelectric material only polarises in the thickness direction, $\mathbf t_3 = \mathbf n$.  Then, the coefficients $e^{abc}$ can be considered as zeros. Thus two terms contribute to the piezoelectric energy, that is
\begin{equation}
\Pi_\text{piezo}(\mathbf A, \mathbf B, \phi^{(1)}, \phi^{(2)}) = - h \int_{\Gamma} e^{3bc} \phi^{(1)} \alpha_{bc} \, \mathrm d \Gamma 
 - \frac{h^3}{6} \int_{\Gamma} e^{3bc} \phi^{(2)} \beta_{bc} \, \mathrm d \Gamma.
\end{equation}
Because $\phi^{(0)}$ does not contribute to the piezoelectric energy, we conveniently set $\phi^{(0)} = 0$. Since the electric field in the surrounding free space is neglected, the electric energy is expressed as
\begin{equation}
\Pi_\text{elec}(\phi^{(1)},\phi^{(2)}) = \frac{h^3}{24} \int_{\Gamma} \kappa^{ab} \phi^{(1)}_{,a} \phi^{(1)}_{,b} \, \mathrm d \Gamma + \frac{h}{2} \int_\Gamma \kappa^{33} [\phi^{(1)}]^2 \, \mathrm d \Gamma +\frac{h^5}{60} \int_{\Gamma} \kappa^{ab} \phi^{(2)}_{,a} \phi^{(2)}_{,b} \, \mathrm d \Gamma + \frac{h^3}{6} \int_\Gamma \kappa^{33} [\phi^{(2)}]^2 \, \mathrm d \Gamma
\end{equation}
%  }
where the dielectric tensor is expressed in the covariant or the local coordinate systems as 
%\begin{equation}
%\kappa^{ij} = \tilde{\kappa}^{ij} |\bar{\mathbf{g}}^i| |\bar{\mathbf{g}}^j|.
%\end{equation}
\begin{equation}
\boldsymbol{\kappa} = \kappa^{ij}  \bar{\mathbf g}_i \otimes \bar{\mathbf g}_j = \tilde{\kappa}^{kl} \mathbf t_{k} \otimes \mathbf t_{l}
% = \tilde{\kappa}^{lm} \boldsymbol{\xi}_l \otimes \boldsymbol{\xi}_m,
\end{equation}
%Thus
%\begin{equation}
%\kappa^{ij} = \tilde{\kappa}^{ij} \cdot \bar g^{ij }.
%\end{equation}
with components related via
\begin{equation}
\kappa^{ij} = \tilde{\kappa}^{kl} [\bar{\mathbf g}^i \cdot \mathbf t_k] [\bar{\mathbf g}^j \cdot \mathbf t_l].
\end{equation}

 \paragraph{$\bullet$ Symmetrically prescribed voltage with electrodes }
Here we assume the shell is electroded on top and bottom surface with constant voltage $V_1$ and $V_2$, respectively. Thus, as the surface potential is constant for all $\mathbf x$ and the following relation must be satisfied
\begin{align}
\phi^{(0)} + \frac{h}{2}  \phi^{(1)}= V_1, \\
\phi^{(0)} - \frac{h}{2}  \phi^{(1)} = V_2.
\end{align}
Thus $\phi^{{0}}$ and $\phi^{(1)}$ are constants and computed as
\begin{align}
\phi^{(0)} = \frac{V_1 + V_2}{2},\\
\phi^{(1)} = \frac{V_1-V_2}{h}.
\end{align}
If the shell is symmetrically electroded with constant voltage, $V_1 = \bar{V}$ and $V_2 = -\bar{V}$, then $\phi^{(0)} \equiv 0$ and $\phi^{(1)} = {2\bar{V}}/{h}$. Equation~\eqref{eq:phi_expansion} thus becomes
\begin{equation}
\phi(\mathbf{r}(\xi,\eta,\zeta)) = \zeta\frac{2\bar{V}}{h} + \left[\zeta^2 - \left[\frac{h}{2}\right]^2\right]   \phi^{(2)}(\mathbf x(\xi,\eta)).
\end{equation}
Eventually, the contravariant coefficients of the electric field simplify to
\begin{equation}
{E}_1 = -\frac{\partial \phi}{\partial \xi} =  -\left[\zeta^2 - \frac{h^2}{4}\right]  \phi^{(2)}_{,\xi}, \quad {E}_2 =-\frac{\partial \phi}{\partial \eta} =  -\left[\zeta^2 - \frac{h^2}{4}\right]  \phi^{(2)}_{,\eta  }, \quad {E}_3 =-\frac{\partial \phi}{\partial \zeta} = -\frac{2\bar{V}}{h} -2\zeta \phi^{(2)}.
\label{eq:E_components_2}
\end{equation}
On substituting expressions~\eqref{eq:strain_decomposation} and~\eqref{eq:E_components_2} into equation~\eqref{eq:piezo_energy_1}, the piezoelectric energy is now expressed as
\begin{equation}
\Pi_\text{piezo}(\mathbf A, \mathbf B, \phi^{(2)}) = 
- \frac{h^3}{6} \int_{\Gamma} e^{3bc} \phi^{(2)} \beta_{bc} \, \mathrm d \Gamma - 2\bar{V}\int_{\Gamma} e^{3bc} \alpha_{bc} \mathrm d \Gamma.
\end{equation}
Furthermore, the electric energy is now expressed as
\begin{equation}
\Pi_\text{elec}(\phi^{(2)}) = \frac{2{\bar{V}}^2}{h}\kappa^{33}\int_{\Gamma}  \, \mathrm d \Gamma + \frac{h^5}{60} \int_{\Gamma} \kappa^{ab} \phi^{(2)}_{,a} \phi^{(2)}_{,b} \, \mathrm d \Gamma + \frac{h^3}{6} \int_\Gamma \kappa^{33} [\phi^{(2)}]^2 \, \mathrm d \Gamma.
\end{equation}
\paragraph{$\bullet$ Short-circuited electrodes }
A special electric condition can be obtained by short-circuiting the electrodes, thus $\bar V = 0$.
The piezoelectric energy is now expressed as
\begin{equation}
\Pi_\text{piezo}( \mathbf B, \phi^{(2)}) = 
% \frac{h^3}{6}\int_{\Gamma} e^{abc} \phi^{(2)}_{,a} \alpha_{bc} \,\mathrm d \Gamma 
% \Pi_\text{piezo}(\mathbf A, \mathbf B, \phi^{(2)}) = 
- \frac{h^3}{6} \int_{\Gamma} e^{3bc} \phi^{(2)} \beta_{bc} \, \mathrm d \Gamma,
\end{equation}
while the corresponding electric energy is given by
\begin{equation}
\Pi_\text{elec}(\phi^{(2)}) = \frac{h^5}{60} \int_{\Gamma} \kappa^{ab} \phi^{(2)}_{,a} \phi^{(2)}_{,b} \, \mathrm d \Gamma + \frac{h^3}{6} \int_\Gamma \kappa^{33} [\phi^{(2)}]^2 \, \mathrm d \Gamma.
\end{equation}
\\
The three electric conditions for the piezoelectric shell are summarised in Table~\ref{tab:three_electric_conditions}.
%{\color{red} We note that  $I$ and $J$ vary from $1$ to $3$ as they are indices in global Cartesian coordinate system, while $i$ and $j$ are indices in local coordinate system. Then, the piezoelectric and dielectric tensor are transformed as
%\begin{align}
%e^{Ijk} &= \sum_i\tilde{e}^{ijk} \tilde v_{i}^{I} \bar a, \\
%\kappa^{IJ} &= \tilde \kappa^{ij} \tilde v_{i}^{I}  \tilde v_{j}^{J}.
%\end{align}
%Where $\tilde{e}^{ijk}$ and $\tilde\kappa^{ij}$ are piezoelectric and dielectric tensor which are considered as material properties. $\bar a$ here is the contravariant metric tensor and $\tilde{\mathbf{v}}_{i}$ denotes the normalised base vectors such that
%\begin{equation}
%\tilde {\mathbf v}_1 = \frac{\bar{\mathbf{a}}_1}{|\bar{\mathbf{a}}_1|}, \quad \tilde {\mathbf v}_2 = \frac{\bar{\mathbf{a}}_2}{|\bar{\mathbf{a}}_2|}, \quad \tilde {\mathbf v}_3 = \bar{\mathbf{a}}_3.
%\end{equation}
%}

\begin{table}[]
\caption{Three different electric conditions applied to the top and bottom surfaces of piezoelectric shells.}
\begin{tabular}{lllll}
\hline
\multirow{2}{*}{Electric conditions} & \multicolumn{3}{l}{Electric functions}  & \multirow{2}{*}{Summary} \\
& $\phi^{(1)}$ & $\phi^{(2)}$ & $\bar{V}$  & \\ \hline
Unelectroded                         & \checkmark   & \checkmark  & \ding{55}  & \begin{tabular}[c]{@{}l@{}}The shell is embedded in free space, the linear potential function $\phi^{(1)}$ \\ is a variable coupled with the membrane strain. The quadratic \\ potential function $\phi^{(2)}$ is a variable coupled with the bending strain.\end{tabular}                                                                                                \\
\\
Prescribed voltage                   & \ding{55}     & \checkmark  & \checkmark & \begin{tabular}[c]{@{}l@{}}The top and bottom surfaces are electroded and a constant potential\\ difference $2\bar{V}$ is symmetrically applied between them. Thus a linear\\ potential is prescribed which induces a global membrane strain.\\ Only the quadratic potential function $\phi^{(2)}$ remains as a variable.\\ If $\bar{V}$ is large, the quadratic coupling term can be ignored and the\\ problem reduces to a one-way coupling.\end{tabular} \\
\\
Short-circuited                      & \ding{55}     & \checkmark  & \ding{55}    & \begin{tabular}[c]{@{}l@{}}The top and bottom surfaces are electroded and short-circuited, $\bar{V} = 0$. \\ Only the quadratic potential function $\phi^{(2)}$ is a variable.\end{tabular}                                                                                                                                                \\ \hline
\end{tabular}
\label{tab:three_electric_conditions}
\end{table}

\subsection{Stress relaxation for thin-shells}
The stress tensor is denoted as $\boldsymbol \sigma = \sigma^{ij} \bar{\mathbf g}_i \otimes \bar{\mathbf g}_j$ with components given by
\begin{equation}
\sigma^{ij} = C^{ijkl}S_{kl} - e^{kij}E_k,
\end{equation}
where $S_{ij}$ denote the components of strain tensor $\mathbf S$. Since the thin shell assumption is adopted in the current work, the dominant stress components are the in-plane terms $\sigma^{ab}$. The Kirchhoff-Love assumption implies the shear stresses and strains are both neglected, thus the $\sigma^{33}$ and $S^{33}$ are the only non-zero out-of-plane components. Stress relaxation is performed by setting $\sigma^{33} = 0$, that is
\begin{equation}
\sigma^{33} = {C}^{33ij} S_{ij} - e^{i33}{E}_i = 0. 
\end{equation}
Since $S_{i3}$ and $S_{3j}$ are $0$, the remaining out-of-plane strain is computed as 
\begin{equation}
S^{33} = -\frac{1}{{C}^{3333}} [ {C}^{33ab}S_{ab} - e^{i33}E_i ].
\end{equation}
The elastic, piezoelectric and dielectric tensors are modified accordingly as
\begin{align}
\hat{C}^{abcd} = {C}^{abcd}  - \frac{ C_{ab33}  C_{33cd}}{ C^{3333}}, \quad 
% \nonumber \\
\hat{e}^{ijk} = {e}^{ijk} - \frac{ e_{i33}  C_{33jk}}{ C^{3333}},
\quad \text{and} \quad
\hat{\kappa}^{ij} = {\kappa}^{ij} + \frac{ e_{i33}  e_{j33}}{ C^{3333}}.
\end{align}
Those modified tensors are used in the following formulation.
\subsection{External energy}
The external energy contains the elastic and dielectric parts expressed as
% \begin{equation}
% \Pi_\text{ext}(\mathbf u, \phi^{(1)},\phi^{(2)}) = \Pi_\text{ext}^\text{ela} (\mathbf u) + \Pi_\text{ext}^\text{elec} (\phi^{(1)},\phi^{(2)}).
% \end{equation}
\begin{equation}
\Pi_\text{ext}(\mathbf u, \phi) = \Pi_\text{ext}^\text{ela} (\mathbf u) + \Pi_\text{ext}^\text{elec} (\phi).
\end{equation}
The external elastic energy is computed as
\begin{align}
\Pi_\text{ext}^\text{ela} (\mathbf u) = h\int_\Gamma  b_i u_i \,\mathrm d \Gamma + h \int_{S_t} \tau_i u_i \, \mathrm d S_t, 
\end{align}
where $b_i$ denotes the components of a body force and $\tau_i$ the components of a prescribed traction. $S_t \in \partial \Gamma$ represents the line where the traction is applied. \\
The external electric energy is only a function of $\phi^{(2)}$ since
\begin{equation}
\Pi_\text{ext}^\text{elec} (\phi^{(2)}) = \frac{h^3}{6} \int_\Gamma q \phi^{(2)} \, \mathrm d \Gamma + \frac{h^3}{6} \int_{S_d} \omega \phi^{(2)} \,\mathrm d S_d,
\label{eq:elec_external}
\end{equation}
where $q$ is the volume charge density and $\omega$ is the surface charge density on the cross-section of the shell. $S_d \in \partial \Gamma$ represents the line where the electric loads are applied. We note that the piezoelectric shell is made of a dielectric material and is thus an insulator. Since its cross-section is very thin, both volume and surface charge are difficult to apply in practical devices. The expression~\eqref{eq:elec_external} is kept in the formulation for the sake of completeness but the contribution neglected the subsequent numerical examples.

%The total potential of the system is expressed explicitly as
%\begin{align}
%\Pi = &\frac{\rho h}{2} \int_\Gamma  \left[\frac{\partial u_i}{\partial t}\right]^2 \,\mathrm d \Gamma +  \int_{\Gamma} \left[\frac{1}{2} \frac{Eh}{1-\nu^2} H^{ijkl}\alpha_{ij}\alpha_{kl} + \frac{1}{2} \frac{Eh^3}{12[1-\nu^2]} H^{ijkl}\beta_{ij}\beta_{kl}\right] \, \mathrm d \Gamma \nonumber \\
%&+ \frac{h^3}{12}\int_{\Gamma} e^{Ijk} (\nabla_\Gamma \varphi)_I \alpha_{jk} \,\mathrm d \Gamma + \frac{h^3}{6} \int_{\Gamma} e^{Ijk} \varphi \bar{n}_I \beta_{jk} \, \mathrm d \Gamma \nonumber \\
%&- \frac{h^5}{160} \int_{\Gamma} \kappa^{IJ} (\nabla_\Gamma \varphi)_I (\nabla_\Gamma \varphi)_J \, \mathrm d \Gamma - \frac{h^3}{6} \int_\Gamma \kappa^{IJ} \varphi^2 \bar{n}_I\bar{n}_J \, \mathrm d \Gamma \nonumber \\
% &-h\int_\Gamma  b_I u_I \,\mathrm d \Gamma - h \int_{S_t} t_I u_I \, \mathrm d S_t  + \frac{h^3}{12} \int_\Gamma q \varphi \, \mathrm d \Gamma + \frac{h^3}{12} \int_{S_d} \omega \varphi \,\mathrm d S_d.
%\end{align}

\subsection{Variational setting}
\label{sec:variational_setting}
Hamilton's principle, ignoring dissipative mechanisms, states that the variation of the action integral of a piezoelectric shell is zero, thus
\begin{equation}
\delta \int_{t_0}^{t_1} L(\mathbf u , \psi, \varphi) \,\mathrm d t = 0,
\label{eq:variation_l}
\end{equation}
where $\psi$ and $\varphi$ are henceforth used to denote $\phi^{(1)}$ and $\phi^{(2)}$ to simplify the notation. $\delta(\bullet)$ represents the variational operator and the Lagrangian is defined as
\begin{equation}
 L(\mathbf u , \psi, \varphi) =  \Pi_\text{kin} (\mathbf u) -\mathfrak{E} (\mathbf u ,\psi, \varphi) + \Pi_\text{ext} (\mathbf u , \varphi).
\end{equation}
Thus Equation~\eqref{eq:variation_l} expands as
\begin{equation}
 \delta  \int_{t_0}^{t_1}  \Pi_\text{kin} (\mathbf u) \, \mathrm d t - \delta  \int_{t_0}^{t_1}  \mathfrak{E} (\mathbf u , \psi,\varphi)  \, \mathrm d t  + \delta  \int_{t_0}^{t_1}  \Pi_\text{ext} (\mathbf u , \varphi)  \, \mathrm d t  = 0,
\label{eq:three_energy_time_integral}
\end{equation}
where the variation of the kinetic and external energy integrals can be expressed as
\begin{equation}
\delta  \int_{t_0}^{t_1}\Pi_\text{kin} (\mathbf u)  \mathrm d t 
%=  \delta  \int_{t_0}^{t_1} \left[ \frac{\rho h}{2} \int_\Gamma  \left[\frac{\partial u_}{\partial t}\right]^2 \,\mathrm d \Gamma \right]  \mathrm d t 
= -  \int_{t_0}^{t_1} \left[ {\rho h} \int_\Gamma   \delta u_i \frac{\partial^2 u_i}{\partial t^2} \,\mathrm d \Gamma \right]  \mathrm d t,
\end{equation}
and
\begin{align}
\delta  \int_{t_0}^{t_1}\Pi_\text{ext} (\mathbf u, \varphi)  \mathrm d t = \int_{t_0}^{t_1} \Bigg[ h\int_\Gamma  b_i \delta u_i \,\mathrm d \Gamma + h \int_{S_t} t_i \delta u_i \, \mathrm d S_t  + \frac{h^3}{6} \int_\Gamma q \delta \varphi \, \mathrm d \Gamma + \frac{h^3}{6} \int_{S_d} \omega \delta \varphi \,\mathrm d S_d \Bigg] \mathrm d t.
\end{align}
The variation of the electric enthalpy for the unelectroded shell is given by 
\begin{align}
\delta  \int_{t_0}^{t_1} \mathfrak{E}(\mathbf u, \psi, \varphi) \, \mathrm d t = &\int_{t_0}^{t_1} \int_{\Gamma} h\left[ \hat{C}^{abcd}\delta\alpha_{ab}\alpha_{cd} + \frac{h^2}{12} \hat{C}^{abcd}\delta\beta_{ab}\beta_{cd}\right] \, \mathrm d \Gamma  \mathrm d t \nonumber \\
&+ \int_{t_0}^{t_1} \left[ h \int_{\Gamma} \hat e^{3bc} \psi \delta \alpha_{bc} \, \mathrm d \Gamma + \frac{h^3}{6} \int_{\Gamma} \hat e^{3bc} \varphi \delta\beta_{bc} \, \mathrm d \Gamma  \right] \mathrm d t \nonumber \\
&+ \int_{t_0}^{t_1} \left[ h \int_{\Gamma} \hat e^{3bc} \delta \psi \alpha_{bc} \, \mathrm d \Gamma + \frac{h^3}{6} \int_{\Gamma} \hat e^{3bc} \delta\varphi \beta_{bc} \, \mathrm d \Gamma   \right] \mathrm d t \nonumber \\
&- \int_{t_0}^{t_1} \left[ \frac{h^3}{12} \int_{\Gamma} \hat \kappa^{ab} \delta\psi_{,a}  \psi_{,b} \, \mathrm d \Gamma + h \int_\Gamma \hat \kappa^{33}\delta \psi \, \psi \, \mathrm d \Gamma\right] \mathrm d t \nonumber \\
&- \int_{t_0}^{t_1} \left[ \frac{h^5}{30} \int_{\Gamma} \hat \kappa^{ab} \delta\varphi_{,a}  \varphi_{,b} \, \mathrm d \Gamma + \frac{h^3}{3} \int_\Gamma \hat \kappa^{33}\delta \varphi \, \varphi \, \mathrm d \Gamma\right] \mathrm d t,
\end{align}
% {\color{red}
% \begin{align}
% \delta  \int_{t_0}^{t_1} \mcal{H}(\mathbf u, \psi)  \mathrm d t = &\int_{t_0}^{t_1} \int_{\Gamma} h\left[ \hat{C}^{abcd}\delta\alpha_{ab}\alpha_{cd} + \frac{h^2}{12} \hat{C}^{abcd}\delta\beta_{ab}\beta_{cd}\right] \, \mathrm d \Gamma  \mathrm d t  + \int_{t_0}^{t_1} \left[ \frac{h^3}{12} \int_{\Gamma}\frac{e^{3ab}e^{3cd}}{\kappa^{33}} \delta\beta_{ab}\beta_{cd} \mathrm d \Gamma \right]\mathrm d t \nonumber \\
% &+ \int_{t_0}^{t_1} \left[ h \int_{\Gamma} \hat e^{3bc} \psi \delta \alpha_{bc} \, \mathrm d \Gamma  \right] \mathrm d t
% + \int_{t_0}^{t_1} \left[ h \int_{\Gamma} \hat e^{3bc} \delta \psi \alpha_{bc} \, \mathrm d \Gamma  \right] \mathrm d t \nonumber \\
% &- \int_{t_0}^{t_1} \left[ \frac{h^3}{12} \int_{\Gamma} \hat \kappa^{ab} \delta\psi_{,a}  \psi_{,b} \, \mathrm d \Gamma + h \int_\Gamma \hat \kappa^{33}\delta \psi \, \psi \, \mathrm d \Gamma\right] \mathrm d t,
% \end{align}
% }
for the symmetrically electroded shell by
\begin{align}
\delta  \int_{t_0}^{t_1} \mathfrak{E}(\mathbf u, \varphi) \, \mathrm d t = &\int_{t_0}^{t_1} \int_{\Gamma} h \left[ \hat{C}^{abcd}\delta\alpha_{ab}\alpha_{cd} + \frac{h^2}{12} \hat{C}^{abcd}\delta\beta_{ab}\beta_{cd}\right] \, \mathrm d \Gamma  \mathrm d t \nonumber \\
&+ \int_{t_0}^{t_1} \left[ \frac{h^3}{6} \int_{\Gamma} \hat e^{3bc} \varphi \delta\beta_{bc} - 2\bar{V}\int_{\Gamma} \hat e^{3bc} \delta \alpha_{bc}  \, \mathrm d \Gamma  \right] \mathrm d t 
+ \int_{t_0}^{t_1} \left[ \frac{h^3}{6} \int_{\Gamma} \hat e^{3bc} \delta\varphi \beta_{bc} \, \mathrm d \Gamma   \right] \mathrm d t \nonumber \\
&- \int_{t_0}^{t_1} \left[ \frac{h^5}{30} \int_{\Gamma} \hat \kappa^{ab} \delta\varphi_{,a}  \varphi_{,b} \, \mathrm d \Gamma + \frac{h^3}{3} \int_\Gamma \hat \kappa^{33}\delta \varphi \, \varphi \, \mathrm d \Gamma\right] \mathrm d t.
\end{align}
To satisfy Equation~\eqref{eq:three_energy_time_integral} for all possible $\delta\mathbf u$, $\delta \psi$ and $\delta\varphi$ (that vanish at the end of the time interval), the weak form of the governing equation for the unelectroded shell follows as
\begin{align}
& {\rho h} \int_\Gamma   \delta u_i \frac{\partial^2 u_i}{\partial t^2} \,\mathrm d \Gamma \nonumber \\
&+\int_{t_0}^{t_1} \int_{\Gamma} h \left[ \hat C^{abcd}\delta\alpha_{ab}\alpha_{cd} + \frac{h^2}{12} \hat C^{abcd}\delta\beta_{ab}\beta_{cd}\right] \, \mathrm d \Gamma  \mathrm d t \nonumber \\
&+ \int_{t_0}^{t_1} \left[ h \int_{\Gamma} \hat e^{3bc} \psi \delta \alpha_{bc} \, \mathrm d \Gamma + \frac{h^3}{6} \int_{\Gamma} \hat e^{3bc} \varphi \delta\beta_{bc} \, \mathrm d \Gamma  \right] \mathrm d t + \int_{t_0}^{t_1} \left[ h \int_{\Gamma} \hat e^{3bc} \delta \psi \alpha_{bc} \, \mathrm d \Gamma + \frac{h^3}{6} \int_{\Gamma} \hat e^{3bc} \delta\varphi \beta_{bc} \, \mathrm d \Gamma   \right] \mathrm d t \nonumber \\
&- \int_{t_0}^{t_1} \left[ \frac{h^3}{12} \int_{\Gamma} \hat \kappa^{ab} \delta\psi_{,a}  \psi_{,b} \, \mathrm d \Gamma + h \int_\Gamma \hat \kappa^{33}\delta \psi \, \psi \, \mathrm d \Gamma\right] \mathrm d t - \int_{t_0}^{t_1} \left[ \frac{h^5}{30} \int_{\Gamma} \hat \kappa^{ab} \delta\varphi_{,a}  \varphi_{,b} \, \mathrm d \Gamma + \frac{h^3}{3} \int_\Gamma \hat \kappa^{33}\delta \varphi \, \varphi \, \mathrm d \Gamma\right] \mathrm d t \nonumber \\
 &+\int_{t_0}^{t_1} \left[-h\int_\Gamma  b_i \delta u_i \,\mathrm d \Gamma - h \int_{S_t} \tau_i \delta u_i \, \mathrm d S_t  - \frac{h^3}{6} \int_\Gamma q \delta \varphi \, \mathrm d \Gamma - \frac{h^3}{6} \int_{S_d} \omega \delta \varphi \,\mathrm d S_d\right] \mathrm d t \nonumber \\
 &= \,0,
\label{eq:weak_form_1}
\end{align}
% {\color{red}
% \begin{align}
% & {\rho h} \int_\Gamma   \delta u_i \frac{\partial^2 u_i}{\partial t^2} \,\mathrm d +\int_{t_0}^{t_1} \int_{\Gamma} h\left[ \hat{C}^{abcd}\delta\alpha_{ab}\alpha_{cd} + \frac{h^2}{12} \hat{C}^{abcd}\delta\beta_{ab}\beta_{cd}\right] \, \mathrm d \Gamma  \mathrm d t  + \int_{t_0}^{t_1} \left[ \frac{h^3}{12} \int_{\Gamma}\frac{e^{3ab}e^{3cd}}{\kappa^{33}} \delta\beta_{ab}\beta_{cd} \mathrm d \Gamma \right]\mathrm d t \nonumber \\
% &+ \int_{t_0}^{t_1} \left[ h \int_{\Gamma} \hat e^{3bc} \psi \delta \alpha_{bc} \, \mathrm d \Gamma  \right] \mathrm d t
% + \int_{t_0}^{t_1} \left[ h \int_{\Gamma} \hat e^{3bc} \delta \psi \alpha_{bc} \, \mathrm d \Gamma  \right] \mathrm d t \nonumber \\
% &- \int_{t_0}^{t_1} \left[ \frac{h^3}{12} \int_{\Gamma} \hat \kappa^{ab} \delta\psi_{,a}  \psi_{,b} \, \mathrm d \Gamma + h \int_\Gamma \hat \kappa^{33}\delta \psi \, \psi \, \mathrm d \Gamma\right] \mathrm d t + \int_{t_0}^{t_1} \left[-h\int_\Gamma  b_i \delta u_i \,\mathrm d \Gamma - h \int_{S_t} \tau_i \delta u_i \, \mathrm d S_t \right] \mathrm d t \nonumber \\
%  &= \,0,
% \label{eq:weak_form_1}
% \end{align}
% }
and for the symmetrically electroded shell as
\begin{align}
& {\rho h} \int_\Gamma   \delta u_i \frac{\partial^2 u_i}{\partial t^2} \,\mathrm d \Gamma \nonumber \\
& +\int_{t_0}^{t_1} \int_{\Gamma} h\left[ \hat C^{abcd}\delta\alpha_{ab}\alpha_{cd} + \frac{h^2}{12} \hat C^{abcd}\delta\beta_{ab}\beta_{cd}\right] \, \mathrm d \Gamma  \mathrm d t \nonumber \\
&+ \int_{t_0}^{t_1} \left[ \frac{h^3}{6} \int_{\Gamma} \hat e^{3bc} \varphi \delta\beta_{bc} - 2\bar{V}\int_{\Gamma} \hat e^{3bc} \delta \alpha_{bc}  \, \mathrm d \Gamma  \right] \mathrm d t 
+ \int_{t_0}^{t_1} \left[ \frac{h^3}{6} \int_{\Gamma} \hat e^{3bc} \delta\varphi \beta_{bc} \, \mathrm d \Gamma   \right] \mathrm d t \nonumber \\
&- \int_{t_0}^{t_1} \left[ \frac{h^5}{30} \int_{\Gamma} \hat \kappa^{ab} \delta\varphi_{,a}  \varphi_{,b} \, \mathrm d \Gamma + \frac{h^3}{3} \int_\Gamma \hat \kappa^{33}\delta \varphi \, \varphi \, \mathrm d \Gamma\right] \mathrm d t
 \nonumber \\
 &+\int_{t_0}^{t_1} \left[-h\int_\Gamma  b_i \delta u_i \,\mathrm d \Gamma - h \int_{S_t} \tau_i \delta u_i \, \mathrm d S_t  - \frac{h^3}{6} \int_\Gamma q \delta \varphi \, \mathrm d \Gamma - \frac{h^3}{6} \int_{S_d} \omega \delta \varphi \,\mathrm d S_d\right] \mathrm d t \nonumber \\
 &= \,0.
\label{eq:weak_form_2}
\end{align}

\subsection{Discretisation and system of equations}
\label{sec:system_equation}
The displacement is discretised using the subdivision surface basis functions as
\begin{equation}
\mathbf u = \sum_{A=0}^{n_{b}-1} N^A \mathbf{U}_A,
\end{equation}
where $n_b$ is the number of basis functions, and $\mathbf U_A$ denotes the $A^{\text{th}}$ nodal coefficients of displacement. Thus the membrane and bending strain components are computed as
\begin{align}
\alpha_{a b} = &  \sum_{A=0}^{n_{b}-1} \frac{1}{2} [ N^A_{,b}  \bar{\mathbf a}_{a} + N^A_{,a} \bar{\mathbf a}_{b}] \cdot\mathbf U_A, \\
\beta_{a b} = &  \sum_{A=0}^{n_{b}-1} \Big[ -N^A_{,{a b}} \bar{\mathbf a}_3 + \frac{1}{\bar{J}} \big[N^A_{,1} [\bar{\mathbf a}_{a , b} \times \bar{\mathbf a}_2 ] + N^A_{,2} [\bar{\mathbf a}_{1} \times \bar{\mathbf a}_{a , b} ] \big] + \frac{\bar{\mathbf a}_3 \cdot \bar{\mathbf a}_{a , b}}{\bar J}\big[N^A_{,1} [\bar{\mathbf a}_{2} \times \bar{\mathbf a}_3 ] + N^A_{,2}  [\bar{\mathbf a}_{3} \times \bar{\mathbf a}_{1} ] \big] \Big] \cdot \mathbf{U}_A.
\end{align}
The electrical potential functions are also discretised using the same basis functions as $\mathbf u$, and expressed as
% \begin{equation}
% \varphi = \sum_{A=0}^{n_b -1} N^A \Phi_A,
% \end{equation}
\begin{align}
\psi = \sum_{A=0}^{n_b -1} N^A \Psi_A, \quad \varphi = \sum_{A=0}^{n_b -1} N^A \Phi_A.
\end{align}
Here $\Psi_A$ and $\Phi_A$ are the $A^{\text{th}}$ nodal coefficients of the potential functions. Following a Bubnov-Galerkin approach, the subdivision surface bases are also used for the trial functions $\delta \mathbf u$ and $\delta \varphi$, and the weak form~\eqref{eq:weak_form_1} follows in matrix format as
\begin{equation}
\begin{bmatrix} 
\ary M &  \ary 0 & \ary 0 \\
 \ary 0 &  \ary 0 & \ary 0 \\
  \ary 0 & \ary 0 & \ary 0
\end{bmatrix}
\begin{bmatrix}
\ddot{\ary u} \\
\ary 0 \\
\ary 0
\end{bmatrix} + 
\begin{bmatrix} 
\ary K & \ary {C}_{u \psi} & \ary{C}_{u \varphi} \\
\ary {C}_{\psi u} & \ary {D}_1 & \ary 0 \\
\ary {C}_{\varphi u} & \ary 0  & \ary {D}_2\\
\end{bmatrix} 
\begin{bmatrix}
\ary u \\
\ary \psi \\
\ary \varphi
\end{bmatrix}=
\begin{bmatrix}
\ary {f}_{u}\\
\ary 0 \\
\ary {f}_{\varphi}
\end{bmatrix},
\label{eq:system_of_equations}
\end{equation}

% {\color{red}
% \begin{equation}
% \begin{bmatrix} 
% \ary M &  \ary 0 \\
%  \ary 0 &  \ary 0
% \end{bmatrix}
% \begin{bmatrix}
% \ddot{\ary u} \\
% \ary 0
% \end{bmatrix} +  
% \begin{bmatrix} 
% \ary K + \ary B & \ary {C_{u \psi}} \\
% \ary {C_{\psi u}} & \ary {D_1}
% \end{bmatrix} 
% \begin{bmatrix}
% \ary u \\
% \ary \psi
% \end{bmatrix}=
% \begin{bmatrix}
% \ary {f_{u}}\\
% \ary {0}
% \end{bmatrix},
% \label{eq:system_of_equations}
% \end{equation}
% }
and Equation~\eqref{eq:weak_form_2} follows in matrix format as
\begin{equation}
\begin{bmatrix} 
\ary M &  \ary 0 \\
 \ary 0 &  \ary 0
\end{bmatrix}
\begin{bmatrix}
\ddot{\ary u} \\
\ary 0
\end{bmatrix} +  
\begin{bmatrix} 
\ary K + \ary P & \ary {C}_{u \varphi} \\
\ary {C}_{\varphi u} & \ary {D}_2
\end{bmatrix} 
\begin{bmatrix}
\ary u \\
\ary \varphi
\end{bmatrix}=
\begin{bmatrix}
\ary {f}_{u}\\
\ary {f}_{\varphi}
\end{bmatrix}.
\label{eq:system_of_equations}
\end{equation}
Here $\ary M$ is the global mass matrix. $\ddot{\ary u}$ is the global acceleration vector. $\ary K$ denotes the global stiffness matrix, $\ary {D}_1$ and $\ary {D}_2$ are the global dielectric system matrices, $\ary{C}_{u\psi}$($\ary{C}_{u\varphi}$) and  $\ary{C}_{\psi u}$($\ary{C}_{\varphi u}$) are the direct and converse piezoelectric coupling matrices, respectively. { The global matrix $\ary P$ has only diagonal entries and takes into account the direct piezoelectric effects caused by the prescribed voltage. } {$\ary u$, $\ary \psi$ and $\ary \varphi$} are the global vectors of displacement, and the first and second order electrical potential coefficients, respectively. $\ary{f}_u$ and $\ary{f}_\varphi$ on the right hand side denote the global structural and electrical load vectors. Note, the system of equations is non-symmetric. For computational efficiency, we modify the system of equations using the Schur complements $\ary{C}_{u \psi} \ary{D}_1^{-1}\ary{C}_{\psi u}$ and $\ary{C}_{u \varphi} \ary{D}_2^{-1}\ary{C}_{\varphi u}$.
%global electrical potential related function using the inverse of $\ary D$ as 
%\begin{equation}
%\ary \Phi = \ary{D}^{-1} [\ary f_\varphi - \ary{C}_{\varphi u} \ary U],
%\end{equation}
%and substitute this expression in the mechanical part of the system of equations~\eqref{eq:system_of_equations}, 
Thus the problem for $\ary u$ becomes
% \begin{equation}
% \ary M  \ddot{\ary u}+[{\ary K - \ary{C}_{u \varphi} \ary{D}^{-1}\ary{C}_{\varphi u}}] \ary u = \ary{f}_u - \ary{C}_{u \varphi} \ary{D}^{-1} \ary{f}_\varphi,
% \end{equation}
\begin{equation}
\ary M  \ddot{\ary u}+[{\ary K - \ary{C}_{u \psi} \ary{D}_1^{-1}\ary{C}_{\psi u} - \ary{C}_{u \varphi} \ary{D}_2^{-1}\ary{C}_{\varphi u}}] \ary u = \ary{f}_u- \ary{C}_{u \varphi} \ary{D}_2^{-1} \ary{f}_\varphi
\end{equation}
% {\color{red}
% \begin{equation}
% \ary M  \ddot{\ary u}+[{\ary K + \ary B - \ary{C_{u \psi}} \ary{D_1}^{-1}\ary{C_{\psi u}}}] \ary u = \ary{f_u}
% \end{equation}
% }
for the unelectroded case and
\begin{equation}
\ary M  \ddot{\ary u}+[{\ary K + \ary P - \ary{C}_{u \varphi} \ary{D}_2^{-1}\ary{C}_{\varphi u}}] \ary u = \ary{f}_u - \ary{C}_{u \varphi} \ary{D}_2^{-1} \ary{f}_\varphi
\end{equation}
for shells with symmetrically prescribed voltage electrodes. Consequently, one defines new global system matrices
\begin{equation}
\ary A = \ary K - \ary{C}_{u \psi} \ary{D}_1^{-1}\ary{C}_{\psi u} - \ary{C}_{u \varphi}\ary{D}_2^{-1}\ary{C}_{\varphi u},
\end{equation}
or
\begin{equation}
\ary A = \ary K + \ary P - \ary{C}_{u \varphi} \ary{D}_2^{-1}\ary{C}_{\varphi u},
\end{equation}
respectively.
The system of equations is thus defined by
\begin{equation}
\ary M  \ddot{\ary u}+\ary A \ary u = \ary{f}_u - \ary{C}_{u \varphi} \ary{D}_2^{-1} \ary{f}_\varphi.
\end{equation}
The problem of a free vibrating piezoelectric shell can be obtained by assuming harmonic motions, and is given by
\begin{equation}
[-\omega^2 \ary M + \ary A] \ary u = \ary{f}_u - \ary{C}_{u \varphi} \ary{D}_2^{-1} \ary{f}_\varphi,
\end{equation}
where $\omega$ is the angular frequency.
For the free vibration analysis, the external mechanical and electrical loads are set to zero, and the system of equation reduces to
\begin{equation}
-\omega^2 \ary M +\ary A = \ary 0.
\label{eq:eigen_system_equation}
\end{equation}
\end{nolinenumbers}
\section{Numerical examples}
\label{sec:numerical}
Three numerical examples are considered. This first is the free vibration of an elastic spherical thin shell which is used to validate the Kirchhoff-Love shell formulation and implementation. Then, the piezoelectric effect for curved shells is investigated using the  \textit{Scordelis-Lo} roof geometry. The final example demonstrates the potential of the formulation by analysing the vibration of piezoelectric shell applications with complex geometry. All numerical results are computed using the open source finite element library deal.II~\cite{BangerthHartmannKanschat2007,ArndtBangerthDavydov-2021-a}.
\subsection{Validation using a elastic spherical shell}
\begin{figure}[!t]
\centering
\begin{subfigure}[b]{0.49\linewidth}
	\centering
 	 \includegraphics[width=0.6\linewidth]{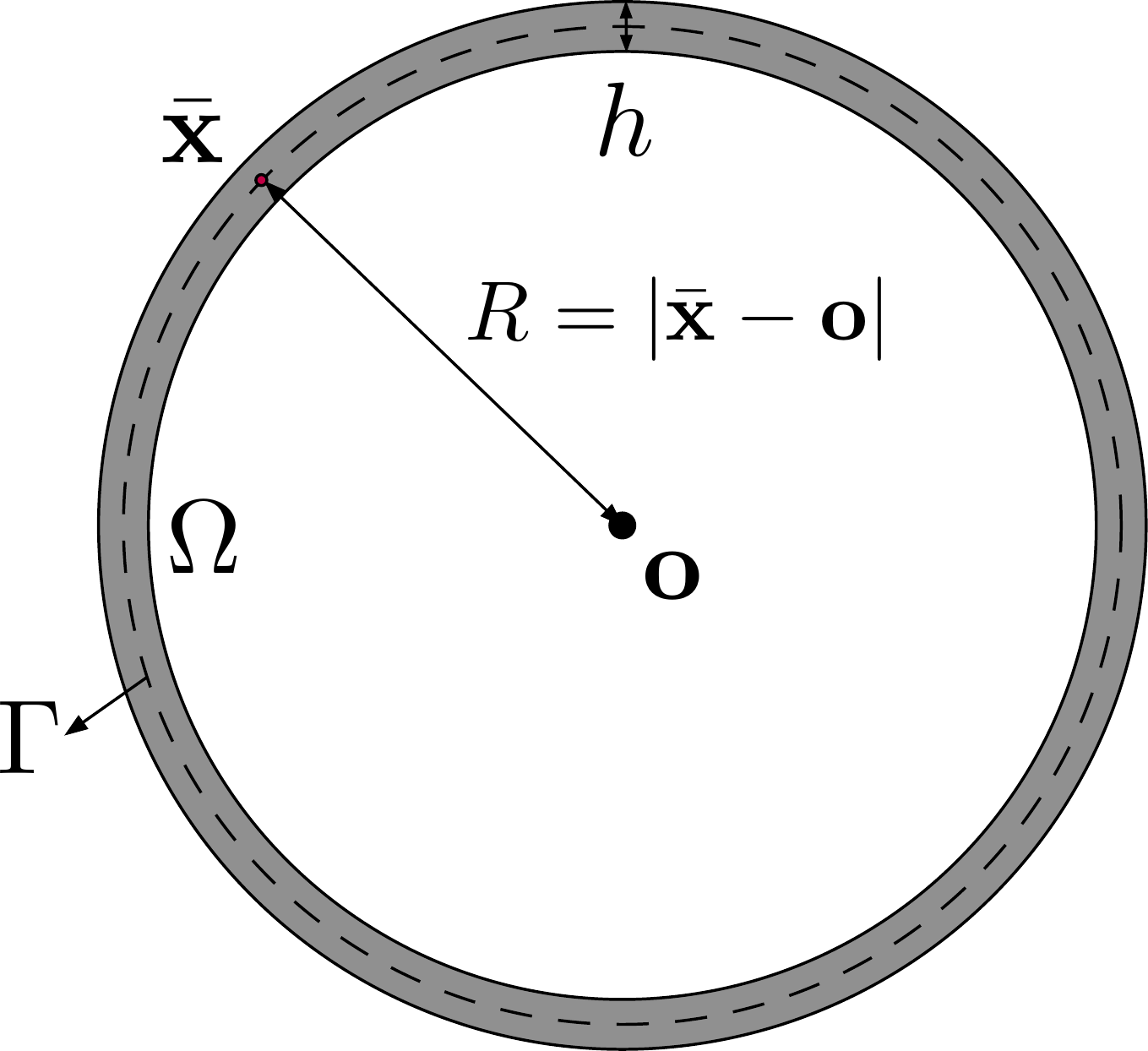}
	\caption{A slice of the spherical shell}
	\label{fig:spherical_shell_definition}
\end{subfigure}
\begin{subfigure}[b]{0.49\linewidth}
	\centering
 	 \includegraphics[width=0.6\linewidth]{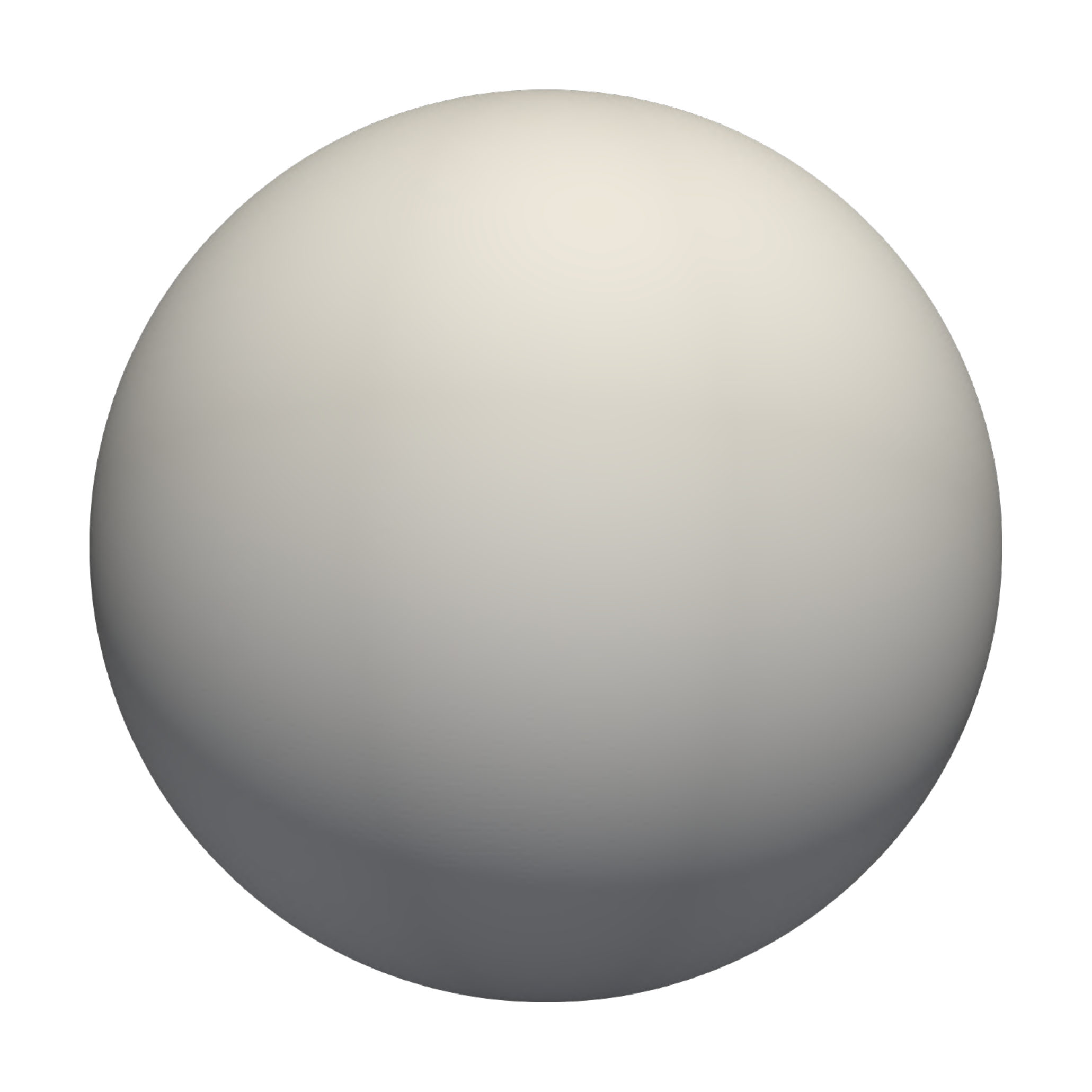}
	\caption{Limit surface}
	\label{fig:sphere_limit}
\end{subfigure}
\begin{subfigure}[b]{0.49\linewidth}
	\centering
  	\includegraphics[width=0.6\linewidth]{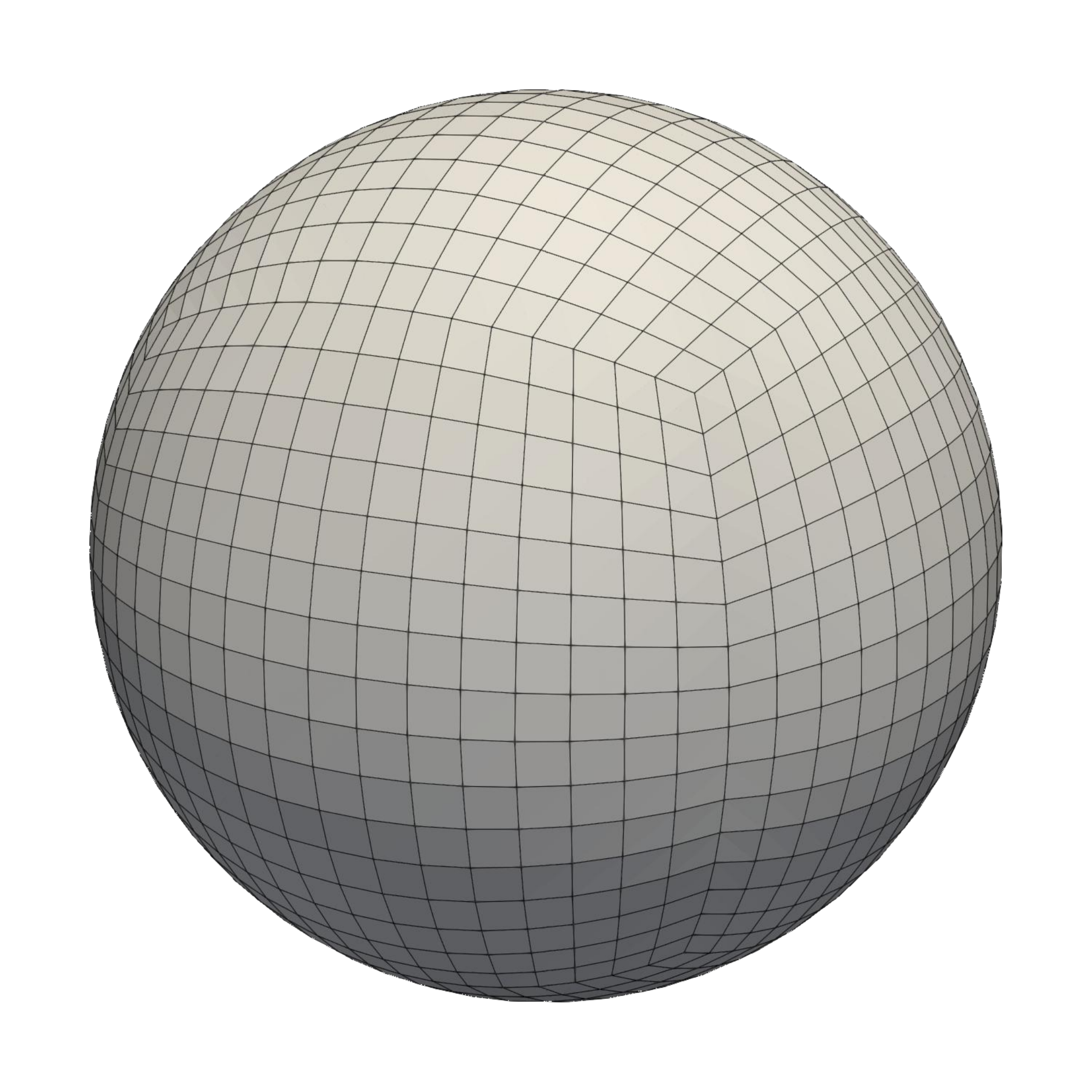}
	\caption{Initial control mesh}
	\label{fig:sphere_mesh}
\end{subfigure}
\begin{subfigure}[b]{0.49\linewidth}
	\centering
  	\includegraphics[width=0.6\linewidth]{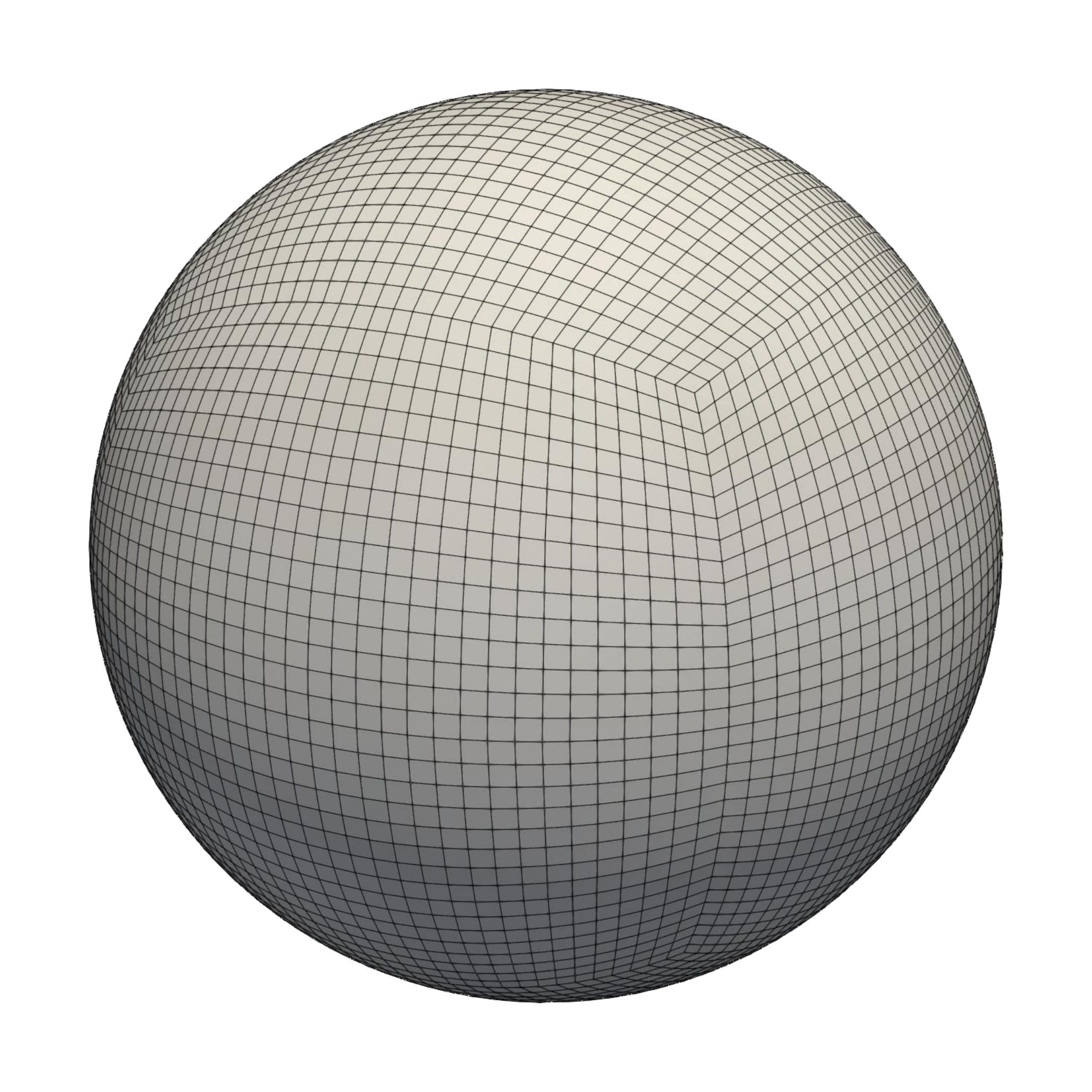}
	\caption{First refinement}
	\label{fig:sphere_mesh_2}
\end{subfigure}
\begin{subfigure}[b]{0.55\linewidth}
	\centering
  	\includegraphics[width=0.6\linewidth]{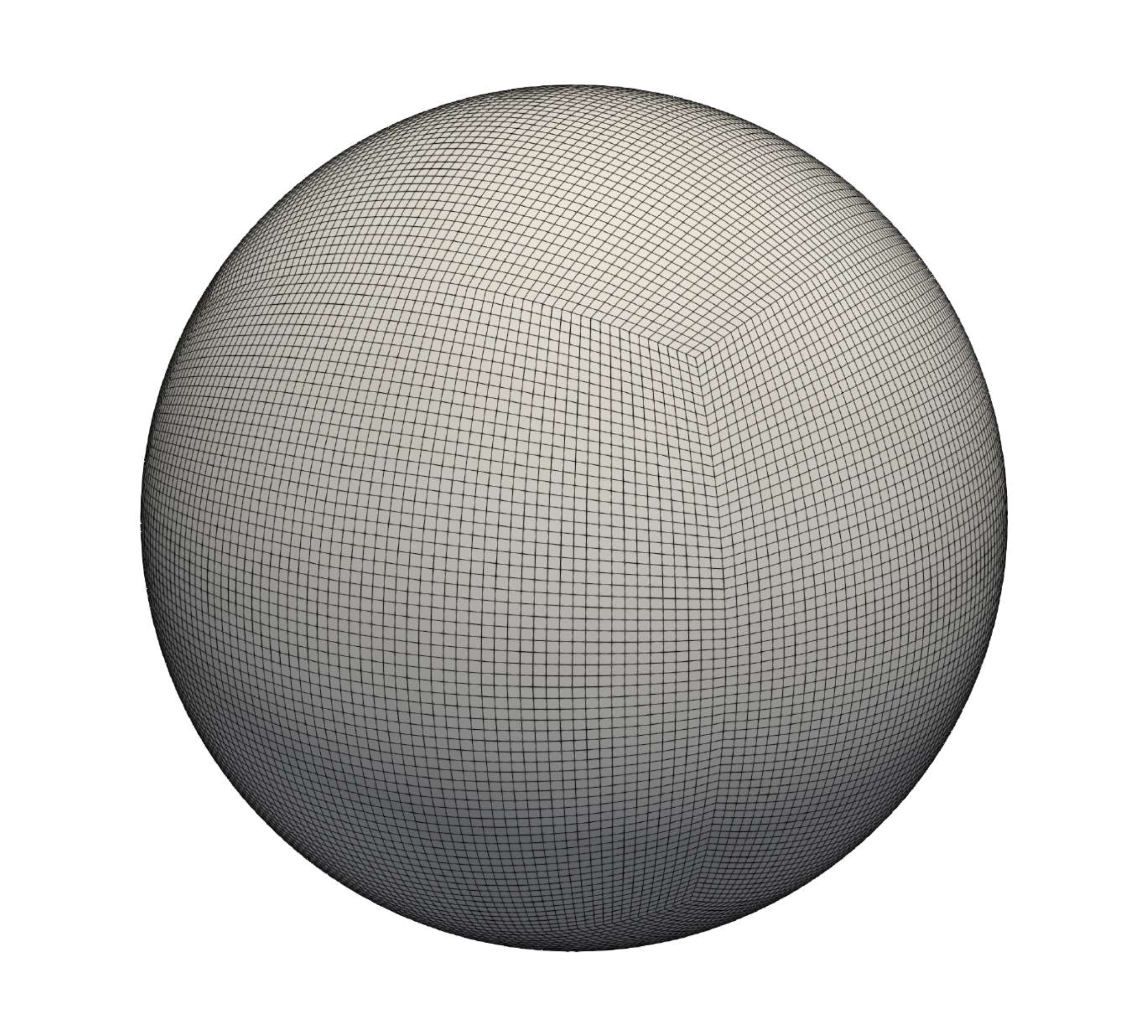}
	\caption{Second refinement}
	\label{fig:sphere_mesh_2}
\end{subfigure}
	\caption{(a) Definition of a spherical shell.  (b) The limiting surface constructed using Catmull-Clark subdivision from (c). (c) A control mesh with $1536$ elements for the mid-surface of the spherical shell. (d) $1^{\text{st}}$ level refined mesh with $6144$ elements using a least square fitting method. (e) $2^{\text{nd}}$ level refined mesh with $24576$ elements}
	\label{fig:sphere}
\end{figure}

The first numerical example is the free vibration analysis of an elastic spherical thin shell which is used to validate the pure elastic Kirchhoff-Love shell formulation. This problem was first examined by~\citet{lamb1882vibrations}. \citet{baker1961axisymmetric} used the membrane theory to examine the axisymmetric modes of a complete spherical shell. The method developed here is based on thin shell elements and can compute both axisymmetric and nonaxisymmetric modes. Figure~\ref{fig:spherical_shell_definition} shows a cross section of the spherical shell domain $\Omega$, which has a uniform thickness $h$ and with $\Gamma$ denoting its mid-surface. The radius $R$ measures the distance between the center of the sphere to the mid-surface. 

If the material is assumed as isotropic, the elastic strain energy density per unit area consists of the membrane and bending parts~\cite{cirakortiz2000} as
\begin{equation}
\widetilde{W}_\text{ela}( \mathbf A, \mathbf B) = \frac{1}{2} \frac{Eh}{1-\nu^2}\left[ [\mathbf A : \mathbf H : \mathbf A] + \frac{h^2}{12} [\mathbf B : \mathbf H : \mathbf B] \right]= \frac{1}{2} \frac{Eh}{1-\nu^2} H^{abcd}\alpha_{ab}\alpha_{cd} + \frac{1}{2} \frac{Eh^3}{12[1-\nu^2]} H^{abcd}\beta_{ab}\beta_{cd} ,
\end{equation}
where $E$ and $\nu$ are the Young's modulus and Poisson's ratio, respectively. %Here $i,j,k,l$ vary from $1$ to $2$. 
$H^{abcd}$ denote the components of the fourth-order tensor $\mathbf H$ computed from the contravariant metric tensors as
\begin{equation}
H^{abcd} = \nu\bar{a}^{ab}\bar{a}^{cd} + \frac{1}{2} [1-\nu][\bar{a}^{ac}\bar{a}^{bd} + \bar{a}^{ad}\bar{a}^{bc}].
\end{equation}

\begin{figure}[!t]
\centering
  \includegraphics[width=\linewidth]{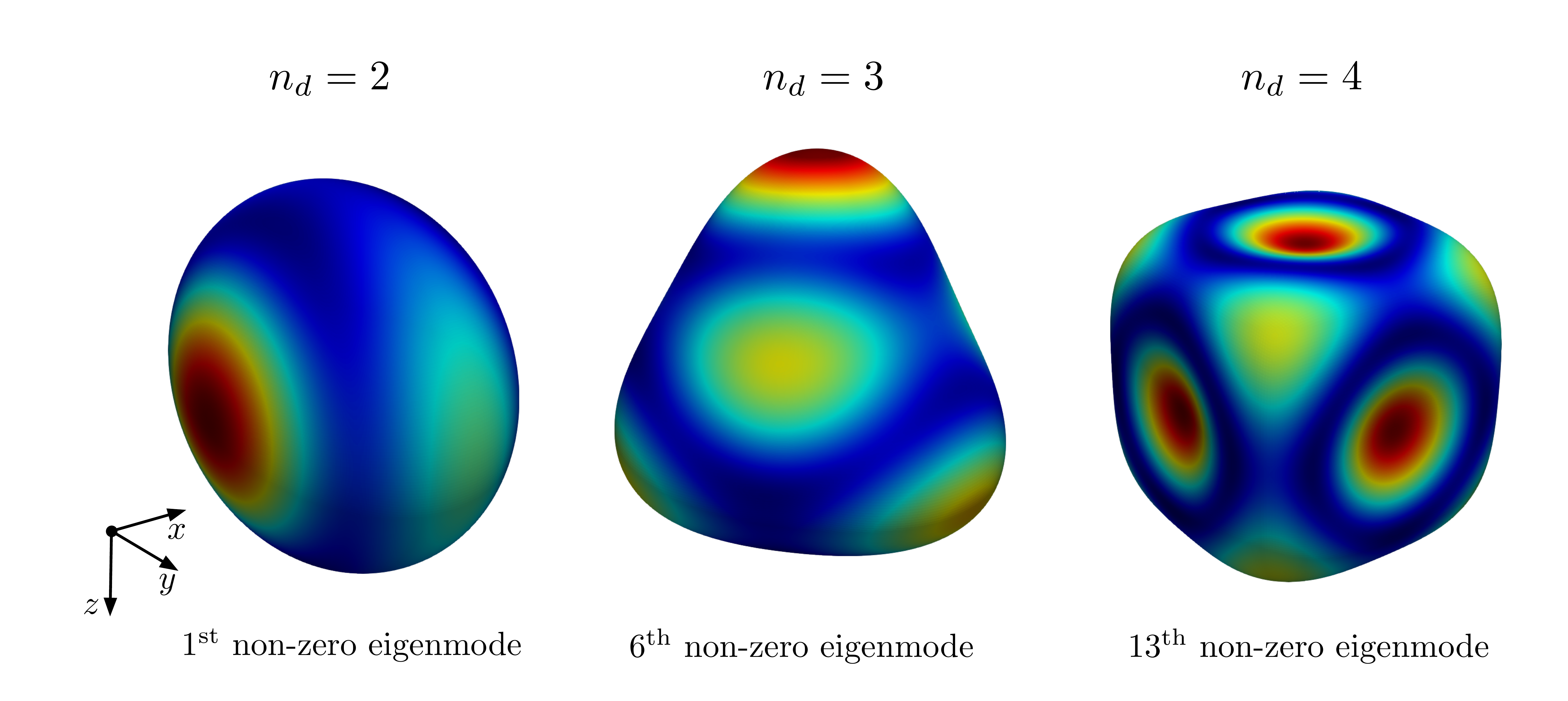}
\caption{Examples of the mode shapes. The color represents the magnitude of displacement $|\mathbf u|$.}
\label{fig:sphere_modes}
\end{figure}

\citet{duffey2007vibrations} provide a comparison of experimental results~\cite{robertson2004modal} with analytical solutions for the problem considered here. The values of the geometric and material parameters are given in Table~\ref{tab:sphere_parameters}. It is worth noting here that they used the imperial system of units in their work. 
Here we aim to simulate the same problem using the proposed method and compare our numerical results to experimental and analytical solutions. 
Since no piezoelectric effect is considered in this problem, the system of equations~\eqref{eq:eigen_system_equation} simplifies to
\begin{equation}
-\omega^2 \ary M + \ary K = \ary 0.
\end{equation}
This system of equations can be solved as an eigenvalue problem where $\omega^2$ is the eigenvalue and the eigenvectors can be used to generate the corresponding eigenmode shapes. The natural frequency is computed as
\begin{equation}
f = \frac{\omega}{2\pi}.
\end{equation}
The vibration modes of the spherical shell can be defined in terms of a polynomial degree $n_d$, where $n_d = 1,2,3,\cdots$. Each polynomial degree corresponds to a $2 n_d + 1$ clustering of eigenvalues with different eigenmodes. $n_d = 1$ corresponds to a rigid body motion and the corresponding eigenvalue equals to $0$. Thus the first non-zero eigenvalue corresponds $n_d = 2$.
Figure~\ref{fig:sphere_mesh} shows the control mesh used to construct the Catmull-Clark subdivision limit surface (Figure~\ref{fig:sphere_limit}) for the mid-surface of a spherical thin shell. The control mesh contains $1536$ elements with $8$ extraordinary vertices. The presence of extraordinary vertices leads to computational errors which can be reduced using an adaptive quadrature scheme~\cite{juttler2016numerical,liu2020assessment}. Two refined meshes with $6144$ and $24576$ elements generated using a least square fitting method are also used for this problem. Table~\ref{tab:sphere_results} shows the numerical results for both the initial and refined control meshes. For $n_d = 2$, the numerically determined natural frequency has only a small error of approximately $0.296\%$ for the initial mesh, $0.087\%$ for the first level refinement and $0.024\%$ for the second level refinement. The numerical error increases as the mode becomes more complex.  For $n_d = 3$, the error is in the range of $(0.180\%, 0.528\%)$ for the initial mesh and $(0.062\%, 0.149\%)$ for the first level refinement and (0.021\%, 0.043\%) for the second level refinement. For $n_d = 4$ the errors are in the range of $(0.159\%, 0.510\%)$ for the initial mesh, $(0.059\%, 0.147\%)$ for the first level refinement, and $(0.020\%, 0.044\%)$ for the second level refinement. The results show clear convergence to the analytical solutions and the deviation of each $n_d$ is reduced after refinement. Figure~\ref{fig:sphere_modes} shows the vibration modes for the $1^{\text{st}}$, $6^{\text{th}}$ and $13^{\text{th}}$ non-zero eigenvalues which corresponding to $n_d = 2, 3 \text{ and } 4$, respectively.
\begin{table}[]
\centering
\caption{Geometry and material parameters for the elastic spherical thin shell.}
\begin{tabular}{@{}lll@{}}
\toprule
Parameter       & Value    \\ \midrule
Radius $R$  & $4.4688$ (in) & $0.1135$ (m)\\
Thickness $h$  & $0.0625$ (in) & $1.5875$ (mm) \\
Young's modulus $E$ & $28\times 10^6$ (psi) & $193.05$ (GPa)\\
Poisson's ratio $\nu$ & $0.28$  & \\
Mass density  $\rho$  & $0.000751$ (lbf-$\text{s}^2$/$\text{in}^4$) & 8025.937(kg/$\text{m}^{3}$) \\
 \bottomrule
\end{tabular}
\label{tab:sphere_parameters}
\end{table}

\begin{table}[]
\centering
\caption{Comparison of numerical results with analytical solutions and experimental results~\cite{robertson2004modal,duffey2007vibrations}}
\begin{tabular}{@{}CCCCCCC@{}}
\toprule
\multirow{3}{*}{$n_d$} & \multirow{3}{*}{Mean experimental } & \multirow{3}{*}{Analytical solutions} & \multicolumn{3}{c}{Numerical results}  \\ \cmidrule(l){4-7} 
                             &                                      &                                      & $Non-zero$  & \multicolumn{2}{c}{$f$ (Hz)}     \\
                                                          & f_e$(Hz)$                                     &    f_a$(Hz)$                                  & $eigenvalue number$         &  $Initial mesh$     & 1^{\text{st}}$Refinement$ & 2^{\text{nd}}$Refinement$\\ \midrule
\multirow{5}{*}{$2$}           & \multirow{5}{*}{$5088$}                & \multirow{5}{*}{$5078$}                & 1                    & 5092.80    & 5082.37 & 5079.23\\
                             &                                      &                                      & 2                    & 5092.80      &  5082.37 & 5079.23\\
                             &                                      &                                      & 3                    & 5093.05       & 5082.41 & 5079.24\\
                             &                                      &                                      & 4                    & 5093.05       &  5082.41 & 5079.24\\
                             &                                      &                                      & 5                    & 5093.05       & 5082.41 & 5079.24\\ \midrule
\multirow{7}{*}{$3$}           & \multirow{7}{*}{$6028$}                & \multirow{7}{*}{$6005$}                & 6                    & 6015.79      &  6008.71 & 6006.26\\
                             &                                      &                                      & 7                    & 6015.79       & 6008.71 & 6006.26\\
                             &                                      &                                      & 8                    & 6015.79       & 6008.71 & 6006.26\\
                             &                                      &                                      & 9                    & 6025.04       & 6010.95 & 6006.81\\
                             &                                      &                                      & 10                   & 6025.04      & 6010.95 & 6006.81\\
                             &                                      &                                      & 11                   & 6025.04       & 6010.95 & 6006.81\\
                             &                                      &                                      & 12                   & 6036.72       & 6013.93 & 6007.57\\ \midrule
\multirow{9}{*}{$4$}           & \multirow{9}{*}{$6379$}                & \multirow{9}{*}{$6378$}                & 13                   & 6388.13      &  6381.75 & 6379.38\\
                             &                                      &                                      & 14                   & 6389.03       & 6381.94 & 6379.42\\
                             &                                      &                                      & 15                   & 6389.03       & 6381.94 & 6379.42\\
                             &                                      &                                      & 16                   & 6389.03       & 6381.94 & 6379.42\\
                             &                                      &                                      & 17                   & 6392.17       & 6382.62 & 6379.59\\
                             &                                      &                                      & 18                   & 6392.17       & 6382.62 & 6379.59\\
                             &                                      &                                      & 19                   & 6410.58       & 6387.35 & 6380.78\\
                             &                                      &                                      & 20                   & 6410.58       & 6387.35 & 6380.78\\
                             &                                      &                                      & 21                   & 6410.58       & 6387.35 & 6380.78\\ 
\bottomrule
\end{tabular}
\label{tab:sphere_results}
\end{table}

\subsection{Piezoelectric effects on the vibration of a \textit{Scordelis-Lo} roof }
\begin{figure}[!t]
\centering
  \includegraphics[width=0.6\linewidth]{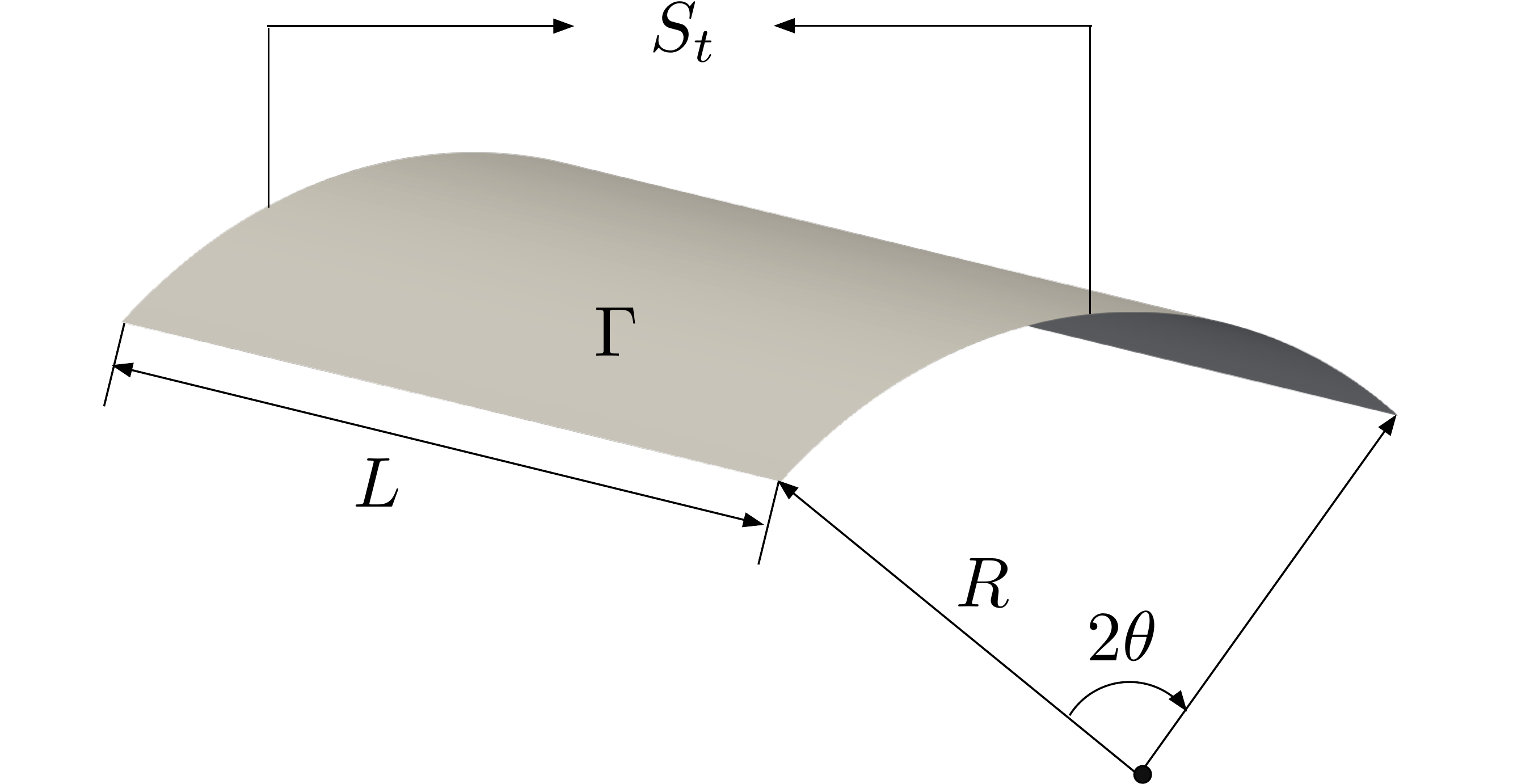}
\caption{ \textit{Scordelis-Lo} roof geometry. Simply supported boundary condition are applied on $S_t$.}
\label{fig:roof_definition}
\end{figure}

The following numerical example is a \textit{Scordelis-Lo} roof, which is commonly used as a benchmark problem for shell formulations. The \textit{Scordelis-Lo} roof is a simple geometry which only requires a structured quadrilateral mesh without extraordinary vertices. It can be considered as a plate curved in one direction. Figure~\ref{fig:roof_definition} shows the geometry of the roof which can be defined using a length $L$, a radius $R$ and an angular parameter $\theta$. 
%This geometry is used to investigate the piezoelectric effect on curved shell. The choice of the geometry and material parameters are given in Table~\ref{tab:roof_parameters}.
We note here that for the well-known benchmark problem~\cite{belytschko1985stress}, the units of parameters are omitted. The geometry parameters are set to $L = 50$, $R = 25$ and $\theta = 40^\circ$.  The two curved edges of the roof are simply supported. The roof has a thickness $h=0.25$ and a self weight of $90$ is applied as a uniformed load in negative $z$ direction. The Young's modulus $E$ for the benchmark problem is $4.32\times 10^8$ and Poisson's ratio $\nu = 0$. The reference solution of the \textit{Scordeli-Lo} roof shell is given by the mid-point vertical displacement $u_z$ of the two free edges and is equal to $0.3024$. Our results converge to $0.3006$. Such a minor difference is also observed in other IGA shell literature~\cite{kiendl2009isogeometric}. 

\begin{figure}[!t]
\centering
\begin{subfigure}[b]{0.15\linewidth}
	\centering
 	 \includegraphics[width=\linewidth]{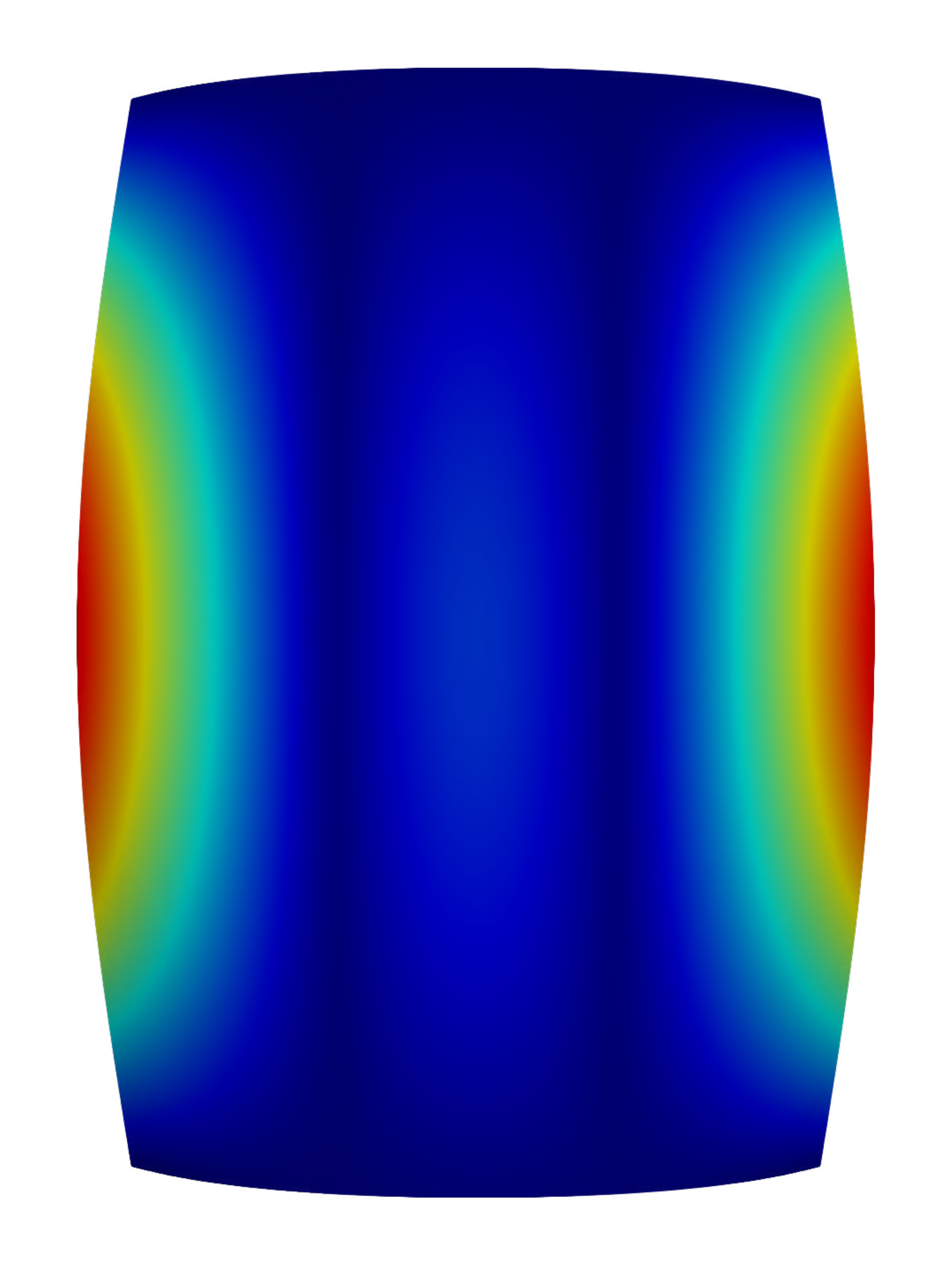}
	\caption{Mode 1:$|\mathbf u|$}
	\label{fig:roof_coupled_disp_1}
\end{subfigure}
\begin{subfigure}[b]{0.15\linewidth}
	\centering
  	\includegraphics[width=\linewidth]{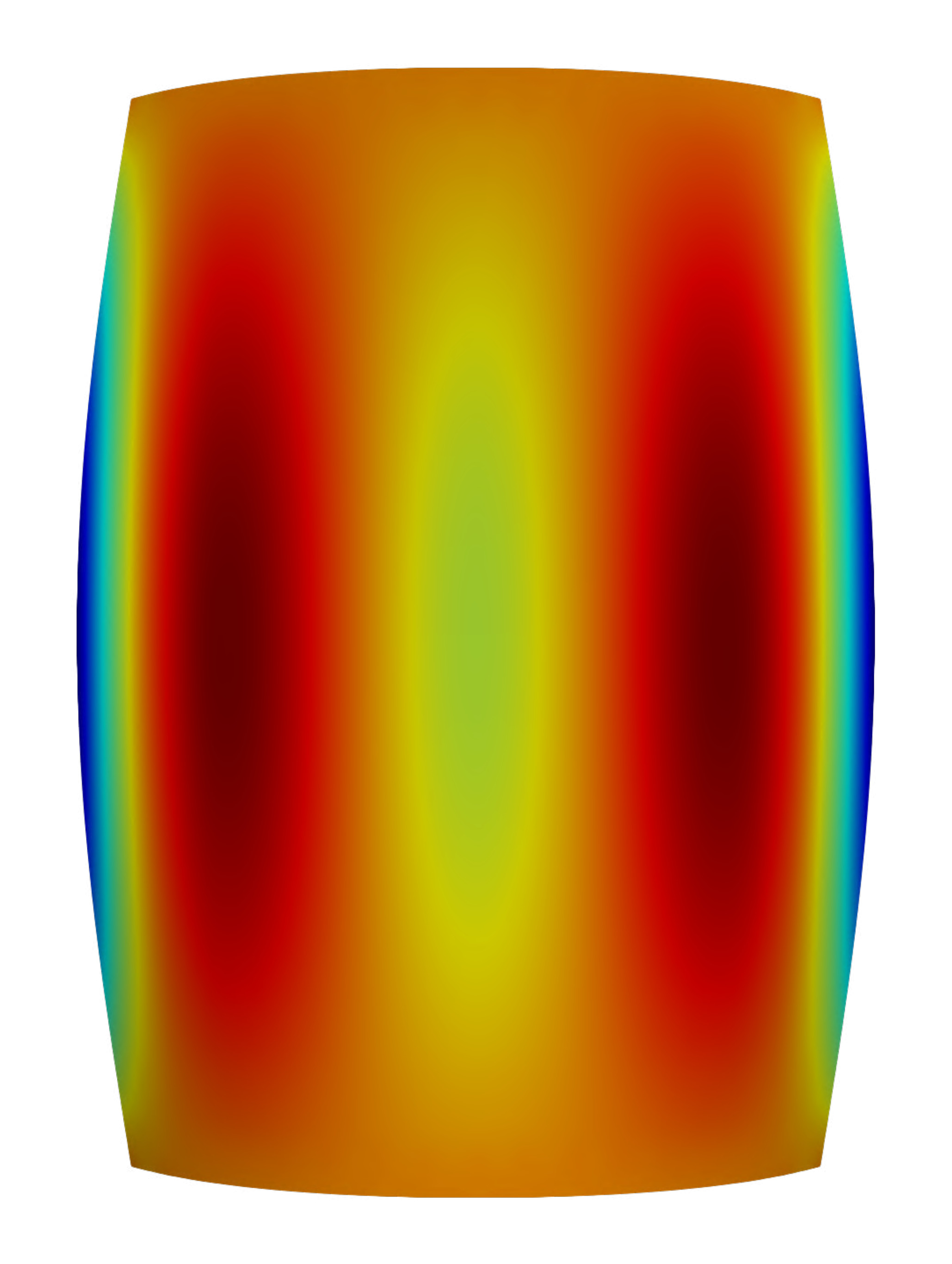}
	\caption{Mode 1: $\psi$}
	\label{fig:roof_coupled_psi_1}
\end{subfigure}
\begin{subfigure}[b]{0.15\linewidth}
	\centering
  	\includegraphics[width=\linewidth]{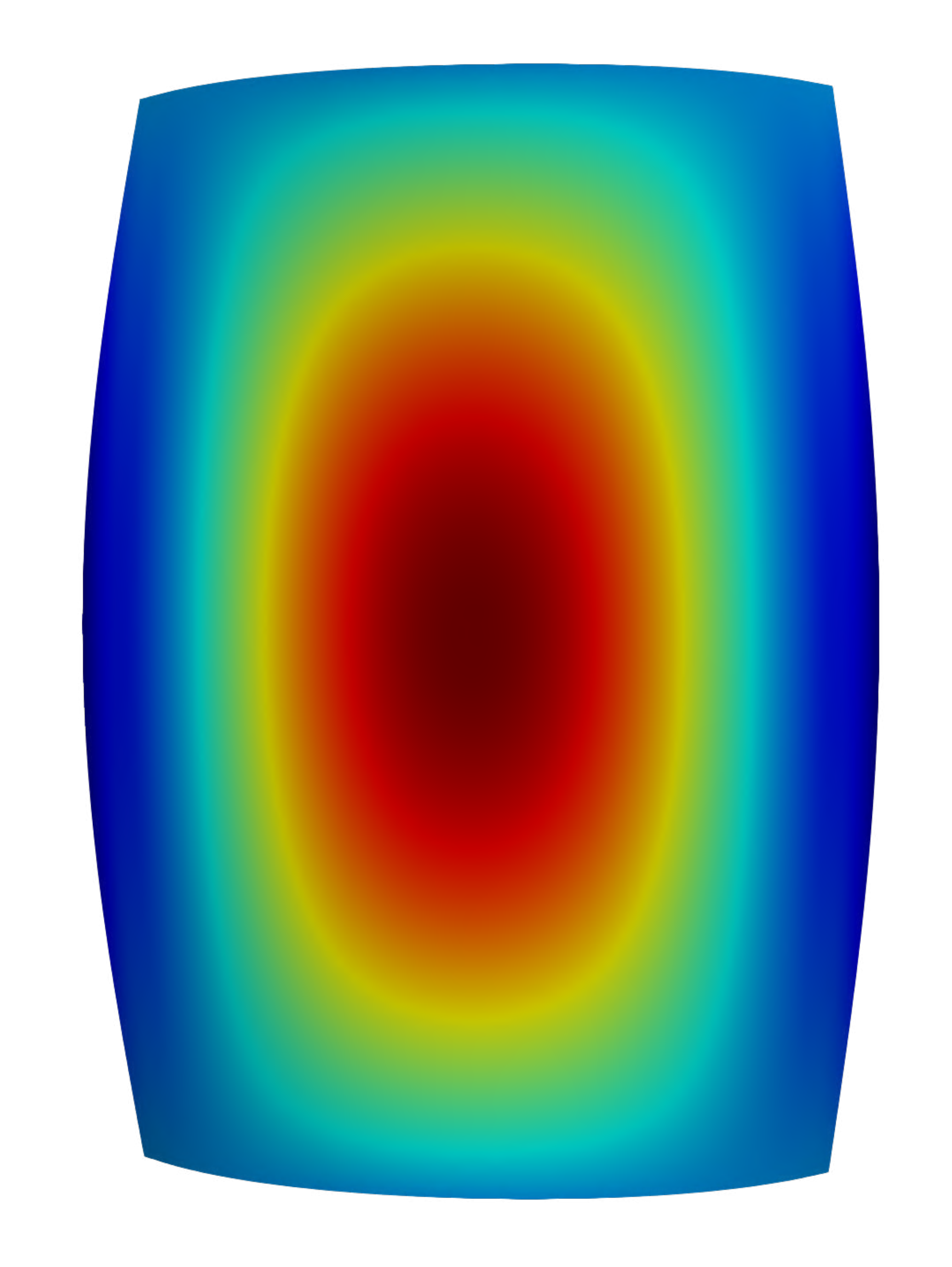}
	\caption{Mode 1: $\varphi$}
	\label{fig:roof_coupled_phi_1}
\end{subfigure}
\begin{subfigure}[b]{0.15\linewidth}
	\centering
 	 \includegraphics[width=\linewidth]{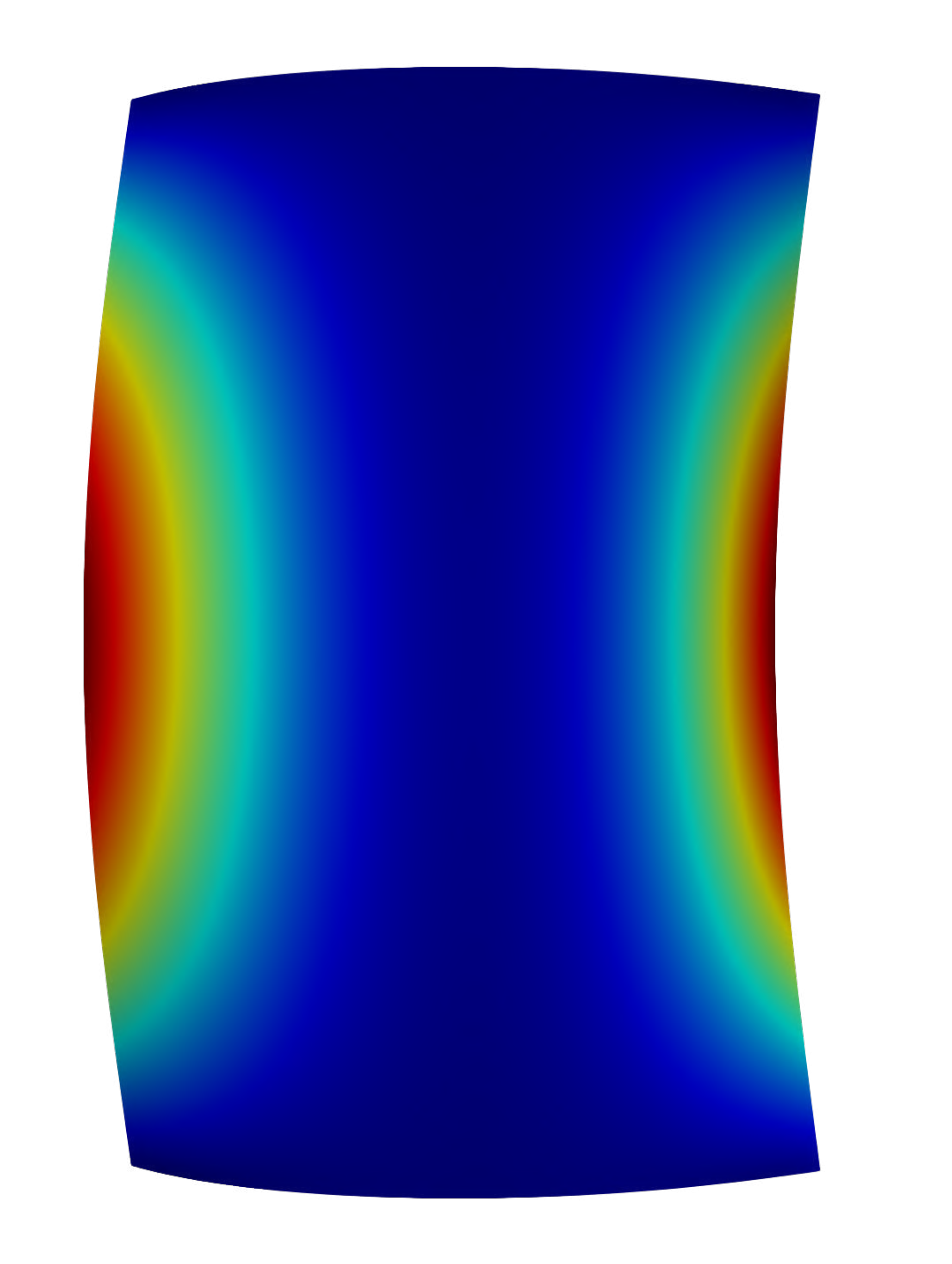}
	\caption{Mode 2:$|\mathbf u|$}
	\label{fig:roof_coupled_disp_2}
\end{subfigure}
\begin{subfigure}[b]{0.15\linewidth}
	\centering
  	\includegraphics[width=\linewidth]{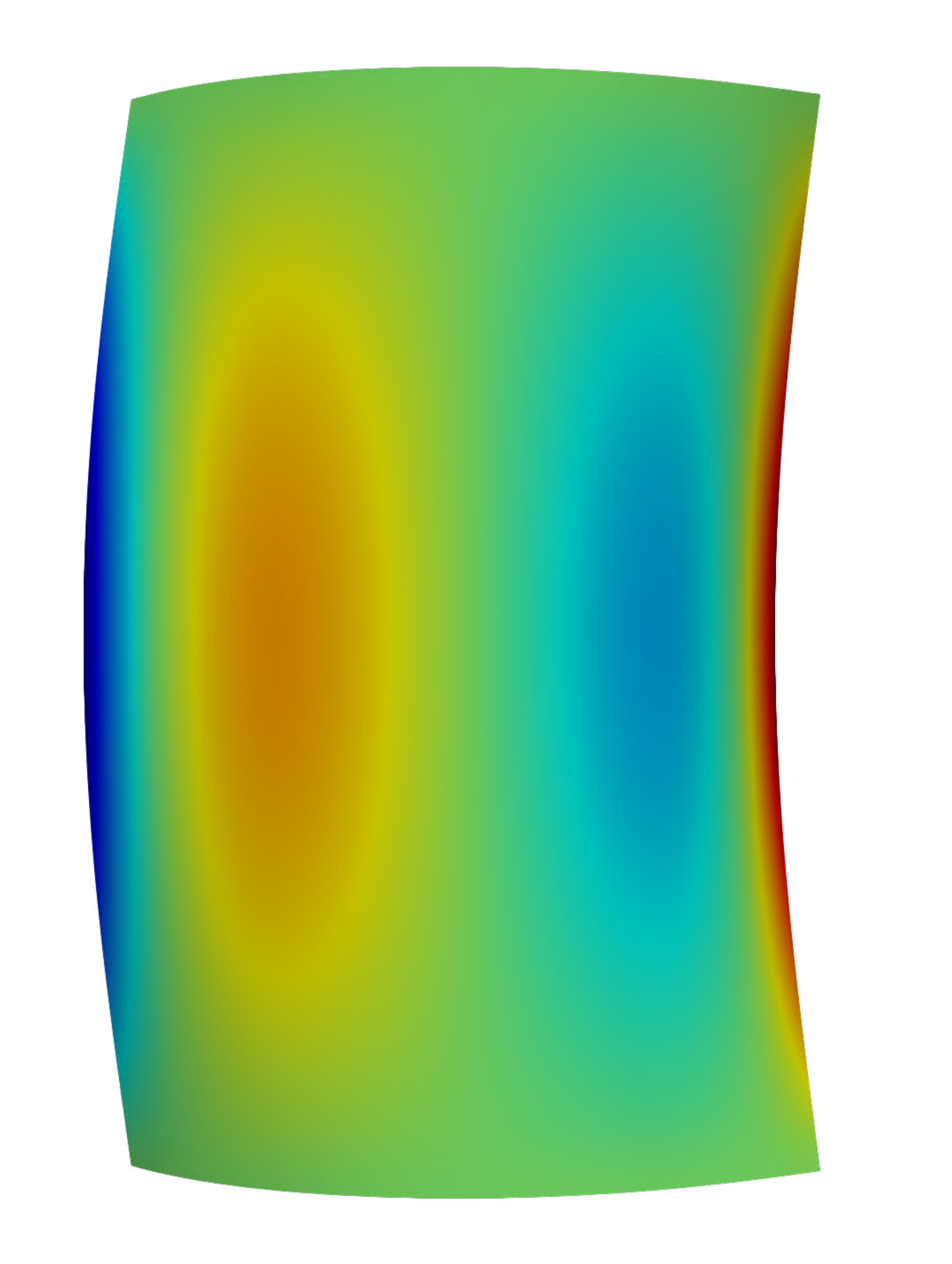}
	\caption{Mode 2: $\psi$}
	\label{fig:roof_coupled_psi_2}
\end{subfigure}
\begin{subfigure}[b]{0.15\linewidth}
	\centering
  	\includegraphics[width=\linewidth]{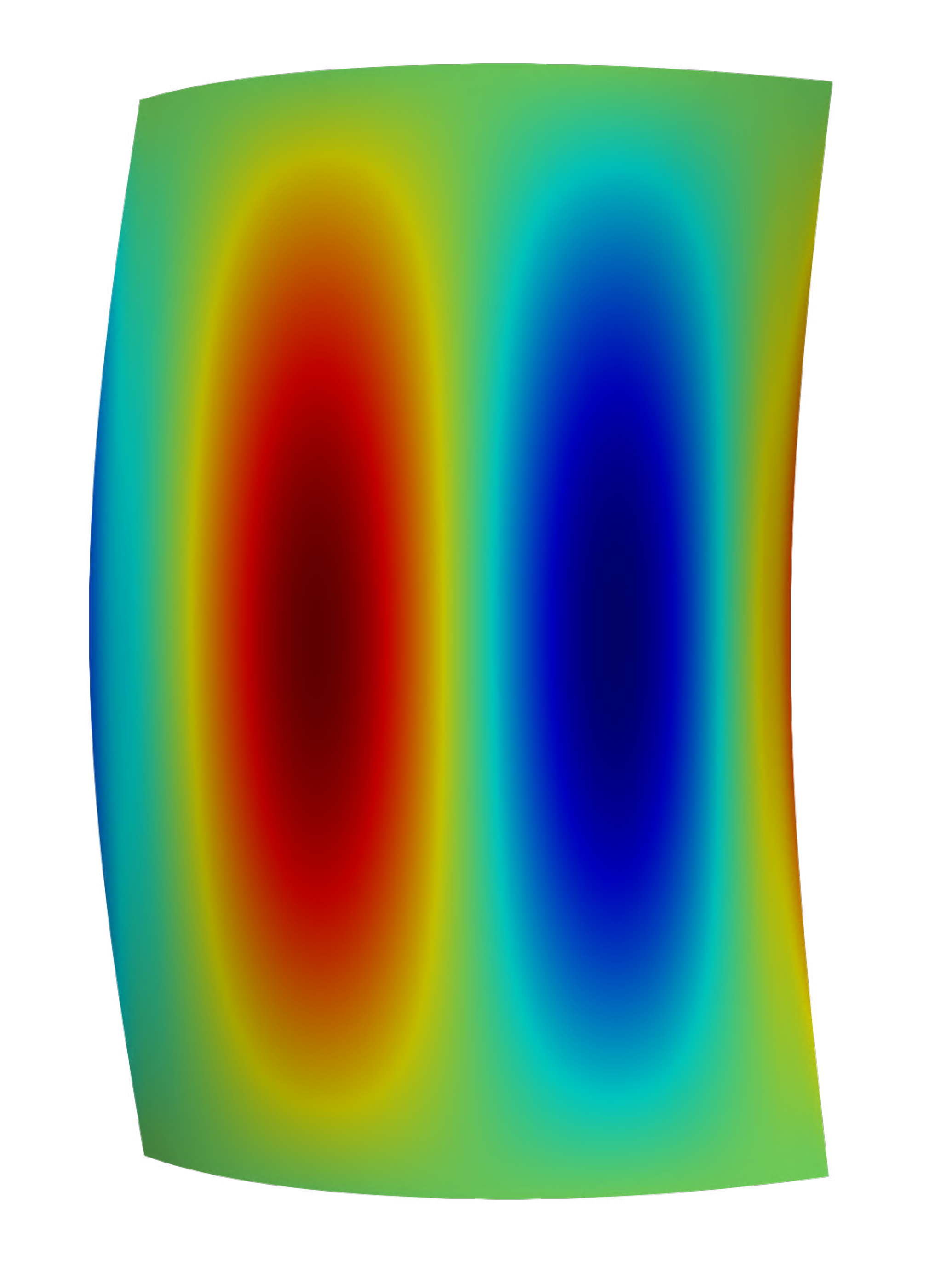}
	\caption{Mode 2: $\varphi$}
	\label{fig:roof_coupled_phi_2}
\end{subfigure}

\begin{subfigure}[b]{0.15\linewidth}
	\centering
 	 \includegraphics[width=\linewidth]{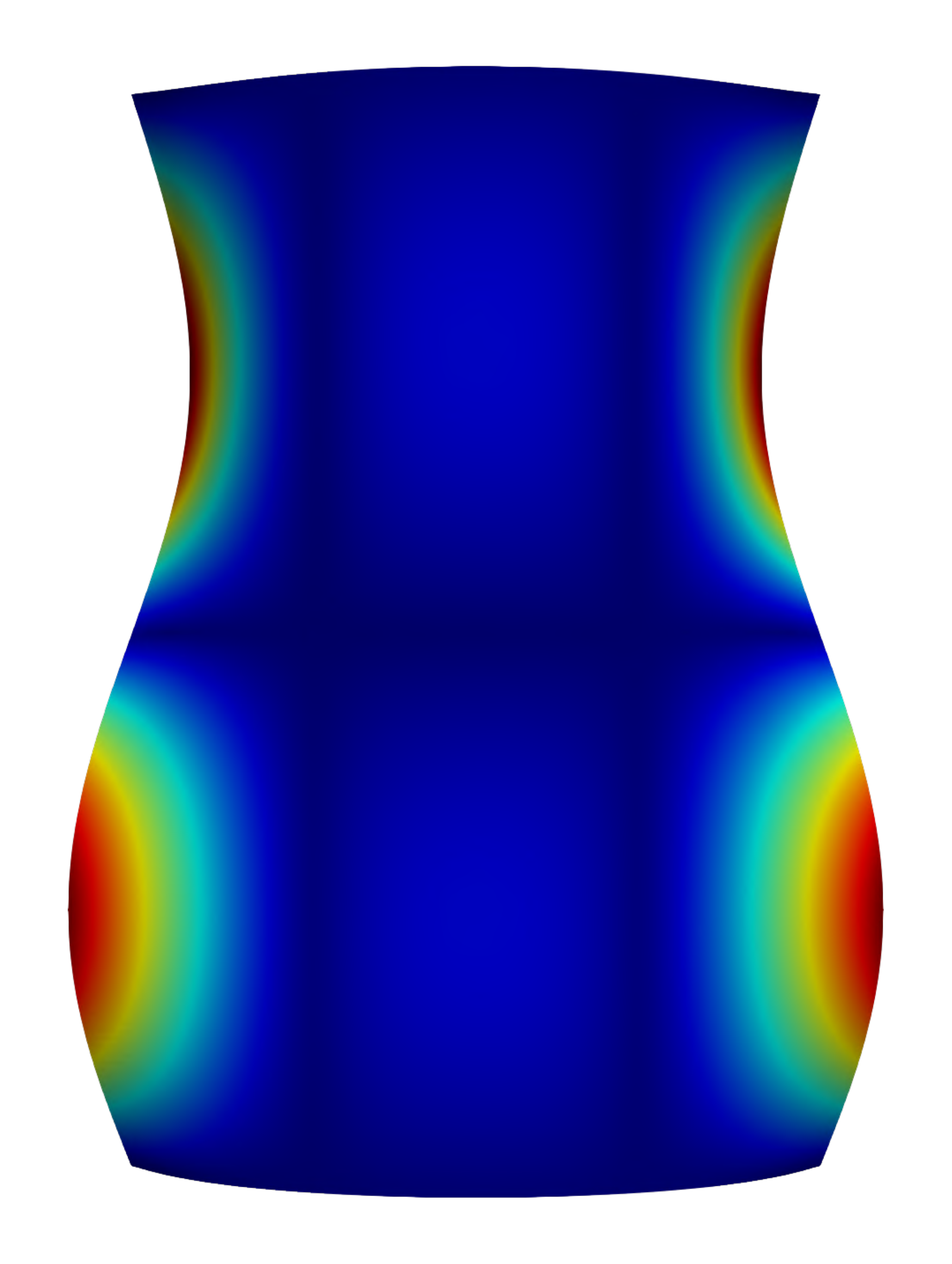}
	\caption{Mode 3:$|\mathbf u|$}
	\label{fig:roof_coupled_disp_3}
\end{subfigure}
\begin{subfigure}[b]{0.15\linewidth}
	\centering
  	\includegraphics[width=\linewidth]{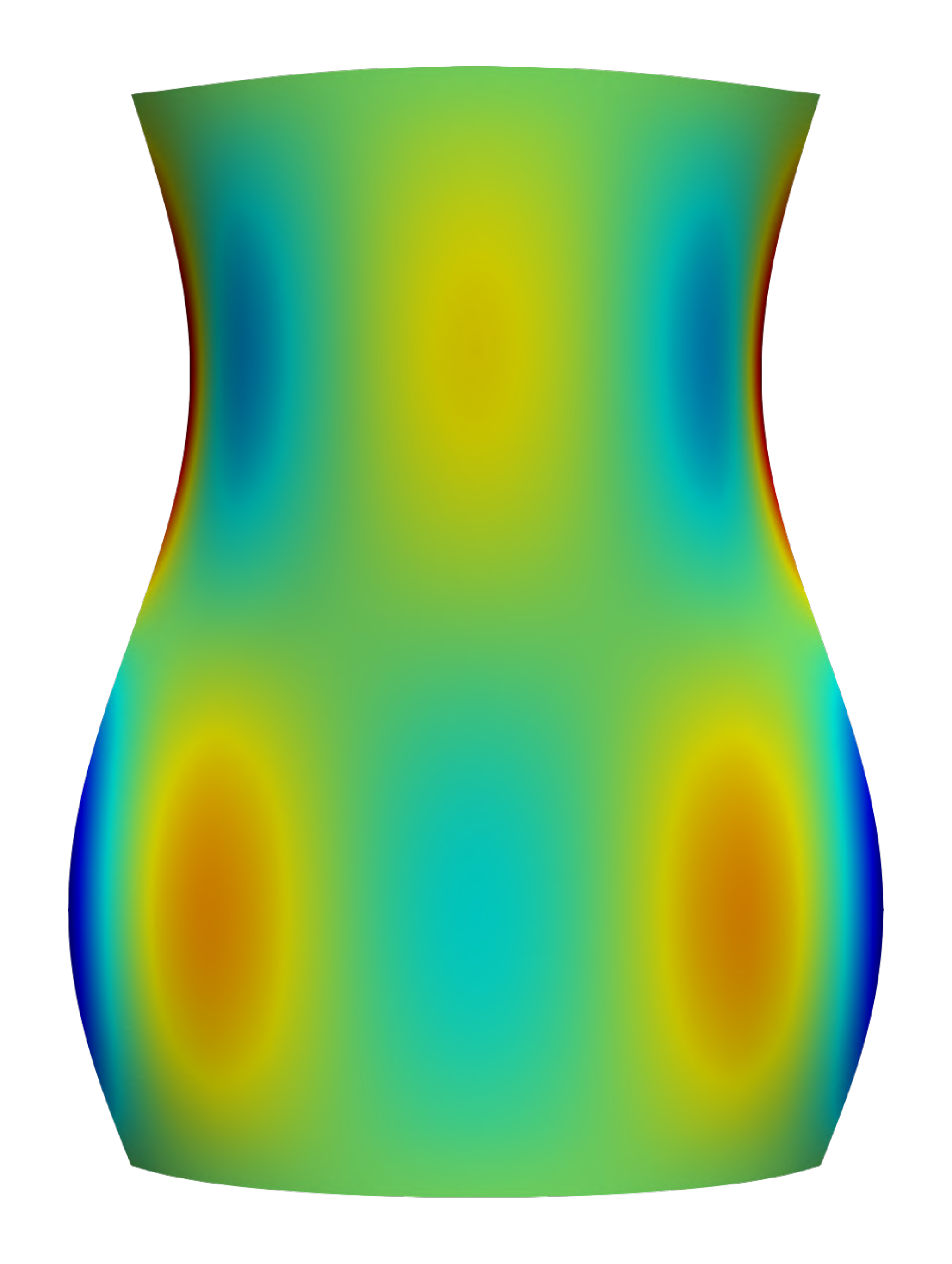}
	\caption{Mode 3: $\psi$}
	\label{fig:roof_coupled_psi_3}
\end{subfigure}
\begin{subfigure}[b]{0.15\linewidth}
	\centering
  	\includegraphics[width=\linewidth]{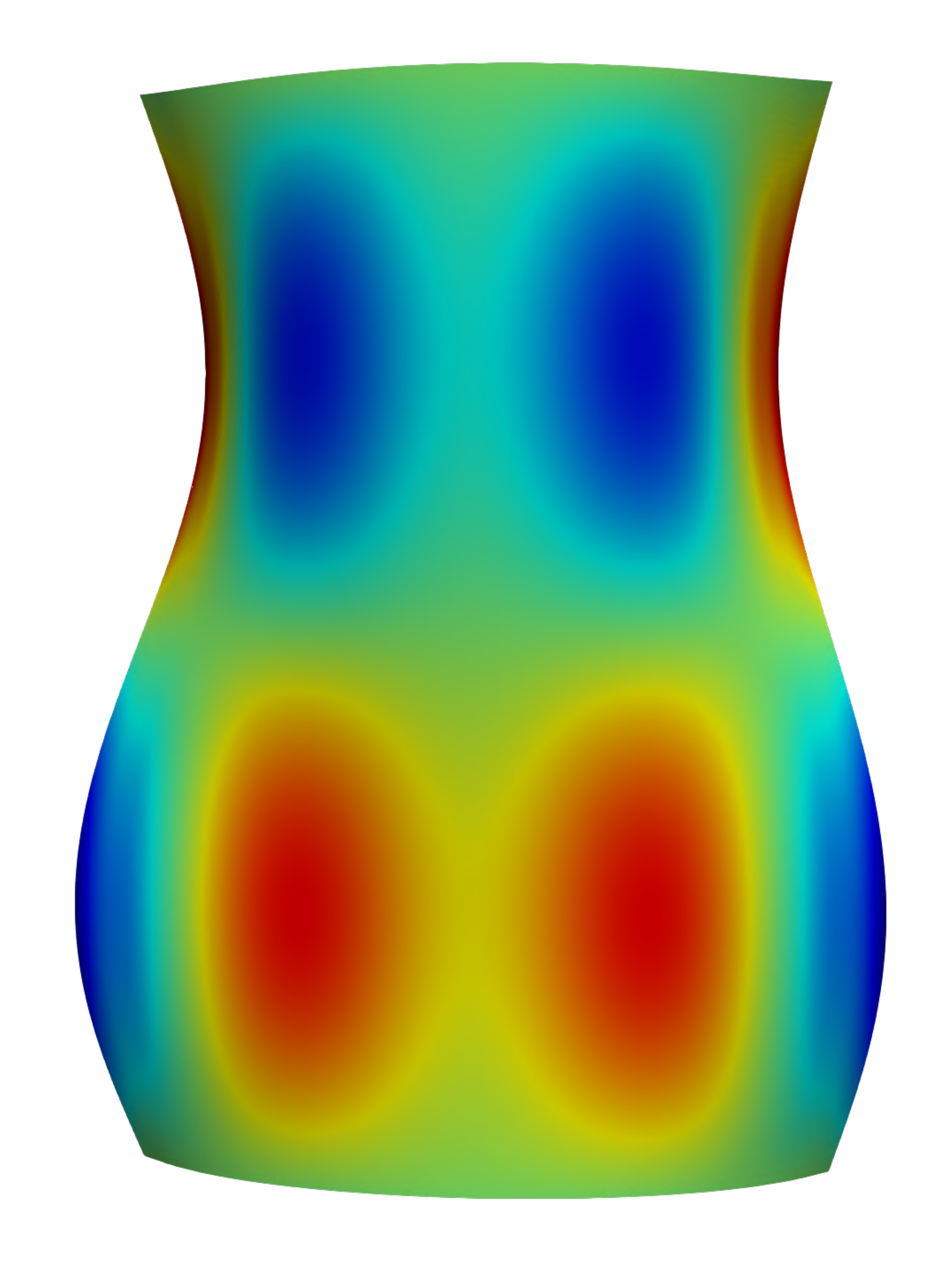}
	\caption{Mode 3: $\varphi$}
	\label{fig:roof_coupled_phi_3}
\end{subfigure}
\begin{subfigure}[b]{0.15\linewidth}
	\centering
 	 \includegraphics[width=\linewidth]{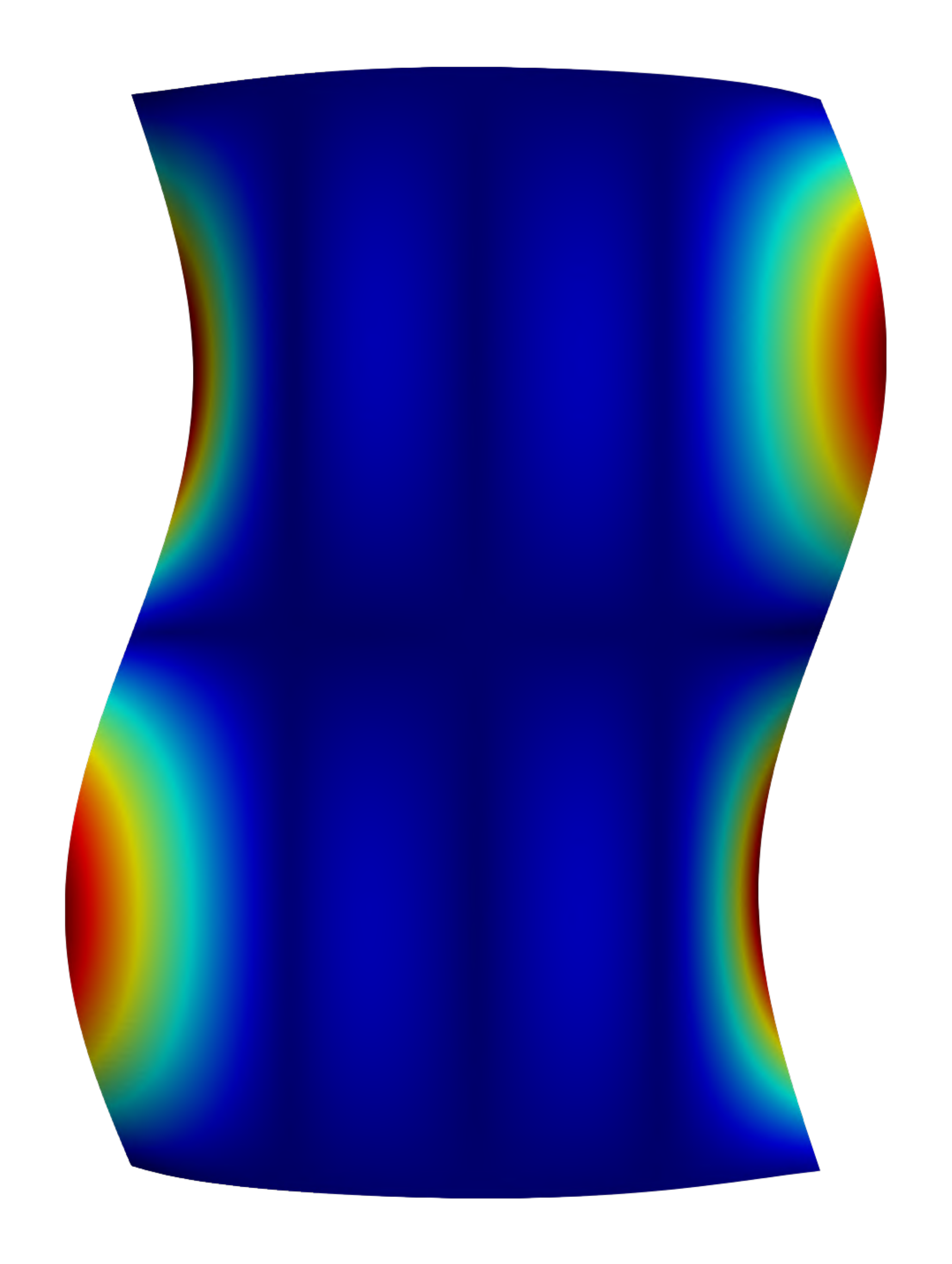}
	\caption{Mode 4:$|\mathbf u|$}
	\label{fig:roof_coupled_disp_4}
\end{subfigure}
\begin{subfigure}[b]{0.15\linewidth}
	\centering
  	\includegraphics[width=\linewidth]{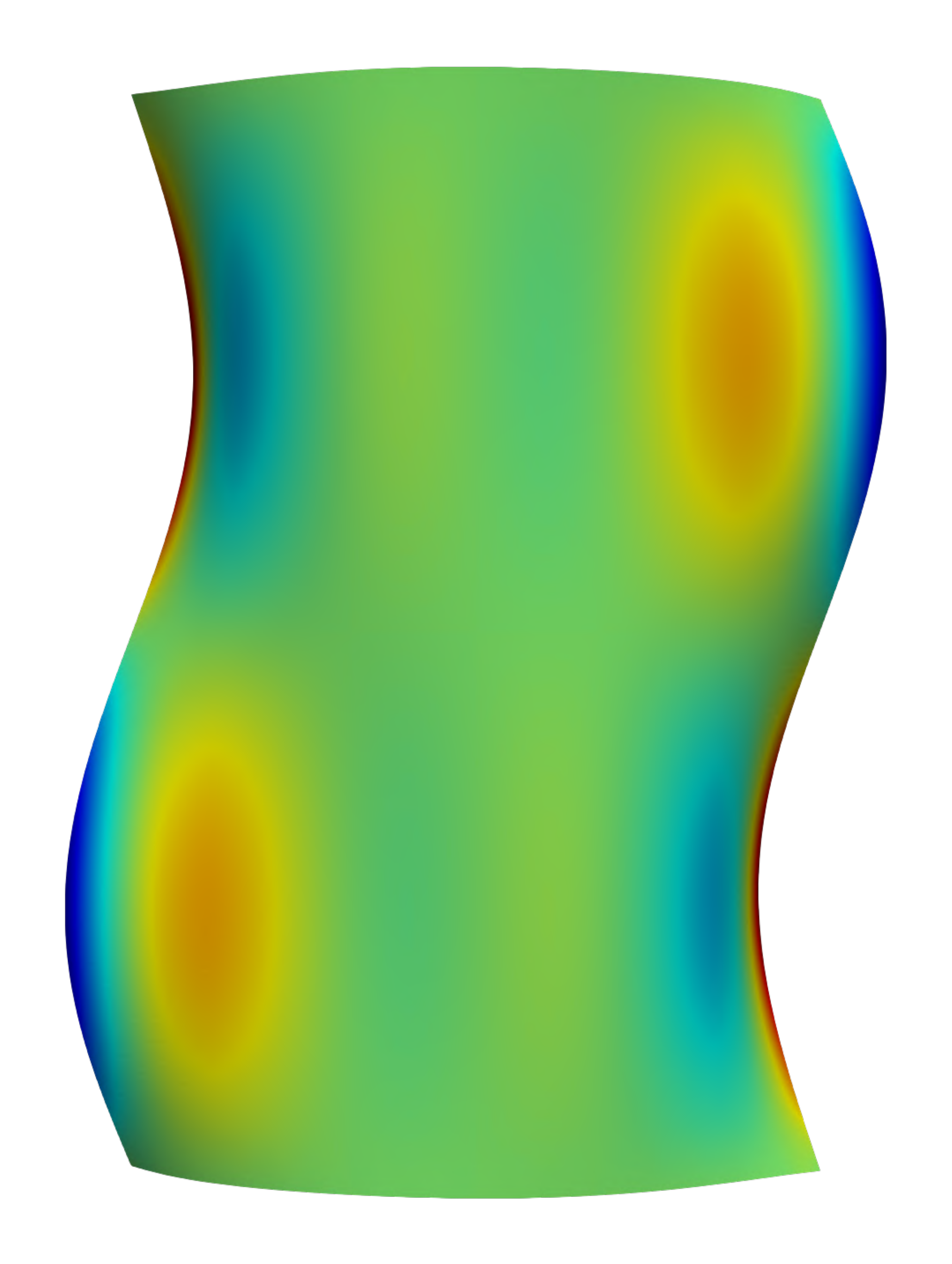}
	\caption{Mode 4: $\psi$}
	\label{fig:roof_coupled_psi_4}
\end{subfigure}
\begin{subfigure}[b]{0.15\linewidth}
	\centering
  	\includegraphics[width=\linewidth]{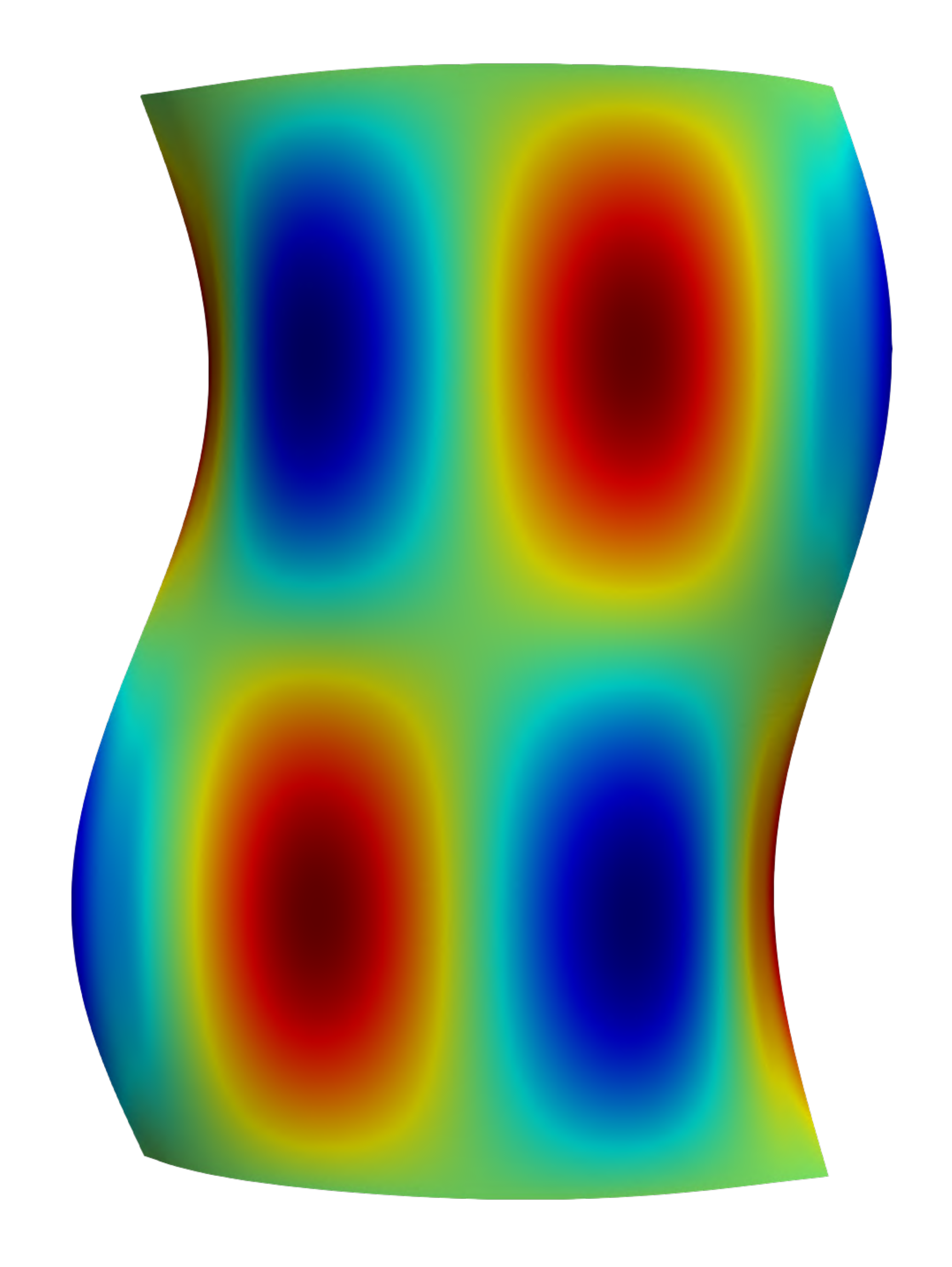}
	\caption{Mode 4: $\varphi$}
	\label{fig:roof_coupled_phi_4}
\end{subfigure}

\begin{subfigure}[b]{0.15\linewidth}
	\centering
 	 \includegraphics[width=\linewidth]{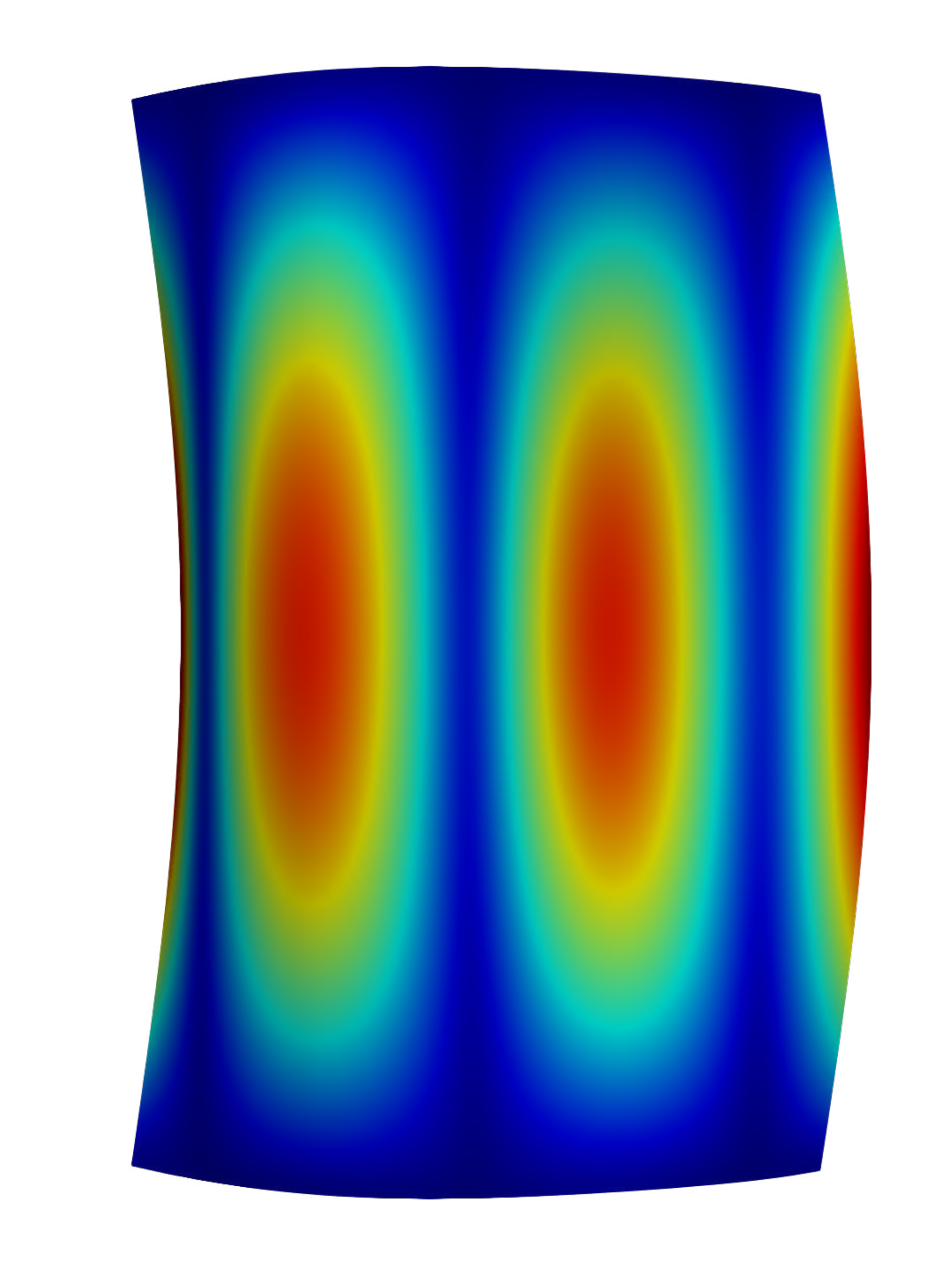}
	\caption{Mode 5:$|\mathbf u|$}
	\label{fig:roof_coupled_disp_5}
\end{subfigure}
\begin{subfigure}[b]{0.15\linewidth}
	\centering
  	\includegraphics[width=\linewidth]{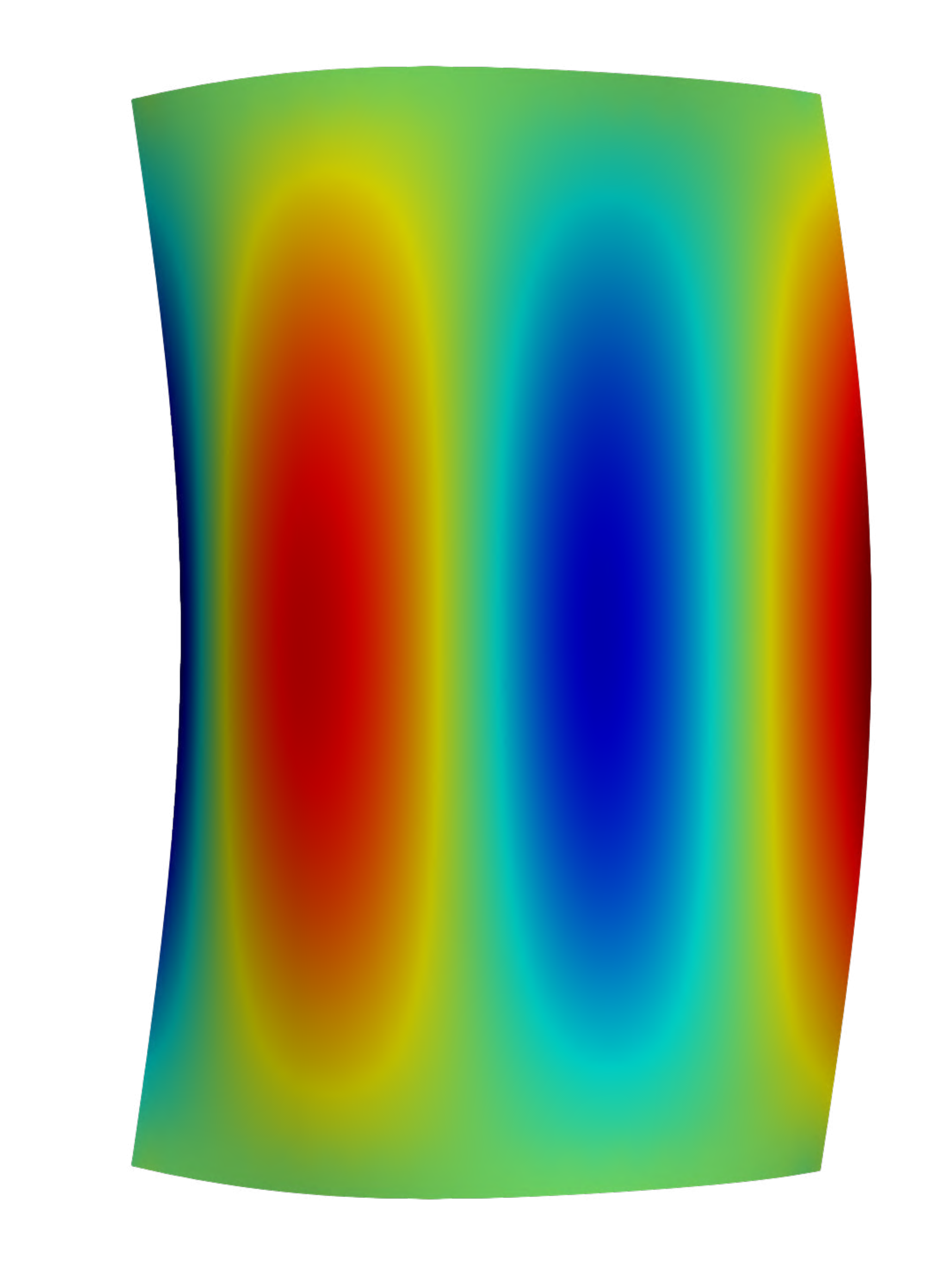}
	\caption{Mode 5: $\psi$}
	\label{fig:roof_coupled_psi_5}
\end{subfigure}
\begin{subfigure}[b]{0.15\linewidth}
	\centering
  	\includegraphics[width=\linewidth]{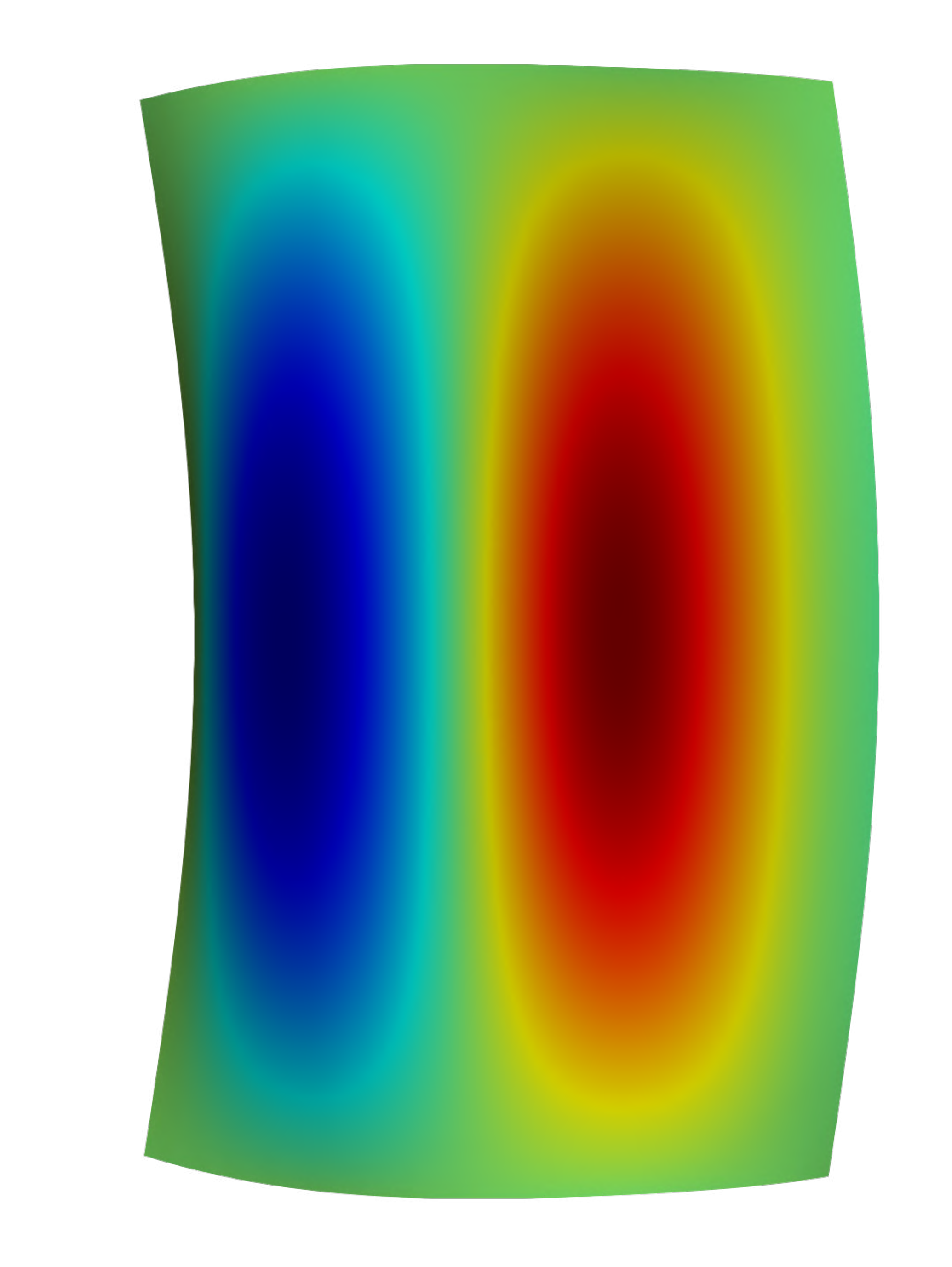}
	\caption{Mode 5: $\varphi$}
	\label{fig:roof_coupled_phi_5}
\end{subfigure}
\begin{subfigure}[b]{0.15\linewidth}
	\centering
 	 \includegraphics[width=\linewidth]{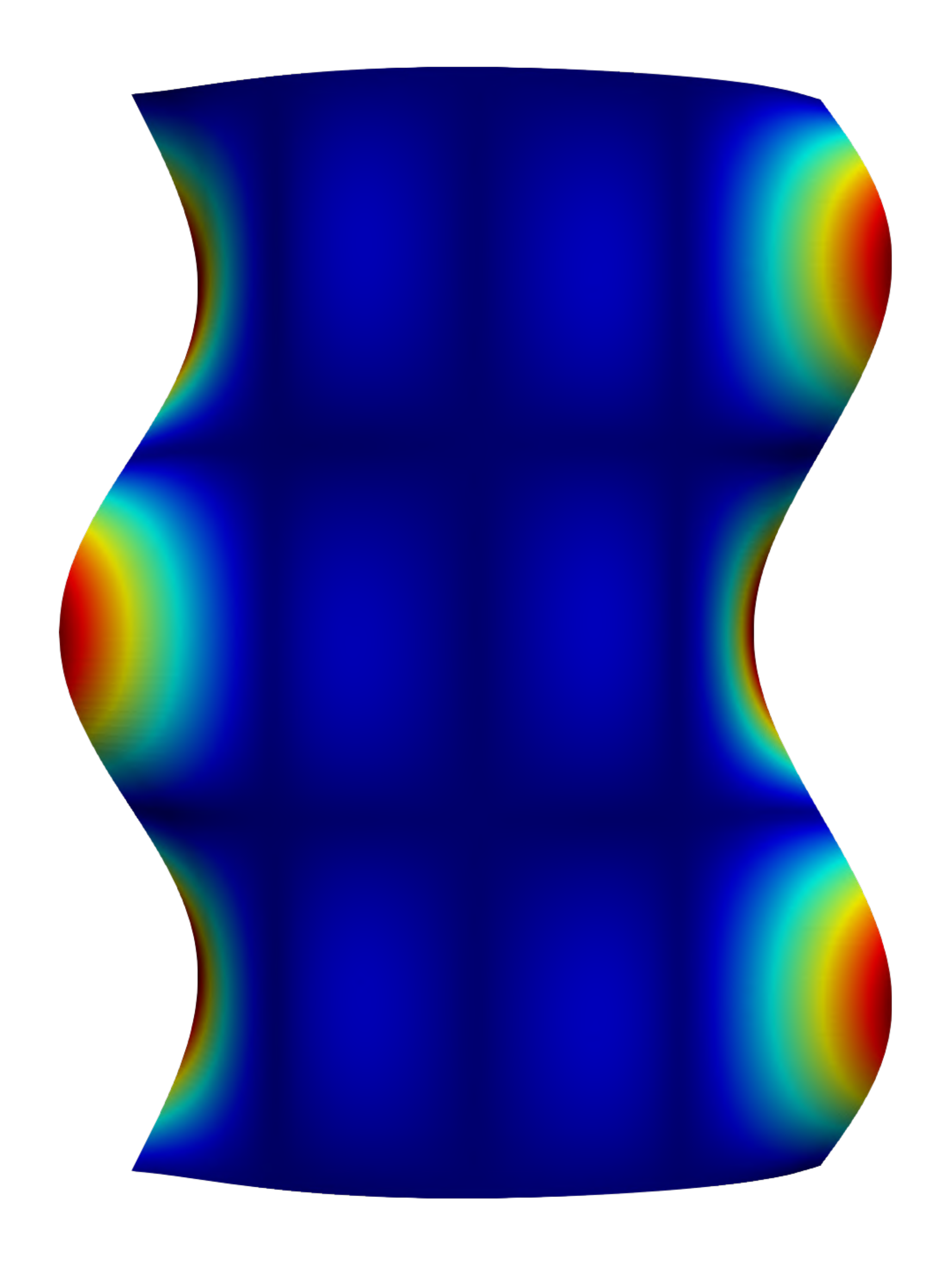}
	\caption{Mode 6:$|\mathbf u|$}
	\label{fig:roof_coupled_disp_6}
\end{subfigure}
\begin{subfigure}[b]{0.15\linewidth}
	\centering
  	\includegraphics[width=\linewidth]{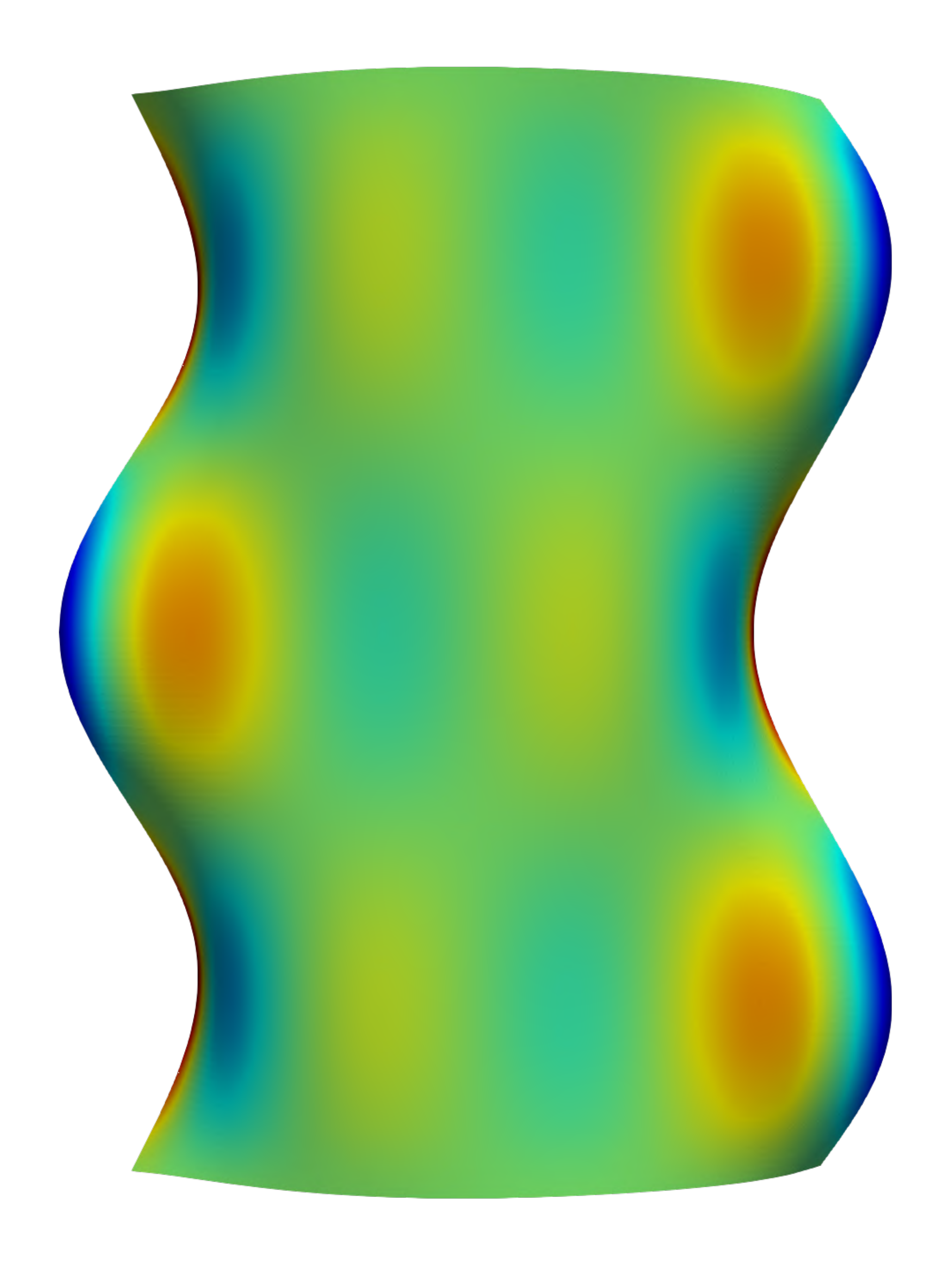}
	\caption{Mode 6: $\psi$}
	\label{fig:roof_coupled_psi_6}
\end{subfigure}
\begin{subfigure}[b]{0.15\linewidth}
	\centering
  	\includegraphics[width=\linewidth]{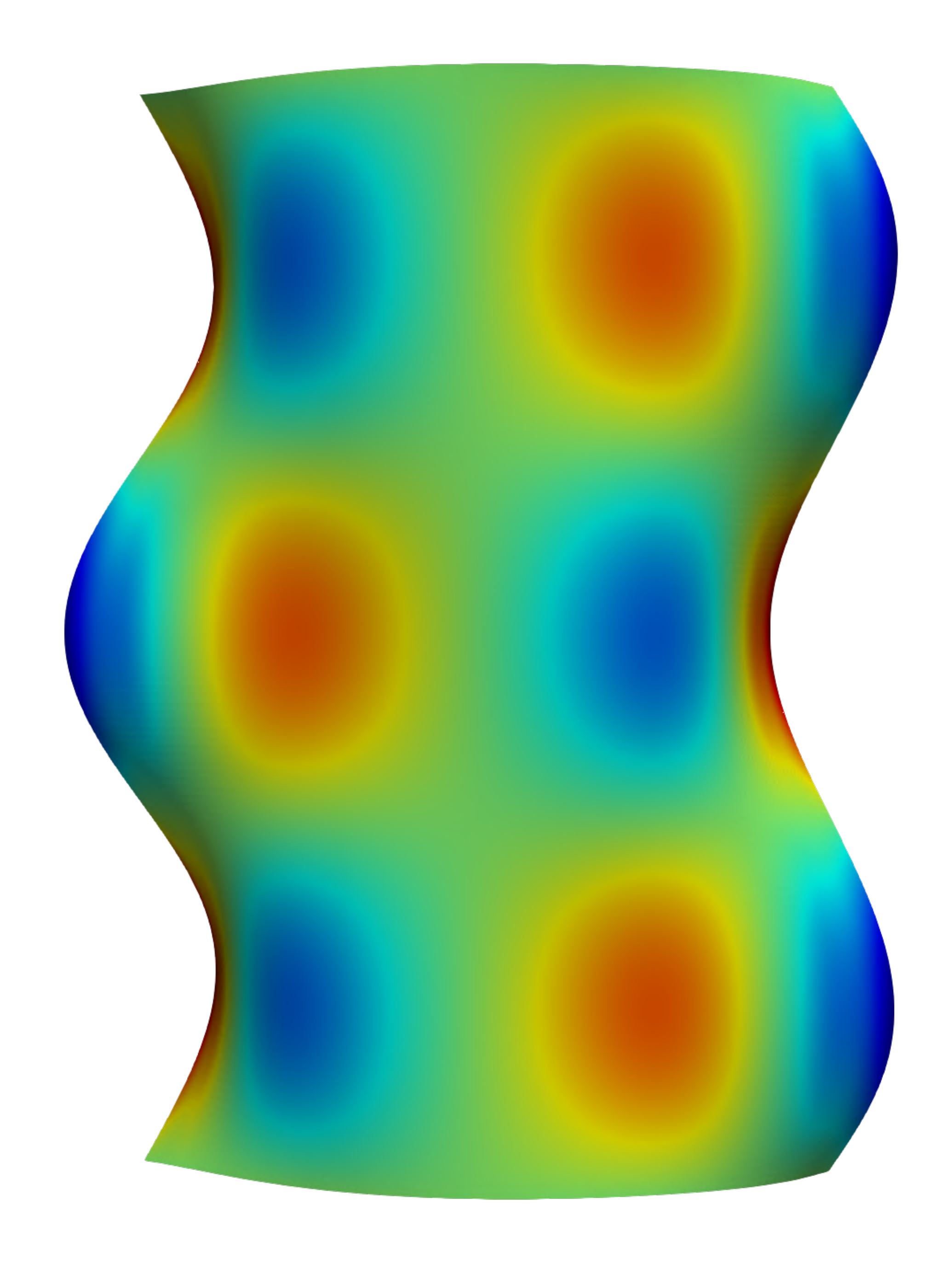}
	\caption{Mode 6: $\varphi$}
	\label{fig:roof_coupled_phi_6}
\end{subfigure}
	\caption{First six vibration modes of the piezoelectric \textit{Scordeli-Lo} roof structure. The magnitude of displacements and the potential functions are plotted on the deformed mid-surface.}
	\label{fig:roof_modes_2}
\end{figure}

The material parameters for the piezoelectric elastic shell are also given in Table~\ref{tab:roof_parameters}. The Benchmark adopted isotropic material but piezoelectric material henceforth are anisotropic. The chosen material BaTiO\textsubscript{3} has a hexagonal crystal system with $6mm$ point group (Hermann–Mauguin notation)~\cite{de2015database}. The piezoelectric tensor $\mathbf e$ has $5$ non-zero components when expressed in Voigt notation~\cite{voigt1928lehrbuch}, which are $e^{31}, e^{32}, e^{33}, e^{15}$ and $e^{24}$. However, since the shell formulation adopts the Kirchhoff-Love and linear elastic assumptions, the components of the strain tensor $S_{13}, S_{23}$ are zero and stress relaxation is used to determine the elastic, piezoelectric and dielectric tensors. The only contributing components in the modified piezoelectric tensor are ${e}^{311}$ and ${e}^{322}$ in the ordinary tensor notation. Figure~\ref{fig:roof_modes_2} shows the first $6$ eigenmodes of a piezoelectric roof-like structure. The magnitude of the displacement and the electric potential functions $\psi$ and $\varphi$ distribution on the piezoelectric shell are plotted. Compared with purely elastic shells, the modal displacements do not exhibit notable change, but the coupling effect will increase the eigenmode frequencies which is known as "piezoelectric stiffening"~\cite{johannsmann2015piezoelectric}. Table~\ref{tab:roof_results} shows the frequency increase of each eigenmode of the short-circuited and unelectroded shells. The increase is more significant for unelectroded shells due to the consideration of the additional linear potential term along the thickness direction.

The coupling effect on the piezoelectric shell with different curvature is also investigated. The arclength $L_{\text{arc}} = 2R\theta$ is held constant. Another two roof-like structure with $\theta = 20^{\circ} (1/9 \pi)$ and $ 60^{\circ} (1/3 \pi)$ are chosen to compare with the original \textit{Scordelis-Lo} roof. The corresponding results are also shown in Table~\ref{tab:roof_results}. The shells with larger curvature have lower frequencies, whereby the rise in frequency is more pronounced for some eigenmodes than for others.

\begin{table}[]
\centering
\caption{Geometry and material parameters of the \textit{Scordelis-Lo} Roof.}
\begin{tabular}{@{}lc@{}}
\toprule
Name        &  BaTiO\textsubscript{3}   \\ \midrule
\textit{Geometry:}\\
Length $L$   & $0.5$m \\
Radius $R$  &  $0.25$m\\
Thickness $h$   &  $2.5 \times 10^{-4} $m\\
Angle $\theta$  &  $20^{\circ},40^{\circ},60^{\circ}$ \\ 
\\
\textit{Material:}\\
Crystal system  & hexagonal (6mm) \\
Mass density  $\rho$   & $5800 \text{ kg}/\text{m}^3$ \\
Elastic constants & \\
$\tilde C^{1111}, \tilde C^{2222}$  & $166$ GPa\\
$\tilde C^{1122}, \tilde C^{2211}$  &$77$ GPa \\
$\tilde C^{3333}$  & $162$ GPa\\
$\tilde C^{1133}, \tilde C^{3311}, \tilde C^{2233}, \tilde C^{3322}$  & $78$ GPa \\
$\tilde C^{1212}, \tilde C^{1221}, \tilde C^{2121}, \tilde C^{2112}$  &$45$ GPa \\

Piezoelectric constants & \\
$\tilde e^{311}, \tilde e^{322}$ & $ -4.4 \text{ C}/\text{m}^2$\\
%$\tilde e^{312}, \tilde e^{321}$ & - & $-1.45031 \text{ C}/\text{m}^2$\\
$\tilde e^{333}$  & $18.6 \text{ C}/\text{m}^2$\\ 
%$\tilde e^{333}$ & $15.1$ \\
%$\tilde e^{113}, \tilde e^{223}$ & $12.7$ \\
Permittivity  & \\
$\tilde\kappa^{11} , \tilde\kappa^{22}$  & $11.2\times 10^{-9}$ C\textsuperscript{2} / (Nm\textsuperscript{2})\\
$\tilde\kappa^{33}$  & $12.6\times 10^{-9}$ C\textsuperscript{2} / (Nm\textsuperscript{2})\\
 \bottomrule
\end{tabular}
\label{tab:roof_parameters}
\end{table}

\begin{table}[]
\centering
\caption{Eigenmode frequencies for the elastic and the piezoelectric roof-like shells with different curvatures.(SC stands for short-circuited shell and UE denotes for unelectroded shell.)}
\begin{tabular}{@{}CCCC|CCC|CCC@{}}
\toprule
\multirow{2}{*}{Mode} & \multicolumn{3}{c|}{$f(\theta = 20^{\circ}, R = 50)$(Hz)}            & \multicolumn{3}{c|}{$f(\theta = 40^{\circ}, R = 25)$(Hz)}               & \multicolumn{3}{c}{$f(\theta = 60^{\circ}, R = 50/3)$(Hz)}               \\
                      & $Elastic$  & $SC$ &$UE$  & $Elastic$  & $SC$ & $UE$ & $Elastic$ & $SC$ &$UE$  \\ \midrule
1                     & 82.32          & 82.45      & 83.68        & 125.31         & 127.85        & 128.62        & 143.64          & 145.89      & 147.27         \\
2                     & 109.49          & 111.95      & 112.39         & 133.59          & 134.33        & 136.42       & 170.50          & 172.88      & 174.91         \\
3                     & 214.49          & 217.09       & 218.57        & 283.21          & 285.42        & 288.56       & 334.96          & 339.55      &  342.36       \\
4                     & 229.78          & 231.52      & 233.90       & 293.62          & 297.83        & 299.95       & 343.63         & 347.06      & 350.77         \\
5                     & 275.65          & 276.81       & 282.35       & 336.33          & 346.90        & 348.83       & 377.23         & 385.11      & 389.34         \\
6                     & 311.73          & 324.00       & 324.45       & 470.41          & 475.82        & 479.50        & 545.56 & 551.02      & 556.07         \\
7                     & 369.39          & 373.67       & 376.50       & 477.89          & 482.63        & 486.82       & 551.55 & 558.05      & 562.57         \\
8                     & 382.86         & 387.51       & 390.19       & 493.33          & 494.79        & 504.94       & 554.37          & 576.33      &   577.58       \\ \bottomrule
\end{tabular}
\label{tab:roof_results}
\end{table}

\subsection{Free vibration of a piezoelectric speaker}
\begin{figure}[!t]
\centering
\begin{subfigure}[b]{0.38\linewidth}
	\centering
 	 \includegraphics[width=\linewidth]{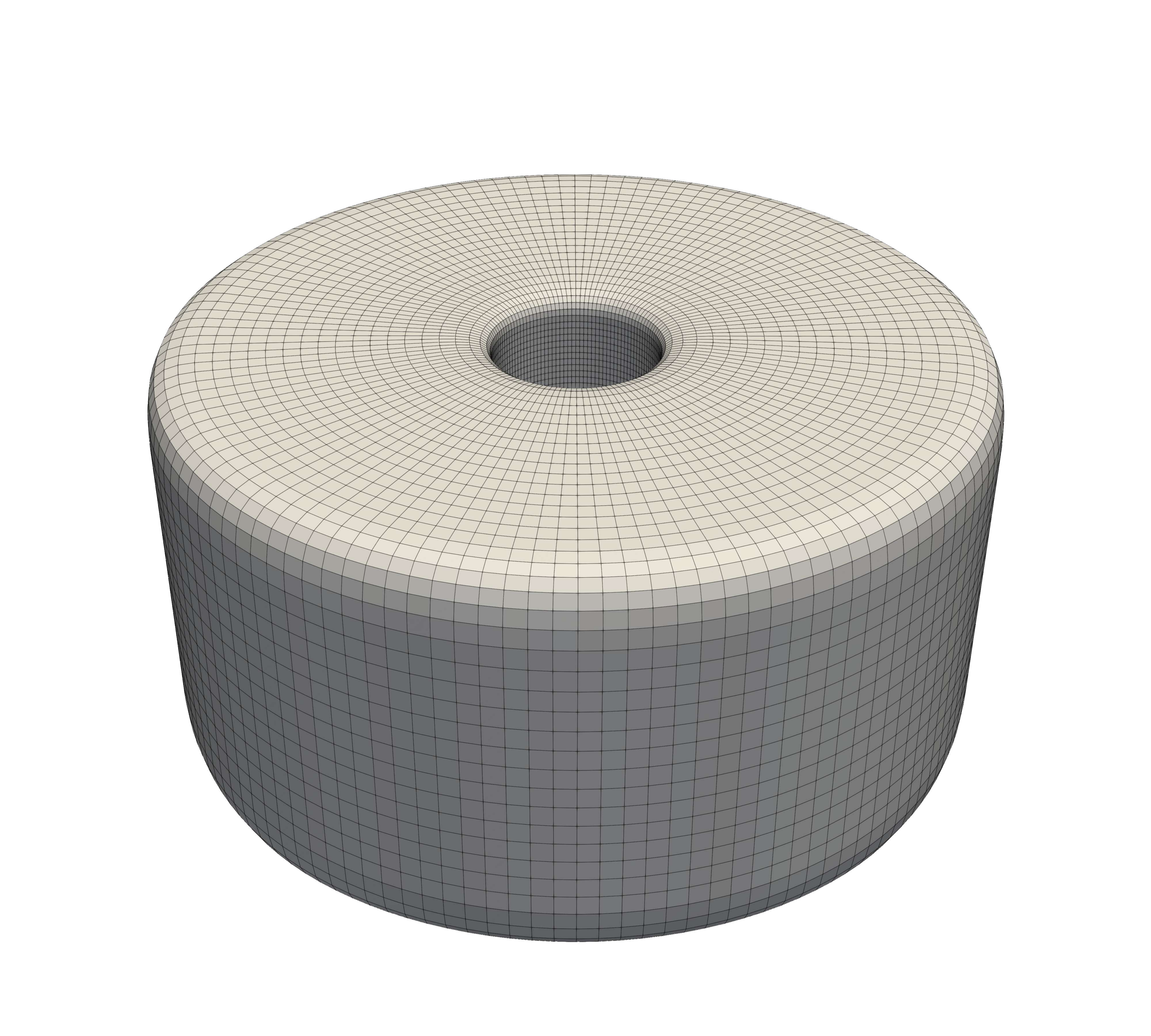}
	\caption{}
	\label{fig:piezo_buzzer_a}
\end{subfigure}
\begin{subfigure}[b]{0.4\linewidth}
	\centering
 	 \includegraphics[width=\linewidth]{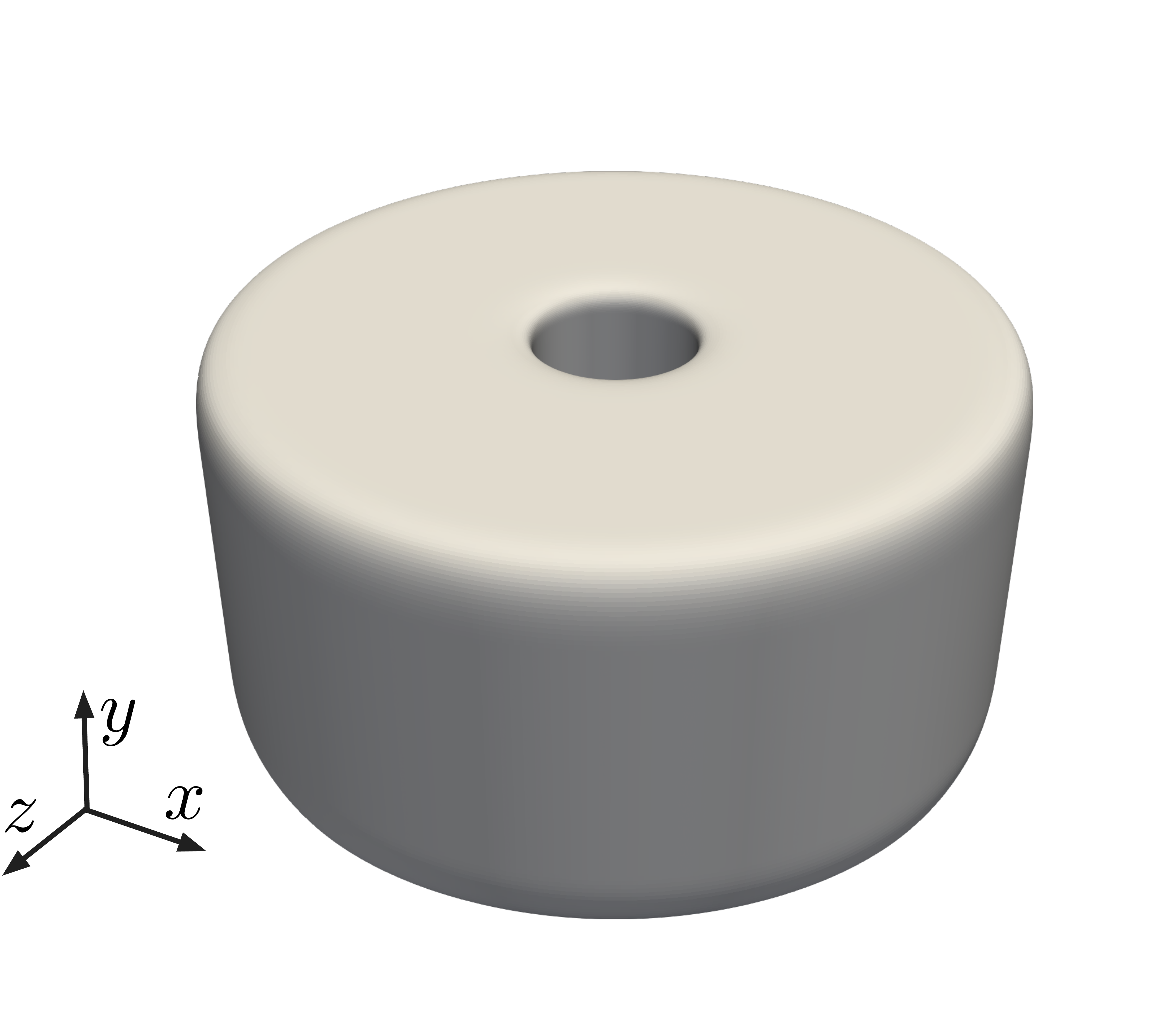}
	\caption{}
	\label{fig:piezo_buzzer_b}
\end{subfigure}
\begin{subfigure}[b]{0.4\linewidth}
	\centering
 	 \includegraphics[width=0.75\linewidth]{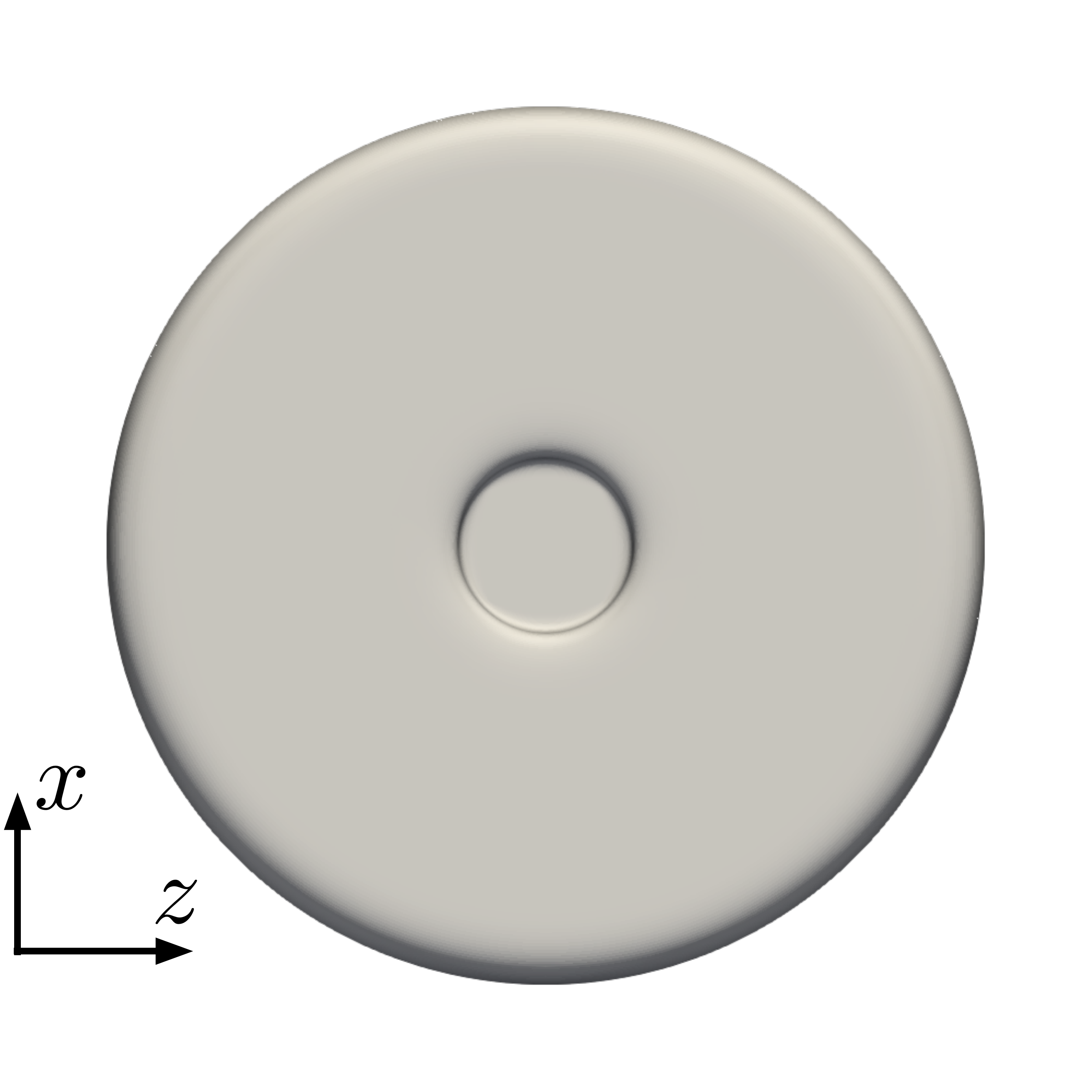}
	\caption{}
	\label{fig:piezo_buzzer_c}
\end{subfigure}
\begin{subfigure}[b]{0.4\linewidth}
	\centering
 	 \includegraphics[width=0.75\linewidth]{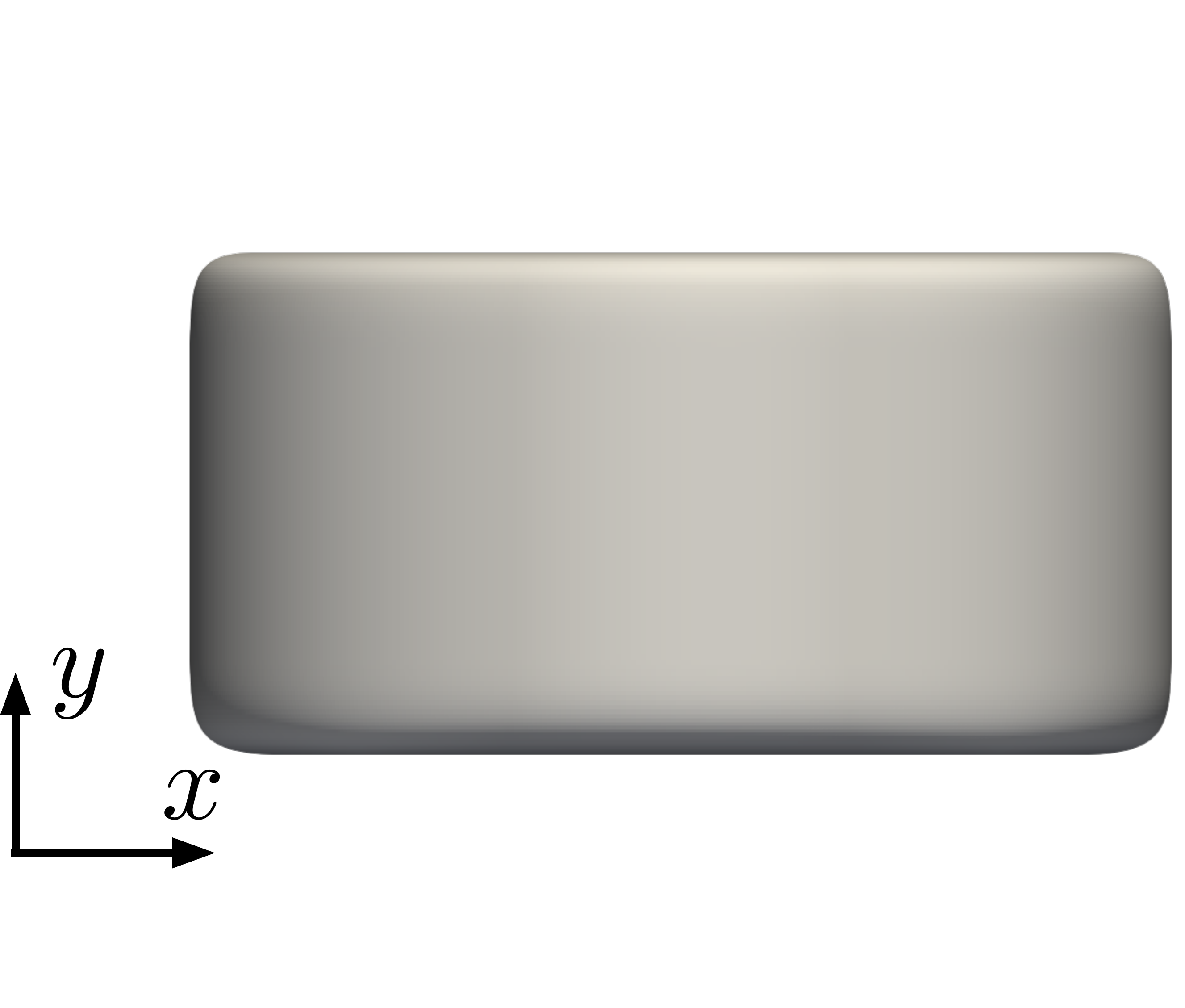}
	\caption{}
	\label{fig:piezo_buzzer_d}
\end{subfigure}

\caption{ A piezoelectric buzzer geometry. (a) is a mesh of the shell with 12288 elements. (b) is the limit subdivision surface constructed using (a). (c) is the top view of the shell and (d) represents both the front and side view of the axisymmetric geometry.}
\label{fig:piezo_buzzer_3d}
\end{figure}

The final example considers a potential application to a piezoelectric speaker made from a single shell. The geometry considered is regenerated from a CAD model of a piezoelectric speaker. It is imported into Autodesk Maya~\cite{maya} for removal of extraneous geometry. A quadrilateral control mesh for the geometry is shown in Figure~\ref{fig:piezo_buzzer_a}. A model based on Catmull–Clark subdivision surface can directly evaluate the smooth limit surface in Figure~\ref{fig:piezo_buzzer_b} using the control mesh. The limiting surface is smooth everywhere. Figure~\ref{fig:piezo_buzzer_c} and~\ref{fig:piezo_buzzer_d} are the top and front view of the geometry. The minimum bounding box for this model is defined by $[x_i^{min}, x_i^{max}]^3 = [-0.0694,0.0694] \times [0, 0.0711] \times [-0.0694,0.0694]$m\textsuperscript{3}. The geometry is axisymmetric about the $y$-axis. The thickness of the shell is $0.002$m. The eigenvalue analysis with no boundary constraint is performed for this example and the same transverse isotropic piezoelectric material BaTiO\textsubscript{3} is chosen. The unelectroded condition is used. 

\begin{figure}[!t]
\centering
\begin{subfigure}[b]{0.3\linewidth}
	\centering
 	 \includegraphics[width=\linewidth]{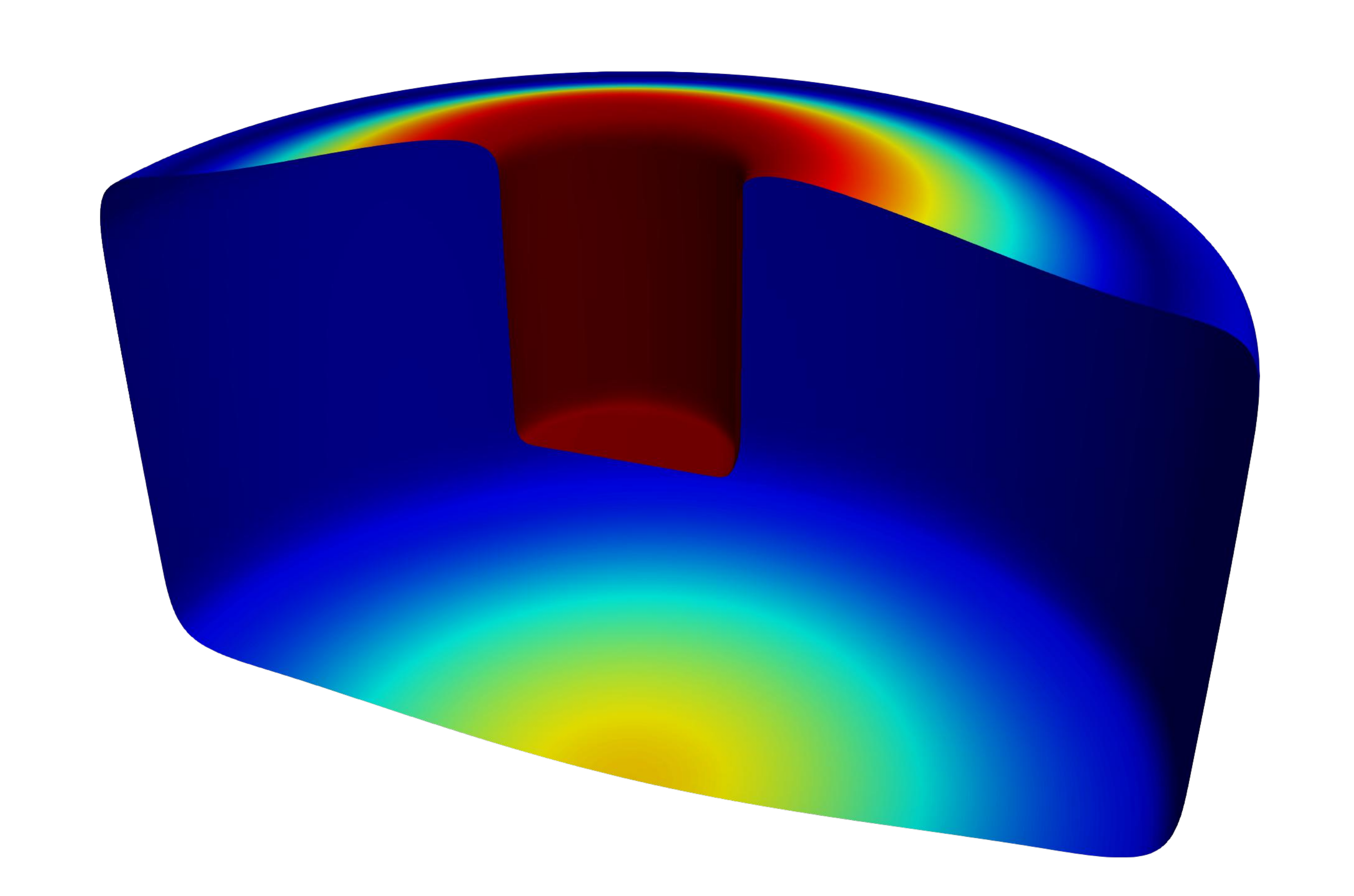}
	\caption{Mode 1: $|\mathbf u|$ on the shell.}
	\label{fig:piezo_buzzer_mode_disp_1}
\end{subfigure}
\begin{subfigure}[b]{0.3\linewidth}
	\centering
 	 \includegraphics[width=\linewidth]{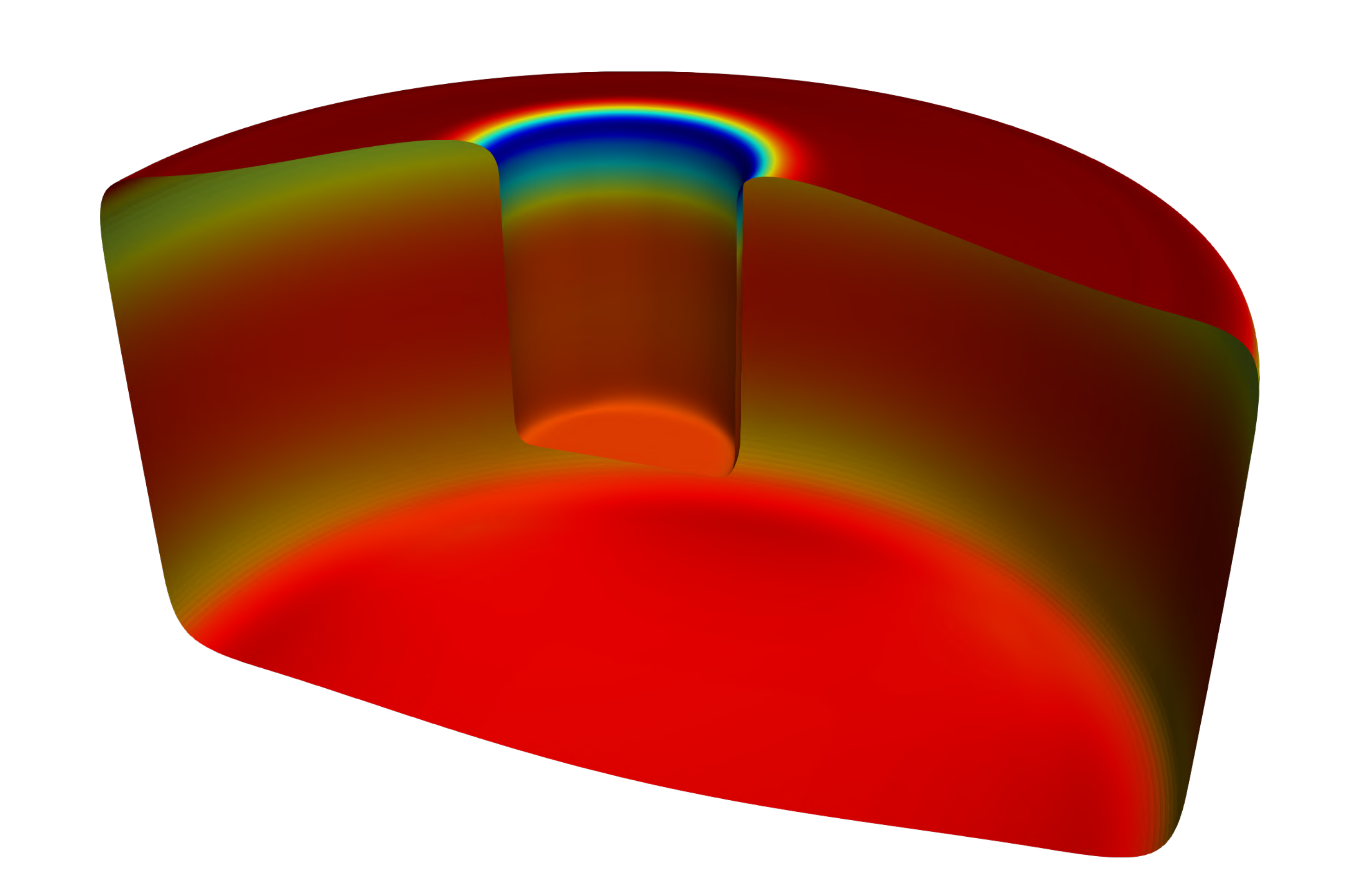}
	\caption{Mode 1: $\psi$ on the shell.}
	\label{fig:piezo_buzzer_mode_psi_1}
\end{subfigure}
\begin{subfigure}[b]{0.3\linewidth}
	\centering
 	 \includegraphics[width=\linewidth]{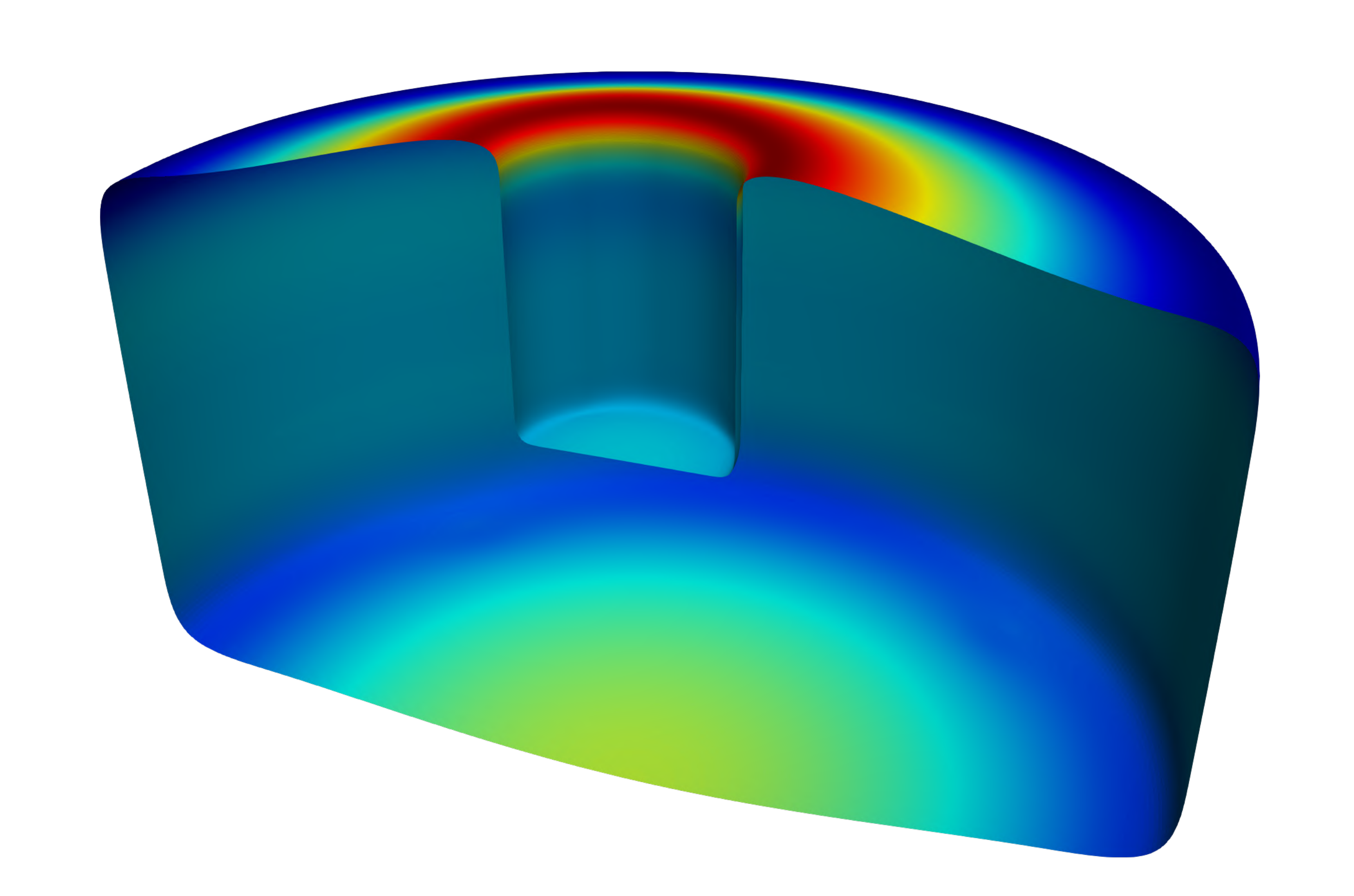}
	\caption{Mode 1: $\varphi$ on the shell.}
	\label{fig:piezo_buzzer_mode_phi_1}
\end{subfigure}

\begin{subfigure}[b]{0.3\linewidth}
	\centering
 	 \includegraphics[width=\linewidth]{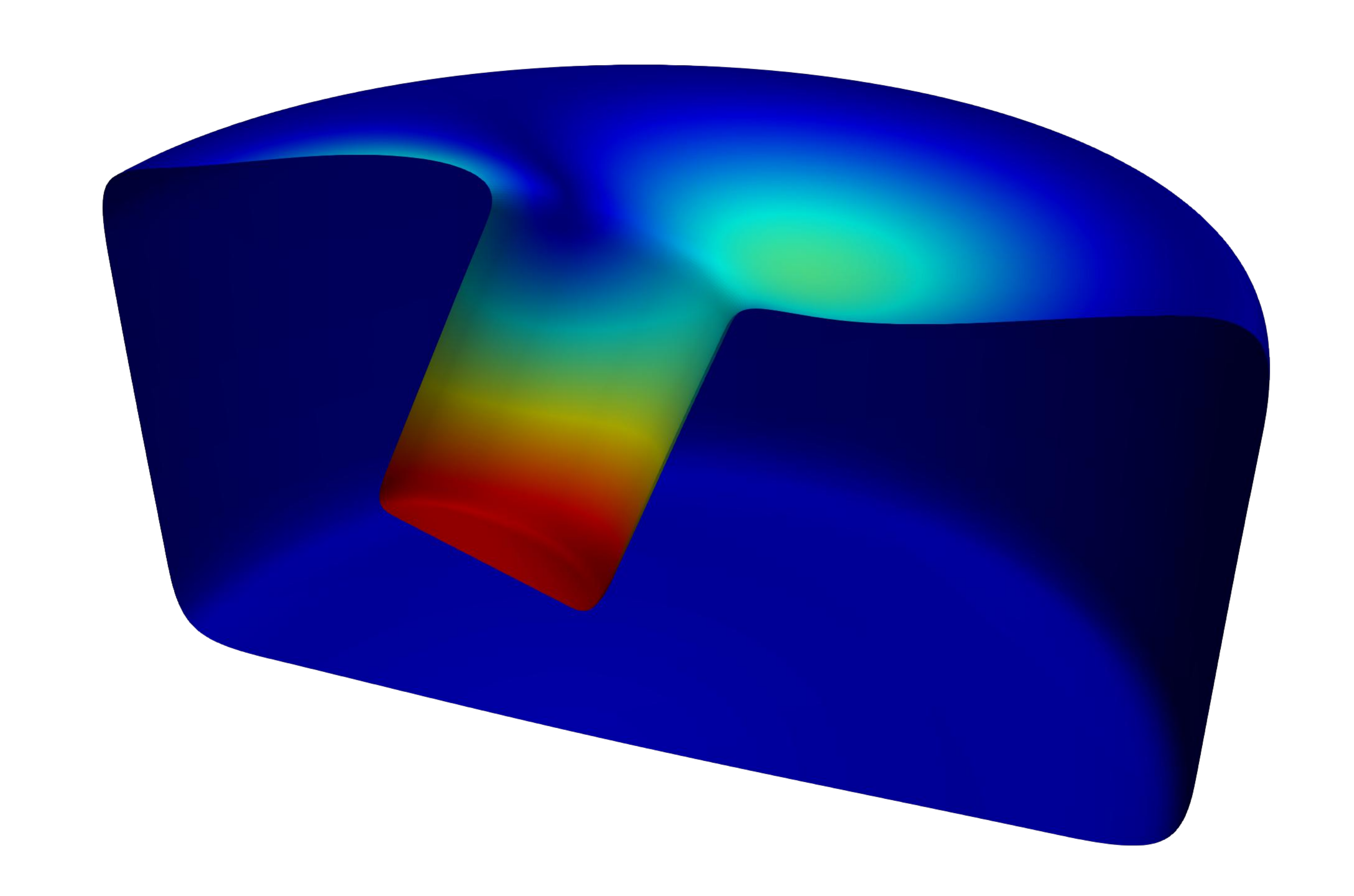}
	\caption{Mode 2: $|\mathbf u|$ on the shell.}
	\label{fig:piezo_buzzer_mode_disp_2}
\end{subfigure}
\begin{subfigure}[b]{0.3\linewidth}
	\centering
 	 \includegraphics[width=\linewidth]{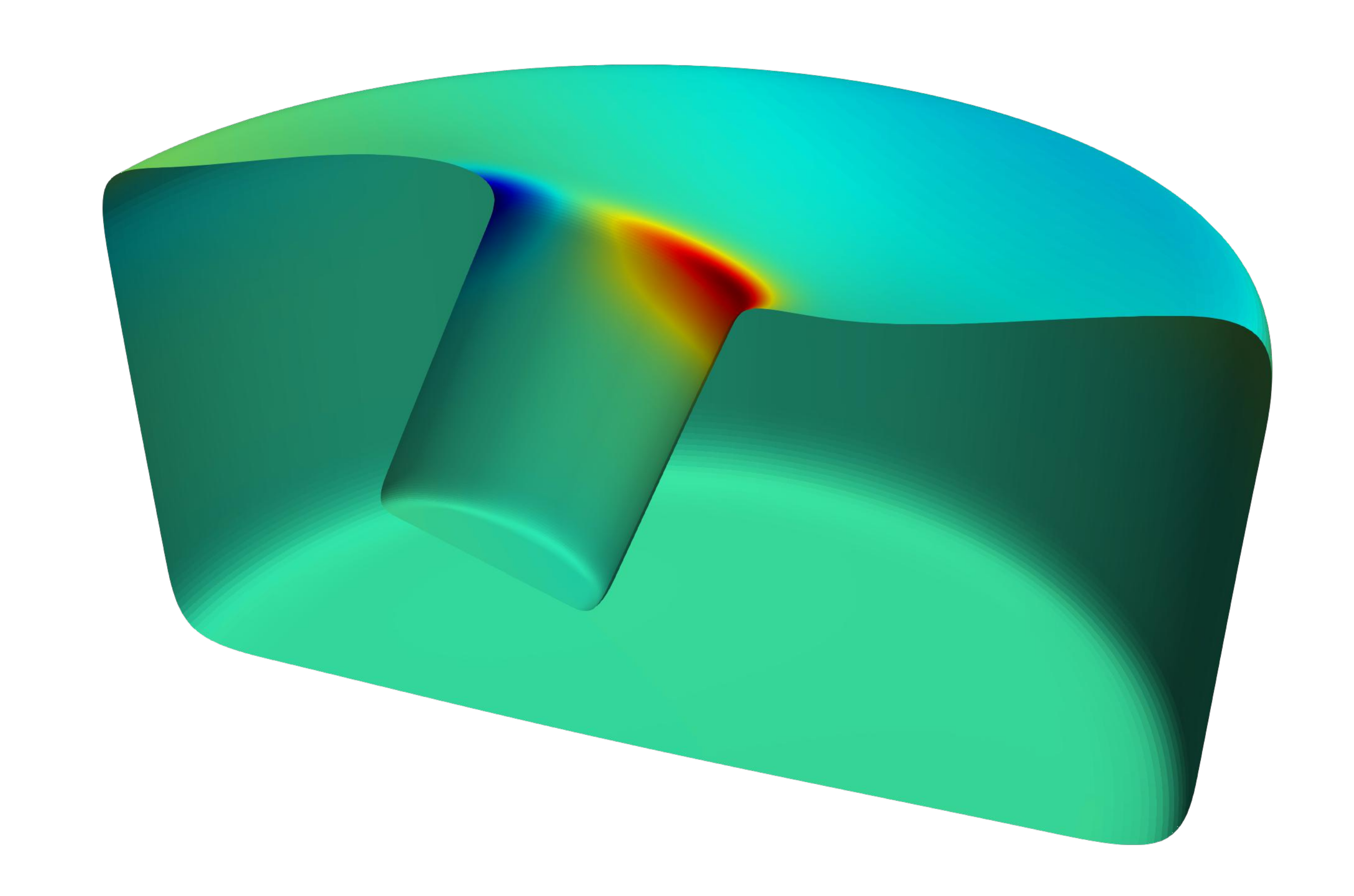}
	\caption{Mode 2: $\psi$ on the shell.}
	\label{fig:piezo_buzzer_mode_psi_2}
\end{subfigure}
\begin{subfigure}[b]{0.3\linewidth}
	\centering
 	 \includegraphics[width=\linewidth]{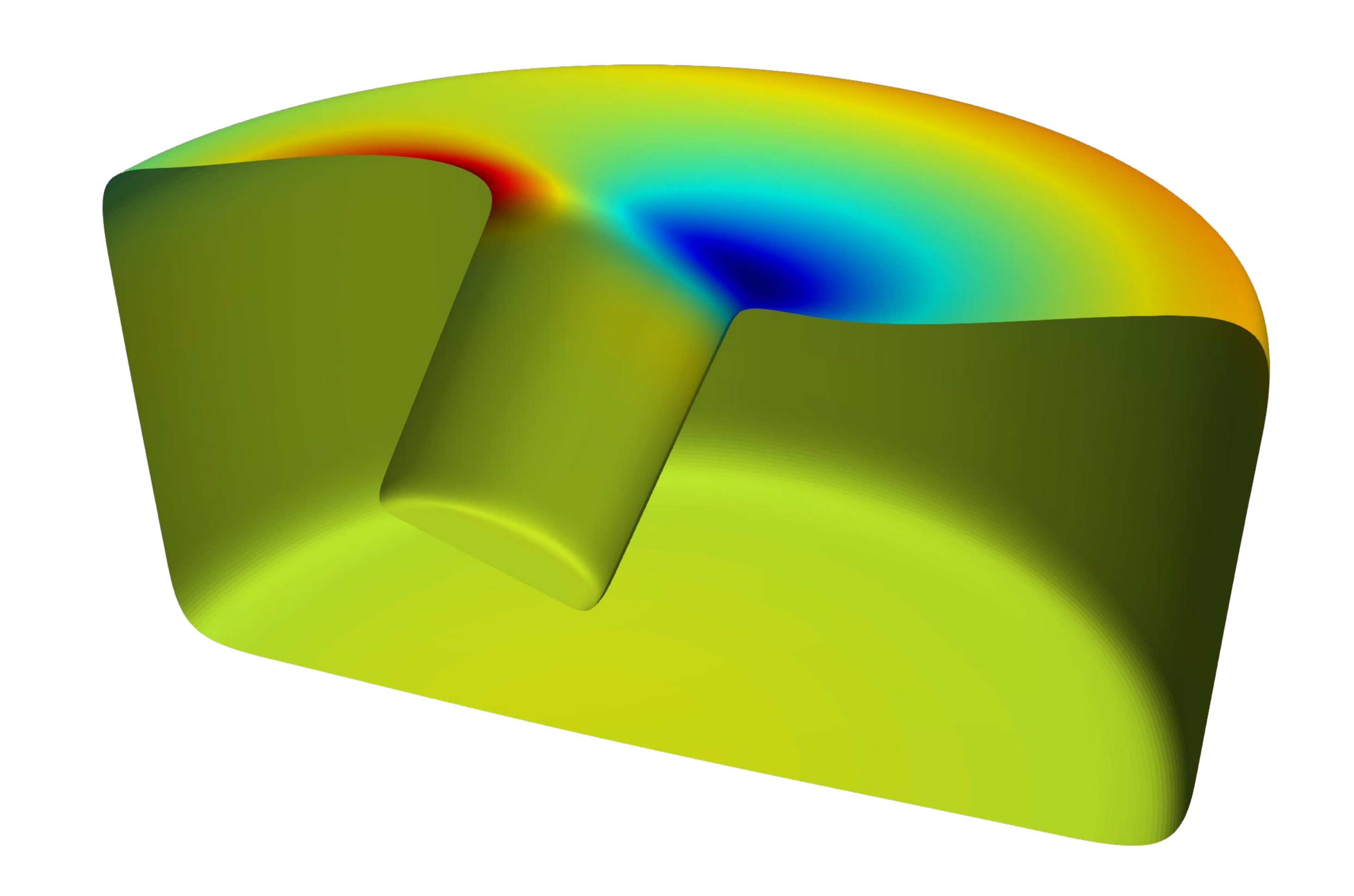}
	\caption{Mode 2: $\varphi$ on the shell.}
	\label{fig:piezo_buzzer_mode_phi_2}
\end{subfigure}

\begin{subfigure}[b]{0.3\linewidth}
	\centering
 	 \includegraphics[width=\linewidth]{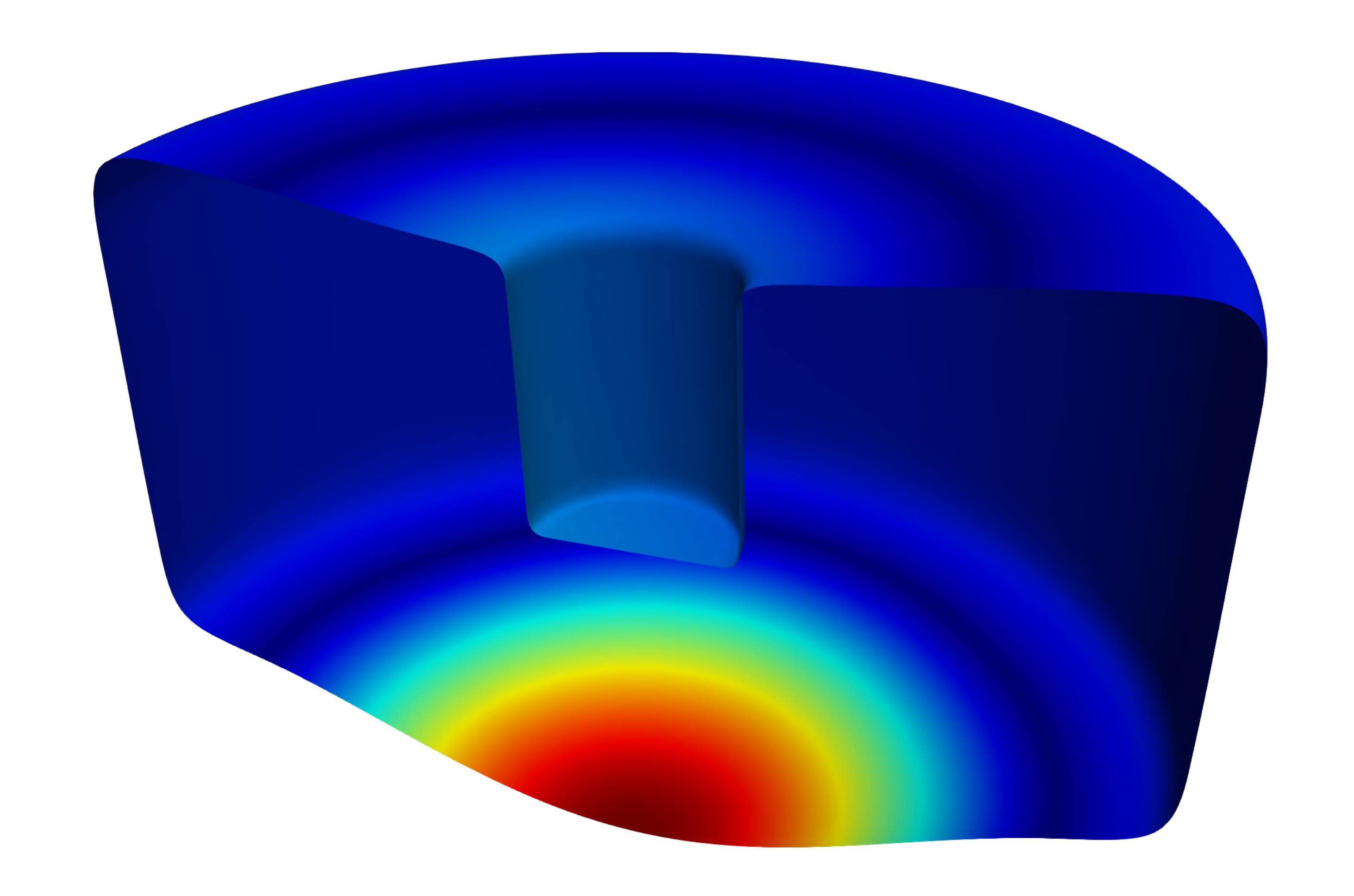}
	\caption{Mode 3: $|\mathbf u|$ on the shell.}
	\label{fig:piezo_buzzer_mode_disp_4}
\end{subfigure}
\begin{subfigure}[b]{0.3\linewidth}
	\centering
 	 \includegraphics[width=\linewidth]{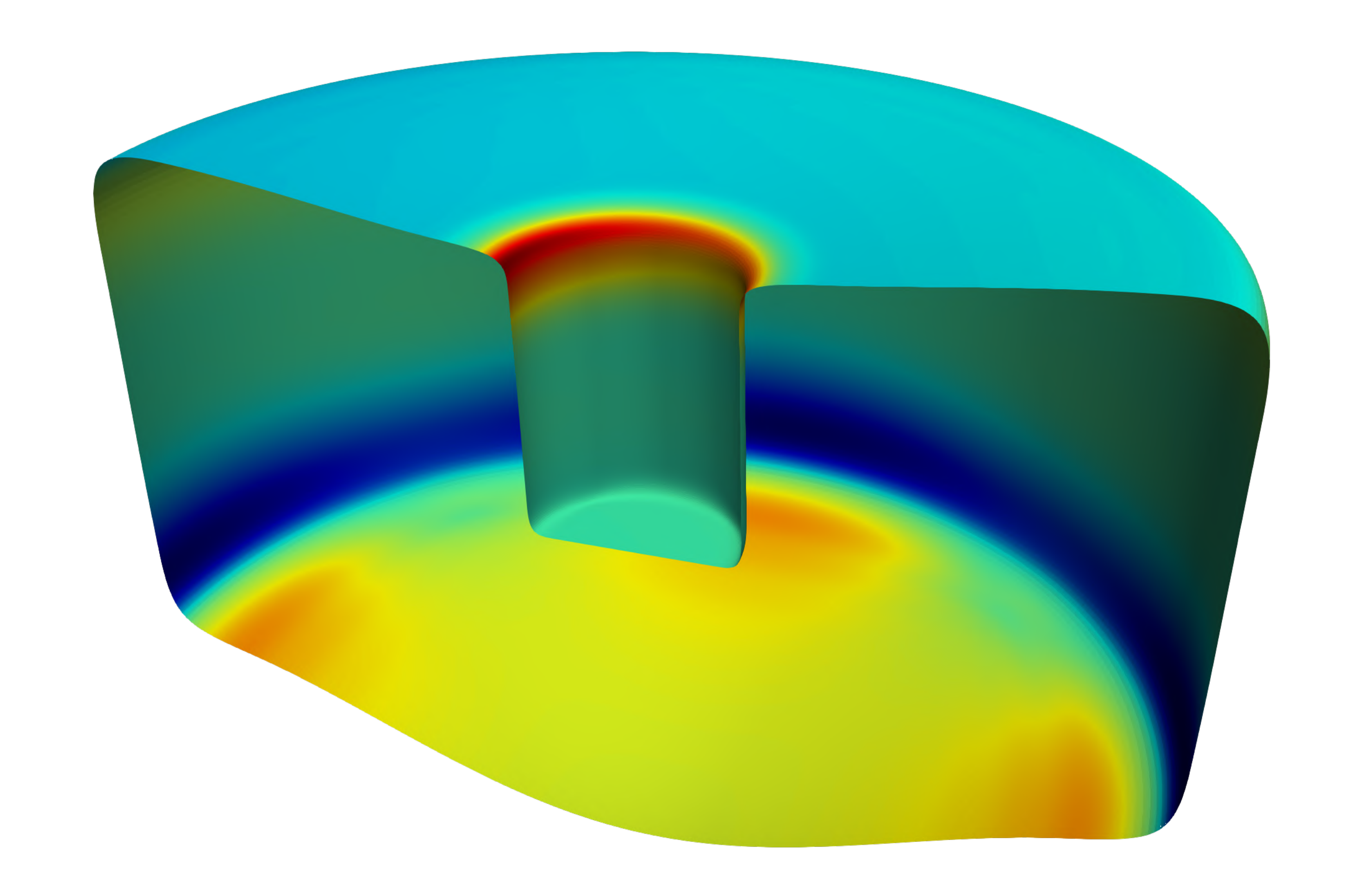}
	\caption{Mode 3: $\psi$ on the shell.}
	\label{fig:piezo_buzzer_mode_psi_4}
\end{subfigure}
\begin{subfigure}[b]{0.3\linewidth}
	\centering
 	 \includegraphics[width=\linewidth]{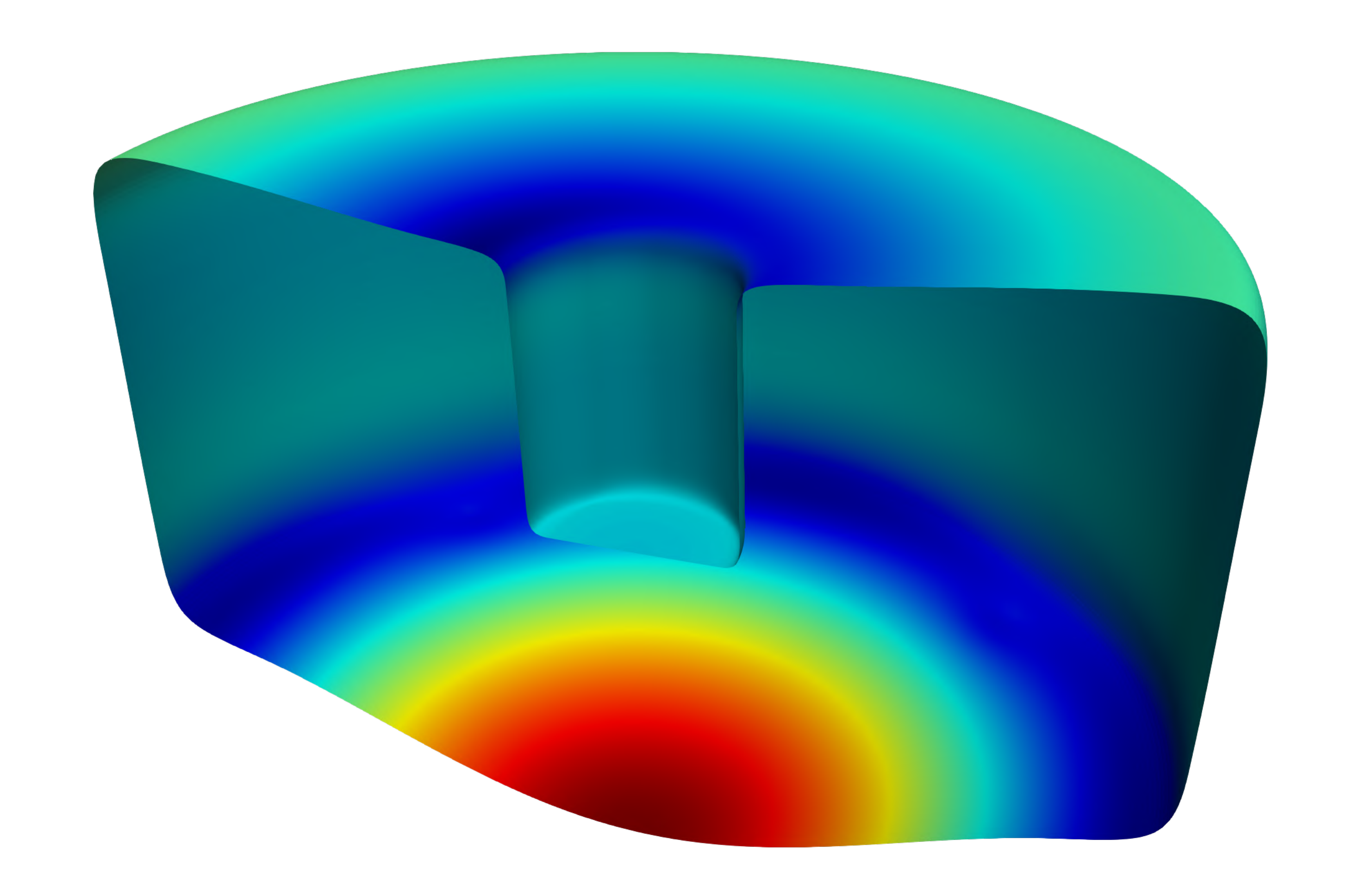}
	\caption{Mode 3: $\varphi$ on the shell.}
	\label{fig:piezo_buzzer_mode_phi_4}
\end{subfigure}

%	\caption{[TO ADD]}
%	\label{fig:pb_modes_1}
%\end{figure}
%
%\begin{figure}
%\centering

\begin{subfigure}[b]{0.3\linewidth}
	\centering
 	 \includegraphics[width=\linewidth]{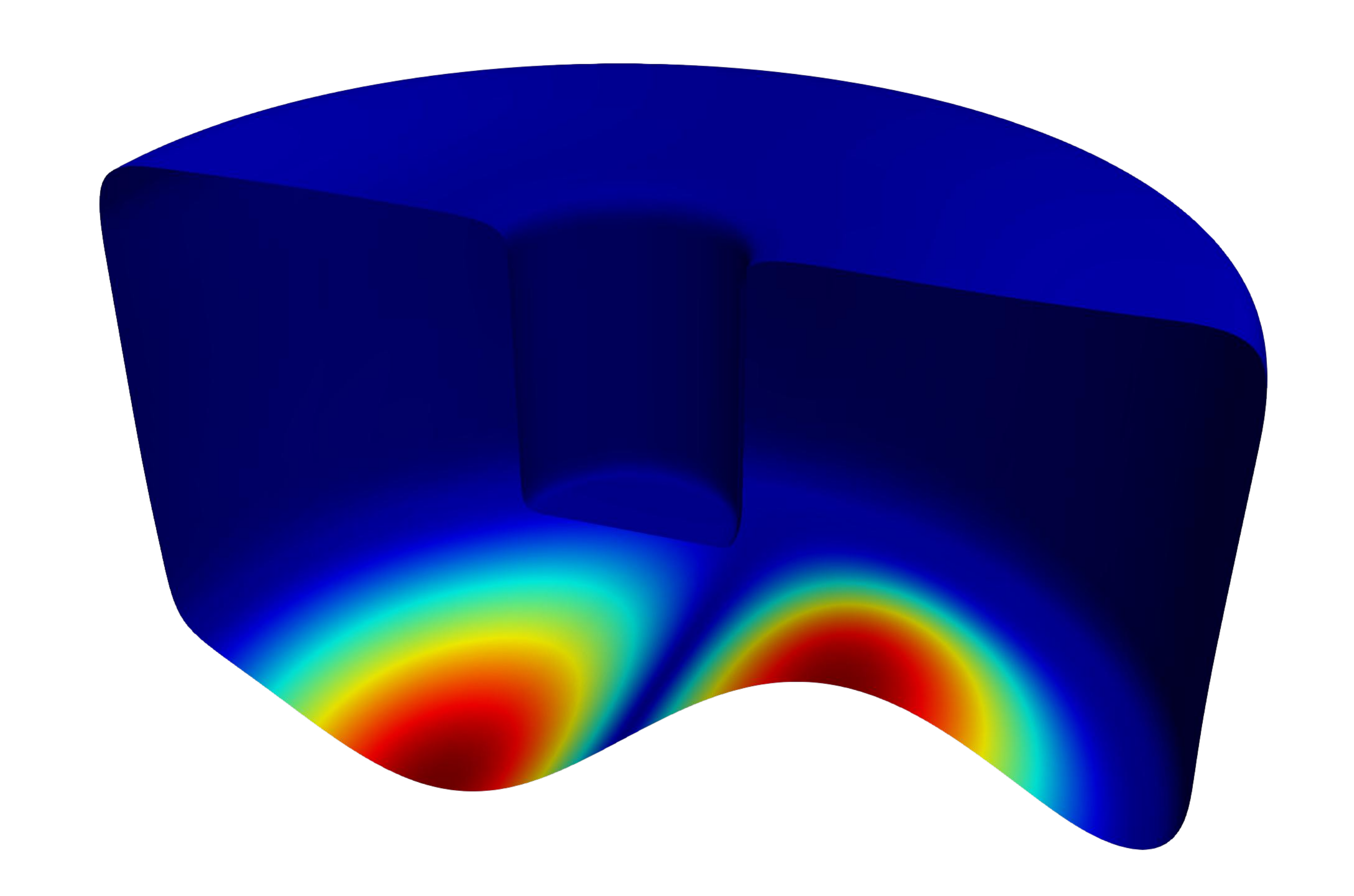}
	\caption{Mode 4: $|\mathbf u|$ on the shell.}
	\label{fig:piezo_buzzer_mode_disp_5}
\end{subfigure}
\begin{subfigure}[b]{0.3\linewidth}
	\centering
 	 \includegraphics[width=\linewidth]{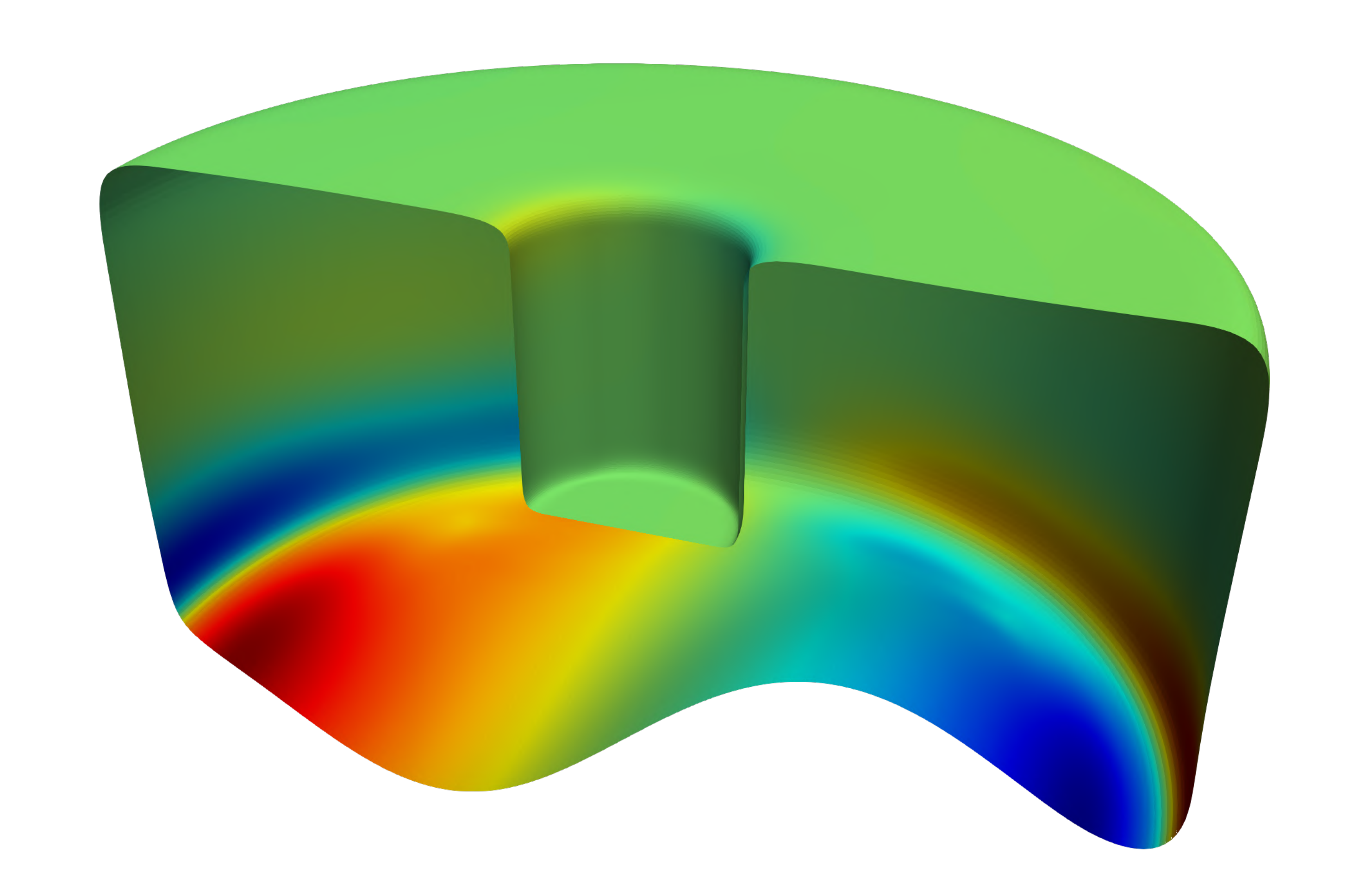}
	\caption{Mode 4: $\psi$ on the shell.}
	\label{fig:piezo_buzzer_mode_psi_5}
\end{subfigure}
\begin{subfigure}[b]{0.3\linewidth}
	\centering
 	 \includegraphics[width=\linewidth]{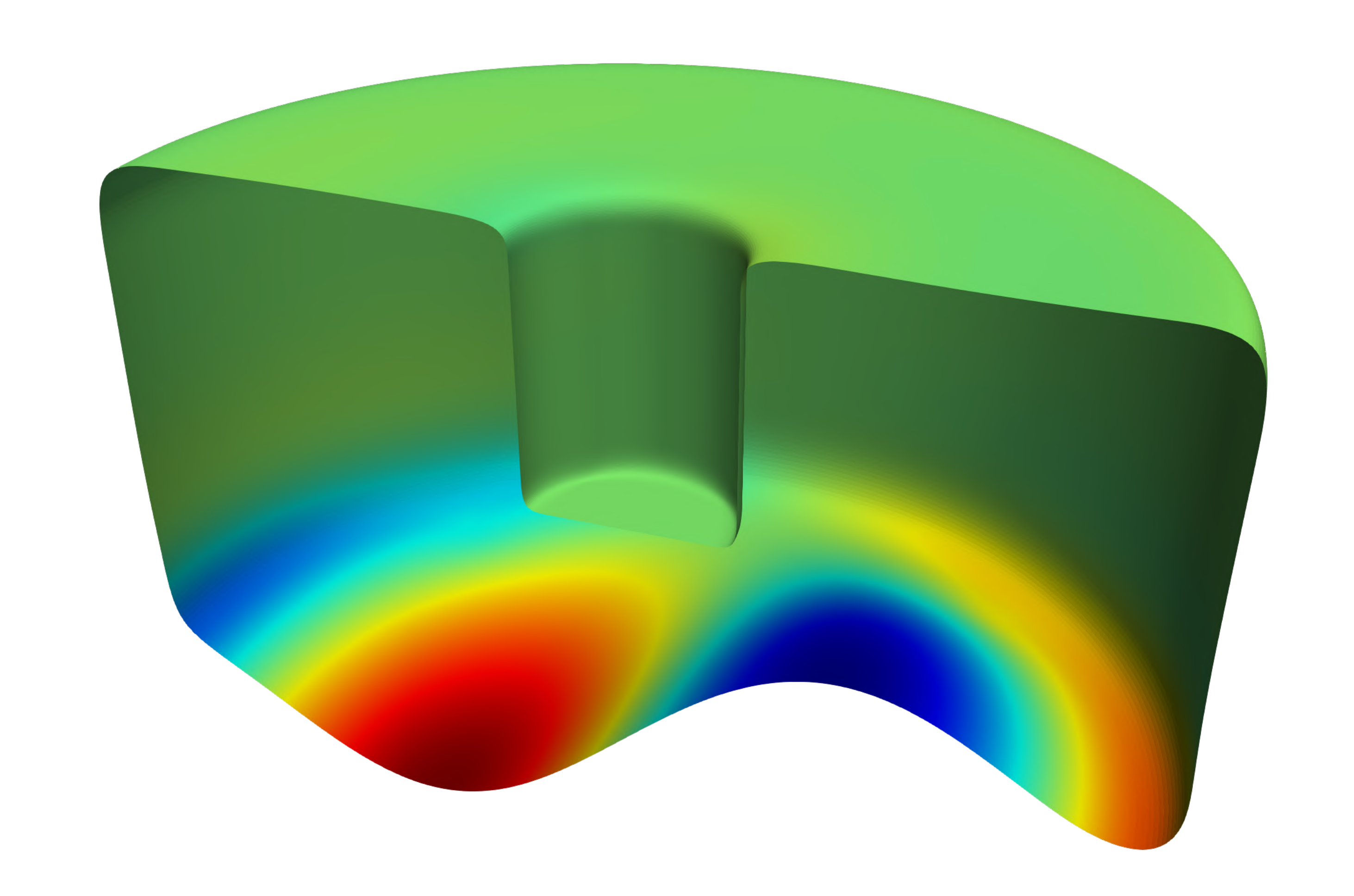}
	\caption{Mode 4: $\varphi$ on the shell.}
	\label{fig:piezo_buzzer_mode_phi_5}
\end{subfigure}

%\begin{subfigure}[b]{0.4\linewidth}
%	\centering
% 	 \includegraphics[width=\linewidth]{piezo_buzzer_mode_disp_8.pdf}
%	\caption{Mode7,8: $|\mathbf u|$ on the shell.}
%	\label{fig:piezo_buzzer_mode_disp_8}
%\end{subfigure}
%\begin{subfigure}[b]{0.4\linewidth}
%	\centering
% 	 \includegraphics[width=\linewidth]{piezo_buzzer_mode_phi_8.pdf}
%	\caption{Mode 7,8: $\varphi$ on the shell.}
%	\label{fig:piezo_buzzer_mode_phi_8}
%\end{subfigure}

\caption{First four vibration modes of the piezoelectric speaker structure. The magnitude of displacement $|\mathbf u| $ and the potential functions $\psi$ and $\varphi$ are plotted on half of the deformed mid-surface of the structure.}
\label{fig:pb_modes}
\end{figure}

Figure~\ref{fig:pb_modes} shows the first four modes of this structure. Modes $1$ and $3$ are axisymmetric. Mode $2$ corresponds to two identical eigenvalues which are the $2^\text{nd}$ and $3^\text{rd}$. Similarly, mode $4$ also relates to the $5^\text{th}$ and $6^\text{th}$ eigenvalues, which are also identical. Table~\ref{tab:pb_result} compares the eigenmode frequency of the piezoelectric shell against a pure elastic shell with approximately a $4\%$ rise in the frequencies for the first four modes.

\begin{table}[]
\centering
\caption{Eigenmode frequencies for the elastic and the piezoelectric speaker.}
\begin{tabular}{@{}CCCCC@{}}
\toprule
n_d & $Eigenvalue No.$ & $Elastic $ f $(Hz)$ & $Coupled $ f $(Hz)$ & $Difference$ $($\%$)$  \\ \midrule
1 & 1    & 769.4   & 800.8  & 4.08\\
2 & 2,3  & 987.4   & 1032.0 & 4.52\\
3 & 4    & 1001.3  & 1045.5 & 4.41\\
4 & 5,6  & 1876.8  & 1959.0  & 4.38\\ \bottomrule
\end{tabular}
\label{tab:pb_result}
\end{table}

\section{Conclusions}
An isogeometric Galerkin method for the vibration analysis of piezoelectric thin shells has been proposed. The shell formulation follows the Kirchhoff-Love hypothesis. Hamilton's variational principle has been adopted to formulate the weak form of the governing equations for the coupled problem and Catmull-Clark subdivision bases have been used for discretising the geometry and physical fields. A Galerkin method has been implemented using the finite element library deal.II. Assuming the piezoelectric shell vibrates harmonically, the problem renders an eigenvalue problem for the system matrix. The vibration of a purely elastic shell has been verified first with a spherical shell benchmark. Then the electromechanical coupling effects of piezoelectric shells with different curvature have been evaluated and compared using curved plates. In general, the natural frequencies of the piezoelectric structure are higher than those of the structure in the absence of the piezoelectric effect. This ``piezoelectric stiffening'' effect is particularly significant for certain modes. Finally, an example has been presented to demonstrate the capability of the proposed method in the design and analysis of piezoelectric shells with complex geometry. 

The main findings of the study are threefold.
First, the effect of piezoelectric coupling for thin shell structures with arbitrary geometries, as applicable to realistic applications generated from CAD, can be described using the isogeometric method. 
Second, the method describes three different types of relevant electrical conditions: no-electrodes, prescribed voltage and short-circuited. The relationship between the strain and electric potential has been made clear.
Third, the change of natural frequency of piezoelectric shells with curvature can be accurately represented. This will provide valuable guidance for the design of piezoelectric energy harvesters.

In future work, the proposed method will be extended to account for large deformation and instabilities of thin shell structures made of electroelastic polymers \cite{liu2020coupled}.

\section*{Acknowledgements}
This work was supported by the UK Engineering and Physical Sciences Research Council grant EP/R008531/1 for the Glasgow Computational Engineering Centre. We also thank for the support from the Royal Society International Exchange Scheme IES/R1/201122.\\
\indent Paul Steinmann gratefully acknowledges financial support for this work by the Deutsche Forschungsgemeinschaft under GRK2495, projects B \& C. We are particularly grateful to Andreas Hegend\"orfer for discussions on the topic of piezoelectric energy harvesters.\\
\indent Yilin Qu also acknowledges the support from the Fundamental Research Funds for the Central Universities (No. xzy022020016).

\bibliographystyle{cas-model2-names}      % basic style, author-year citations
\bibliography{cas-refs}
\end{document}